\documentclass[12pt]{article}
\usepackage{amsmath}
\usepackage{amssymb}
\usepackage[ctr1,eng]{myarticl}

\usepackage{ifthen}

\newcommand{\aNdbasis }{\boldsymbol{\mathrm{f}}}

\newcommand{\Aut }{\mathrm{Aut}}

\newcommand{\cB }{\mathcal{B}}

\newcommand{\cG }{\mathcal{G}}

\newcommand{\Dchaintwo}[4]{
\rule[-3\unitlength]{0pt}{8\unitlength}
\begin{picture}(14,5)(0,3)
\put(1,2){\ifthenelse{\equal{#1}{l}}{\circle*{2}}{\circle{2}}}
\put(2,2){\line(1,0){10}}
\put(13,2){\ifthenelse{\equal{#1}{r}}{\circle*{2}}{\circle{2}}}
\put(1,5){\makebox[0pt]{\scriptsize #2}}
\put(7,4){\makebox[0pt]{\scriptsize #3}}
\put(13,5){\makebox[0pt]{\scriptsize #4}}
\end{picture}}
\newcommand{\Dchainthree}[6]{
\rule[-3\unitlength]{0pt}{8\unitlength}
\begin{picture}(26,5)(0,3)
\put(1,2){\ifthenelse{\equal{#1}{l}}{\circle*{2}}{\circle{2}}}
\put(2,2){\line(1,0){10}}
\put(13,2){\ifthenelse{\equal{#1}{m}}{\circle*{2}}{\circle{2}}}
\put(14,2){\line(1,0){10}}
\put(25,2){\ifthenelse{\equal{#1}{r}}{\circle*{2}}{\circle{2}}}
\put(1,5){\makebox[0pt]{\scriptsize #2}}
\put(7,4){\makebox[0pt]{\scriptsize #3}}
\put(13,5){\makebox[0pt]{\scriptsize #4}}
\put(19,4){\makebox[0pt]{\scriptsize #5}}
\put(25,5){\makebox[0pt]{\scriptsize #6}}
\end{picture}}
\newcommand{\Dtriangle}[7]{
\rule[-3\unitlength]{0pt}{12\unitlength}
\begin{picture}(18,7)(0,3)
\put(4,4){\ifthenelse{\equal{#1}{l}}{\circle*{2}}{\circle{2}}}
\put(5,4){\line(1,0){8}}
\put(14,4){\ifthenelse{\equal{#1}{r}}{\circle*{2}}{\circle{2}}}
\put(4.4472,4.8944){\line(1,2){4.1056}}
\put(9,14){\ifthenelse{\equal{#1}{t}}{\circle*{2}}{\circle{2}}}
\put(13.5528,4.8944){\line(-1,2){4.1056}}
\put(2,3){\makebox[0pt][r]{\scriptsize #2}}
\put(9,16){\makebox[0pt]{\scriptsize #3}}
\put(16,3){\makebox[0pt][l]{\scriptsize #4}}
\put(6,9){\makebox[0pt][r]{\scriptsize #5}}
\put(12.5,9){\makebox[0pt][l]{\scriptsize #6}}
\put(9,1){\makebox[0pt]{\scriptsize #7}}
\end{picture}}
\newcommand{\cDtriangle}[7]{
\rule[-10\unitlength]{0pt}{\unitlength}
\begin{picture}(18,7)(0,10)
\put(4,4){\ifthenelse{\equal{#1}{l}}{\circle*{2}}{\circle{2}}}
\put(5,4){\line(1,0){8}}
\put(14,4){\ifthenelse{\equal{#1}{r}}{\circle*{2}}{\circle{2}}}
\put(4.4472,4.8944){\line(1,2){4.1056}}
\put(9,14){\ifthenelse{\equal{#1}{t}}{\circle*{2}}{\circle{2}}}
\put(13.5528,4.8944){\line(-1,2){4.1056}}
\put(2,3){\makebox[0pt][r]{\scriptsize #2}}
\put(9,17){\makebox[0pt]{\scriptsize #3}}
\put(16,3){\makebox[0pt][l]{\scriptsize #4}}
\put(6,9){\makebox[0pt][r]{\scriptsize #5}}
\put(12.5,9){\makebox[0pt][l]{\scriptsize #6}}
\put(9,1){\makebox[0pt]{\scriptsize #7}}
\end{picture}}

\newcommand{\extWBG}{W^\mathrm{ext}}
\newcommand{\gDd }{\mathcal{D}}
\newcommand{\gr }{\mathrm{gr}}

\newcommand{\longDchaintwo}[4]{
\rule[-3\unitlength]{0pt}{8\unitlength}
\begin{picture}(18,5)(0,3)
\put(1,2){\ifthenelse{\equal{#1}{l}}{\circle*{2}}{\circle{2}}}
\put(2,2){\line(1,0){14}}
\put(17,2){\ifthenelse{\equal{#1}{r}}{\circle*{2}}{\circle{2}}}
\put(1,5){\makebox[0pt]{\scriptsize #2}}
\put(9,4){\makebox[0pt]{\scriptsize #3}}
\put(17,5){\makebox[0pt]{\scriptsize #4}}
\end{picture}}

\newcommand{\Ndbasis }{\boldsymbol{\mathrm{e}}}

\newcommand{\ndN }{\mathbb{N}}

\newcommand{\ndZ }{\mathbb{Z}}
\newcommand{\op }{^\mathrm{op}}

\newcommand{\proots }{\boldsymbol{\Delta }^+}

\newcommand{\restr }[1]{\!\upharpoonright _{#1}}

\newcommand{\roots }{\boldsymbol{\Delta }}

\newcommand{\WBg }{\mathcal{G}}
\newcommand{\YD }{Yetter--Drinfel'd }

\newcommand{\Dchainfour}[8]{
\rule[-3\unitlength]{0pt}{5\unitlength}
\begin{picture}(38,5)(0,3)
\put(1,2){\ifthenelse{\equal{#1}{1}}{\circle*{2}}{\circle{2}}}
\put(2,2){\line(1,0){10}}
\put(13,2){\ifthenelse{\equal{#1}{2}}{\circle*{2}}{\circle{2}}}
\put(14,2){\line(1,0){10}}
\put(25,2){\ifthenelse{\equal{#1}{3}}{\circle*{2}}{\circle{2}}}
\put(26,2){\line(1,0){10}}
\put(37,2){\ifthenelse{\equal{#1}{4}}{\circle*{2}}{\circle{2}}}
\put(1,5){\makebox[0pt]{\scriptsize #2}}
\put(7,4){\makebox[0pt]{\scriptsize #3}}
\put(13,5){\makebox[0pt]{\scriptsize #4}}
\put(19,4){\makebox[0pt]{\scriptsize #5}}
\put(25,5){\makebox[0pt]{\scriptsize #6}}
\put(31,4){\makebox[0pt]{\scriptsize #7}}
\put(37,5){\makebox[0pt]{\scriptsize #8}}
\end{picture}}
\newcommand{\Dchainfive}[9]{
\rule[-4\unitlength]{0pt}{5\unitlength}
\begin{picture}(50,4)(0,3)
\put(1,1){\circle{2}}
\put(2,1){\line(1,0){10}}
\put(13,1){\circle{2}}
\put(14,1){\line(1,0){10}}
\put(25,1){\circle{2}}
\put(26,1){\line(1,0){10}}
\put(37,1){\circle{2}}
\put(38,1){\line(1,0){10}}
\put(49,1){\circle{2}}
\put(1,4){\makebox[0pt]{\scriptsize #1}}
\put(7,3){\makebox[0pt]{\scriptsize #2}}
\put(13,4){\makebox[0pt]{\scriptsize #3}}
\put(19,3){\makebox[0pt]{\scriptsize #4}}
\put(25,4){\makebox[0pt]{\scriptsize #5}}
\put(31,3){\makebox[0pt]{\scriptsize #6}}
\put(37,4){\makebox[0pt]{\scriptsize #7}}
\put(43,3){\makebox[0pt]{\scriptsize #8}}
\put(49,4){\makebox[0pt]{\scriptsize #9}}
\end{picture}}
\newcommand{\bigDchainfour}[8]{
\rule[-3\unitlength]{0pt}{5\unitlength}
\begin{picture}(42,5)(0,3)
\put(1,2){\ifthenelse{\equal{#1}{1}}{\circle*{2}}{\circle{2}}}
\put(2,2){\line(1,0){14}}
\put(17,2){\ifthenelse{\equal{#1}{2}}{\circle*{2}}{\circle{2}}}
\put(18,2){\line(1,0){10}}
\put(29,2){\ifthenelse{\equal{#1}{3}}{\circle*{2}}{\circle{2}}}
\put(30,2){\line(1,0){10}}
\put(41,2){\ifthenelse{\equal{#1}{4}}{\circle*{2}}{\circle{2}}}
\put(1,5){\makebox[0pt]{\scriptsize #2}}
\put(9,4){\makebox[0pt]{\scriptsize #3}}
\put(17,5){\makebox[0pt]{\scriptsize #4}}
\put(23,4){\makebox[0pt]{\scriptsize #5}}
\put(29,5){\makebox[0pt]{\scriptsize #6}}
\put(35,4){\makebox[0pt]{\scriptsize #7}}
\put(41,5){\makebox[0pt]{\scriptsize #8}}
\end{picture}}

\newcommand{\Dthreefork}[8]{
\rule[-9\unitlength]{0pt}{12\unitlength}
\begin{picture}(28,12)(0,9)
\put(2,10){\ifthenelse{\equal{#1}{l}}{\circle*{2}}{\circle{2}}}
\put(3,10){\line(1,0){10}}
\put(14,10){\ifthenelse{\equal{#1}{m}}{\circle*{2}}{\circle{2}}}
\put(15,10){\line(1,1){7}}
\put(15,10){\line(1,-1){7}}
\put(22,18){\ifthenelse{\equal{#1}{t}}{\circle*{2}}{\circle{2}}}
\put(22,2){\ifthenelse{\equal{#1}{b}}{\circle*{2}}{\circle{2}}}
\put(2,12){\makebox[0pt]{\scriptsize #2}}
\put(8,11){\makebox[0pt]{\scriptsize #3}}
\put(14,12){\makebox[0pt]{\scriptsize #4}}
\put(19,16){\makebox[0pt][r]{\scriptsize #5}}
\put(19,4){\makebox[0pt][r]{\scriptsize #6}}
\put(24,17){\makebox[0pt][l]{\scriptsize #7}}
\put(24,2){\makebox[0pt][l]{\scriptsize #8}}
\end{picture}}

\newcommand{\Drightofway}[9]{
\rule[-9\unitlength]{0pt}{12\unitlength}
\begin{picture}(28,12)(0,9)
\put(2,10){\ifthenelse{\equal{#1}{l}}{\circle*{2}}{\circle{2}}}
\put(3,10){\line(1,0){10}}
\put(14,10){\ifthenelse{\equal{#1}{m}}{\circle*{2}}{\circle{2}}}
\put(15,10){\line(1,1){7}}
\put(15,10){\line(1,-1){7}}
\put(22,18){\ifthenelse{\equal{#1}{t}}{\circle*{2}}{\circle{2}}}
\put(22,2){\ifthenelse{\equal{#1}{b}}{\circle*{2}}{\circle{2}}}
\put(22,3){\line(0,1){14}}
\put(2,12){\makebox[0pt]{\scriptsize #2}}
\put(8,11){\makebox[0pt]{\scriptsize #3}}
\put(14,12){\makebox[0pt]{\scriptsize #4}}
\put(19,16){\makebox[0pt][r]{\scriptsize #5}}
\put(19,4){\makebox[0pt][r]{\scriptsize #6}}
\put(24,18){\makebox[0pt][l]{\scriptsize #7}}
\put(23,10){\makebox[0pt][l]{\scriptsize #8}}
\put(24,1){\makebox[0pt][l]{\scriptsize #9}}
\end{picture}}

\title{Classification of arithmetic root systems}
\author{I.~Heckenberger}

\begin{document}

\setlength{\unitlength}{1mm}

\maketitle

\begin{abstract}
Arithmetic root systems are
invariants of Nichols algebras of diagonal type
with a certain finiteness property. 
They can also be considered as generalizations of ordinary root systems
with rich structure and many new examples.
On the other hand, Nichols algebras
are fundamental objects in the construction of quantized enveloping algebras,
in the noncommutative differential geometry of quantum groups,
and in the classification of pointed Hopf algebras by the lifting method
of Andruskiewitsch and Schneider.
In the present paper arithmetic root systems are classified in full generality.
As a byproduct many new finite dimensional pointed Hopf algebras are obtained.


\end{abstract}

\section{Introduction}

The theory of Nichols algebras is relatively young, but it
is affected by various research areas of mathematics and theoretical physics.
It is dominated
and motivated by Hopf algebra theory in the following way.
Let $H$ be a Hopf algebra with
coradical filtration $H_0\subset H_1\subset
\ldots $ such that $H_0$ is a Hopf subalgebra of $H$.
Let $\gr \,H$ denote the $\ndN _0$-graded Hopf algebra
$\bigoplus _iH_i/H_{i-1}$.
Then $H$ possesses a rich invariant, namely the subalgebra $\cB (V)\subset
\gr \,H$ generated by the vector space $V$ of $H_0$-coinvariants of $H_1/H_0$.
It is called a \textit{Nichols algebra} \cite{a-AndrSchn98} in commemoration
to W.~Nichols who started to study these objects
systematically \cite{a-Nichols78}. Nichols algebras
can be described in many different ways
\cite{a-Schauen96}, \cite{a-AndrGr99}, \cite{a-Rosso98}.
The importance of such algebras was detected and pointed out in many papers by
Andruskiewitsch and Schneider, see for example \cite{a-AndrSchn98} and
\cite{inp-Andr02}. Their structure was enlightened,
mainly in the case when $H_0$ is the group algebra of a finite group, among
others in \cite{inp-Grana99},
\cite{a-AndrSchn00}, \cite{inp-MilSchn00}, \cite{a-AndrGr03}.
Nichols algebras were used by Andruskiewitsch and Schneider
to start a very promising
program \cite{a-AndrSchn98}
to classify pointed Hopf algebras with certain finiteness properties.
This so called lifting method was already succesfully performed for finite
dimensional pointed Hopf algebras where the nilpotency order of the
elements of $H_1/H_0$ is bigger than 7 \cite{a-AndrSchn05p}.

Nichols algebras appear in a natural way also in the
construction of quantized Kac--Moody algebras
\cite[Sect.~8.2.1]{b-KS} \cite[Sect.\,3.2.9]{b-Joseph}
and their $\ndZ _2$-graded variants \cite{a-KhorTol91}.
Using a particular Nichols algebra
Yamane \cite{a-Yamane06} described
a $\ndZ /3\ndZ $-graded quantum group which has a representation theory
fitting into the general picture.
Nichols algebras
are also natural objects in the theory of covariant differential
calculus on quantum groups initiated by Woronowicz \cite{a-Woro2}.
Further, Bazlov \cite{a-Bazlov06} proved that the cohomology ring of a
flag variety can be considered as a subalgebra of a particular
Nichols algebra of nonabelian group type.
In contrast to these various interesting aspects of the subject
one still does not know too much about the structure of
Nichols algebras in general.

Color Lie algebras \cite{b-Bahturin92} are generalizations of Lie algebras.
For the study of their structure methods are developed which are useful
also to analyze Nichols algebras. For example,
Kharchenko \cite{a-Khar99} proved
that any Hopf algebra generated by skew-primitive and group-like elements
has a restricted Poincar\'e--Birkhoff-Witt basis.
Here ``restricted''
means that the possible powers of the root vectors can be bounded by an integer
number. However
in contrast to color Lie algebras this bound may be different from 2. Therefore
significant new classes of examples can be expected. However note that
under some hypotheses all examples are deformations of the upper triangular
part of a semisimple Lie algebra \cite{a-Rosso98} \cite{a-AndrSchn00}.
Kharchenko's results
apply in particular to Nichols algebras $\cB (V)$ of diagonal type,
that is when $V$ is a direct sum of 1-dimensional \YD modules over $H_0$.
Motivated by the close relation to Lie theory, to any such Nichols algebra
a Weyl groupoid and an arithmetic root system were
associated \cite{a-Heck04c}. These constructions were used
in \cite{a-Heck04d} and \cite{a-Heck05a}
to determine Nichols algebras of rank 2 and 3
with a finite set of PBW generators
without dealing with the complicated defining relations of $\cB (V)$.
In this paper the classification is performed for Nichols algebras of
diagonal type and of
arbitrary (finite)
rank following the ideas in \cite{a-Heck04d} and \cite{a-Heck05a}.
The main results are collected in Theorems~\ref{t-classrank4}
and \ref{t-classrank>4}.
Their proof was enabled by the additional development of an efficient
technique to recognize that a given triple is an arithmetic root system.
Together with the classification results for rank 2 and rank 3 Nichols algebras
of diagonal type the first half of Question~5.9 of
Andruskiewitsch~\cite{inp-Andr02} is answered.

With this classification result a large amount of algebras are detected
which did not appear previously in the literature.
In particular, infinite series of finite dimensional $\ndZ /3\ndZ $-graded
algebras are listed. However also many open problems remained unsolved, for
example to clarify the relationship
of these algebras to super Lie algebras and related structures, and
to describe them with help of generators and relations
\cite[Question 5.9]{inp-Andr02}.

In the present paper the notations and conventions in \cite{a-Heck04e}
and \cite{a-Heck05a} are followed and several results from these papers
will be used.

\section{On the finiteness of the Weyl groupoid}
\label{sec-fWg}

Let $k$ be a field of characteristic zero, $d\in \ndN $,
$\chi $ a bicharacter on $\ndZ ^d$ with values in $k^*=k\setminus \{0\}$,
and $E=\{\Ndbasis _1,\ldots ,\Ndbasis _d\}$ a basis of $\ndZ ^d$. 
Set $q_{ij}:=\chi (\Ndbasis _i,\Ndbasis _j)$ for all $i,j\in \{1,\ldots ,d\}$.

In \cite {a-Heck04c} the Weyl groupoid $W_{\chi ,E}$ associated to $\chi $ and
$E$ was defined, see also \cite{a-Heck04e} for the related definition
of $\extWBG _{\chi ,E}$.
If $W_{\chi ,E}$ (or equivalently $\extWBG _{\chi ,E}$)
is full and finite then one obtains an arithmetic
root system $(\roots ,\chi ,E)$ \cite{a-Heck04e},
where $\roots \subset \ndZ ^d$ is a certain
finite subset such that $\roots =-\roots $.
In order to classify arithmetic root systems one has to solve two problems. First,
one has to detect Weyl groupoids which are not full or not finite. This can
be done effectively using subsystems of arithmetic root systems, see
\cite{a-Heck05a}. Second, one has to be able to check the finiteness of full and
finite Weyl groupoids. To do the latter in \cite{a-Heck05a} the group
$\WBg _{\chi ,E}$ was introduced. However in general it is not easy to
determine the structure of this group. As an alternative approach
Proposition~\ref{s-Wfinite} will be proved.
This allows to conclude the finiteness
of $W_{\chi ,E}$ from the finiteness of standard subgroupoids of
$\extWBG _{\chi ,E}$ and from the existence
of one element with a special property.

In the rest of this section, if not stated otherwise,
let $W_{\chi ,E}$ be an arbitrary Weyl groupoid.
Set $\roots =\bigcup \{F\,|\,
(\id ,F)\in W_{\chi ,E}\}$. Note that $\roots $ is finite if and only if
$W_{\chi ,E}$ is finite.
For $(\id ,F)\in W_{\chi ,E}$ let $\roots _F^+$
denote the set $\roots _F^+:=\roots \cap \ndN _0F$.
The proof of \cite[Proposition\,1]{a-Heck04e} shows the following.

\begin{satz}\label{s-posroots}
Let $F$ be a basis of $\ndZ ^d$ such that $(\id ,F)\in W_{\chi ,E}$.
Then $\roots =\roots _F^+\cup -\roots _F^+$.
\end{satz}

\begin{satz}\label{s-s(roots+)}
Let $F$ be a basis of $\ndZ ^d$ such that $(\id ,F)\in W_{\chi ,E}$
and let $\aNdbasis _0\in F$. Then $\roots ^+_{s_{\aNdbasis _0,F}(F)}=
\roots ^+_{F}\cup \{-\aNdbasis _0\}\setminus \{\aNdbasis _0\}$. In particular,
$\roots _F^+\cap -\roots _E^+$ is finite.
\end{satz}

\begin{bew}
By definition of $s_{\aNdbasis _0,F}$
one has $s_{\aNdbasis _0,F}(\aNdbasis )-\aNdbasis \in \ndZ \aNdbasis _0$
for all $\aNdbasis \in F$. Hence if $f=\sum _{\aNdbasis \in F}a_{\aNdbasis }
\aNdbasis $ for some $f\in \ndZ ^d$, $a_{\aNdbasis }\in \ndZ $, then
$f-\sum _{\aNdbasis \in F}a_{\aNdbasis }s_{\aNdbasis _0,F}(\aNdbasis )
\in \ndZ \aNdbasis _0$. By Proposition~\ref{s-posroots} one has
$\roots ^+_{s_{\aNdbasis _0,F}(F)}\subset \ndN _0F\cup -\ndN _0F$.
Thus if $f\in \roots ^+_{s_{\aNdbasis _0,F}(F)}$ and
$f\notin \ndZ \aNdbasis _0$ then $f\in \ndN _0F\cap \roots =\roots ^+_F$.
Since $F$ is a basis of $\ndZ ^d$,
relation $f\in \ndZ \aNdbasis _0$ implies that $f=m\aNdbasis _0$ with $m^2=1$.
Therefore equation $s_{\aNdbasis _0,F}(\aNdbasis _0)=-\aNdbasis _0$
gives that
$\roots ^+_{s_{\aNdbasis _0,F}(F)}\subset 
\roots ^+_{F}\cup \{-\aNdbasis _0\}\setminus \{\aNdbasis _0\}$.
Multiplication with $-1$ yields
$-\roots ^+_{s_{\aNdbasis _0,F}(F)}\subset 
-\roots ^+_{F}\cup \{\aNdbasis _0\}\setminus \{-\aNdbasis _0\}$,
and hence Proposition~1 gives equality in the above relation.
The second part of the claim of the proposition follows from the
fact that any element of $W_{\chi ,E}$ can be written as a finite product
of elements of the form $(s_{\aNdbasis ,F},F)$.
\end{bew}

Assume that $(\id ,F)\in \extWBG _{\chi ,E}$ and let $F'$ be a subset
of $F$.
For any $\aNdbasis \in \ndZ ^d$ let $[\aNdbasis ]_{F'}$ denote the
equivalence class of $\aNdbasis $ in $\ndZ ^d/\ndZ F'$.
For a subset $E'$ of $E$ set
\begin{align}
\extWBG _{\chi ,E'\subset E}&=
\{(T,F)\in \extWBG _{\chi ,E}\,\big|\,
\text{there exists $(T',E)\in \extWBG _{\chi ,E}$ such that} \notag \\
\label{eq-extWBG'}
&T'(E)=F\text{ and }[\Ndbasis ]_{E'}=[T'(\Ndbasis )]_{E'}=[TT'(\Ndbasis )]_{E'}
\text{ for all }\Ndbasis \in E\}
\end{align}
and $W_{\chi ,E'\subset E}=
\extWBG _{\chi ,E'\subset E}\cap W_{\chi ,E}$.
Obviously, $\extWBG _{\chi ,E'\subset E}$ is a subgroupoid of
$\extWBG _{\chi ,E}$ and $W_{\chi ,E'\subset E}$ is a subgroupoid of
$W_{\chi ,E}$.

\begin{lemma}\label{l-T(F)=F}
Assume that $(T,F)\in \extWBG _{\chi ,E}$ and $F'\subset F$.
If $T(F)\cap \ndZ F'=F'$ and $[T(\aNdbasis )]_{F'}=[\aNdbasis ]_{F'}$
for all $\aNdbasis \in F\setminus F'$, then $T(F)=F$.
Moreover, if additionally
$T(\aNdbasis )=\aNdbasis $ for all $\aNdbasis \in F'$ then
$T=\id $.
\end{lemma}

\begin{bew}
By assumption, $T(F\setminus F')\cap \ndZ F'=\emptyset $ and
$T(F)\cap \ndZ F'=F'$, and hence $T(F')=F'$. In particular, since $T\in
\Aut _\ndZ (\ndZ ^d)$, $T$ induces an automorphism of $\ndZ ^d/\ndZ F'$.
Since $(\id ,F),(\id ,T(F))\in \extWBG _{\chi ,E}$, one has
$T(F)\subset \ndN _0F\cup -\ndN _0F$ and
$F\subset \ndN _0T(F)\cup -\ndN _0T(F)$. Thus equation $[T(\aNdbasis )]_{F'}
=[\aNdbasis ]_{F'}$ implies that for all $\aNdbasis \in F\setminus F'$
relations $T(\aNdbasis )-\aNdbasis \in \ndN _0F'$ and
$\aNdbasis -T(\aNdbasis )\in \ndN _0T(F')=\ndN _0F'$ hold. Therefore
$T(\aNdbasis )=\aNdbasis $ for all $\aNdbasis \in F\setminus F'$.
\end{bew}

\begin{satz}\label{s-subgroupoid}
Suppose that $\extWBG _{\chi ,E}$ is full.
Let $E'\not=\emptyset $ be a subset
of $E$ and $d'$ the number of elements of $E'$.
Let $\chi '$ denote the bicharacter
on $\ndZ ^{d'}=\ndZ E'\subset \ndZ ^d$
such that $\chi '(\Ndbasis ,\Ndbasis ')=\chi (\Ndbasis ,
\Ndbasis ')$ for all $\Ndbasis ,\Ndbasis '\in E'$.
Then the map
$\Phi ^\mathrm{ext}: \extWBG _{\chi ,E'\subset E}\to \extWBG _{\chi ',E'}$
defined by
$(T,F)\mapsto (T\restr{\ndZ E'},F\cap \ndZ E')$ is an isomorphism.
\end{satz}

\begin{bew}
Suppose that $(T,F)\in \extWBG _{\chi ,E'\subset E}$.
By definition of $\extWBG _{\chi ,E'\subset E}$
the set $F\cap \ndZ E'$ is a basis of $\ndZ E'$. In the same way one
obtains that
$T(F\cap \ndZ E')=T(F)\cap \ndZ E'$ is a basis of $\ndZ E'$.
Thus $T\restr{\ndZ E'}\in \Aut _\ndZ (\ndZ E')$, and hence the map
$\widetilde{\Phi ^\mathrm{ext}}:\extWBG _{\chi ,E'\subset E}
\to \widetilde{W_{d'}}$, where
\begin{align*}
\widetilde{W_{d'}}=\{(T',F')\,|\,
T'\in \Aut _\ndZ (\ndZ E'),F'\text{ is a basis of }\ndZ E'\},
\end{align*}
given by $(T,F)\mapsto (T\restr{\ndZ E'},F\cap \ndZ E')$, is a well-defined
map of groupoids.

First it will be proved that the map $\widetilde{\Phi ^\mathrm{ext}}$
is injective. If
$(T_1,F_1),(T_2,F_2)\in \extWBG _{\chi ,E'\subset E}$ then there exist
$(T,F_1),\in \extWBG _{\chi ,E}$ such that $T(F_1)=F_2$ and
$[\aNdbasis ]_{E'}=[T(\aNdbasis )]_{E'}$ for all $\aNdbasis \in F_1$.
Assume now that
$\widetilde{\Phi ^\mathrm{ext}}((T_1,F_1))=
\widetilde{\Phi ^\mathrm{ext}}((T_2,F_2))$. Then
$F_1\cap \ndZ E'=F_2\cap \ndZ E'$ and hence Lemma~\ref{l-T(F)=F} implies that
$F_1=F_2$. Moreover, $(T_2^{-1}T_1,F_1)$ with $F'=F_1\cap \ndZ E'$
satisfies the assumptions of the second part of Lemma~\ref{l-T(F)=F}
and therefore one has $T_1=T_2$.

It remains to show that the image of
$\widetilde{\Phi ^\mathrm{ext}}$ is equal to $\extWBG _{\chi ',E'}$.
First of all $\extWBG _{\chi ',E'}$ is contained in
$\widetilde{\Phi ^\mathrm{ext}}
(\extWBG _{\chi ,E'\subset E})$. Indeed, if $(T',E)\in \extWBG _{\chi ,E}$,
$[T'(\Ndbasis )]_{E'}=[\Ndbasis ]_{E'}$ for all $\Ndbasis \in E$, and
$\Ndbasis '\in E'$, then $(s_{T'(\Ndbasis '),T'(E)},T'(E))\in
\extWBG _{\chi ,E'\subset E}$ since $\extWBG _{\chi ,E}$ is full. Therefore
$\widetilde{\Phi ^\mathrm{ext}}((s_{T'(\Ndbasis '),T'(E)},T'(E)))=
(s_{T'(\Ndbasis '),T'(E')},T'(E'))$.

Take now $(T,F)\in \extWBG _{\chi ,E'\subset E}$ and let $(T',E)$ be in
$\extWBG _{\chi ,E}$ satisfying the conditions in
Equation~(\ref{eq-extWBG'}). Then equation $(T,F)=(TT',E)\circ (T',E)^{-1}$
holds, and
$(TT',E),(T',E)\in \extWBG _{\chi ,E'\subset E}$.
Therefore it is sufficient to show that 
\begin{center}(*)\qquad 
$\widetilde{\Phi ^\mathrm{ext}}((T,E))\in \extWBG _{\chi ',E'}$
whenever $(T,E)\in \extWBG _{\chi ,E'\subset E}$.
\end{center}
By Proposition~\ref{s-s(roots+)} the set $\roots ^+_{T(E)}\cap -\roots ^+_E$
is finite. In particular, if $T(E')\cap -\roots ^+_E\not=\emptyset $
then (*) holds for $(T,E)$ if and only if it holds for
$(s_{\aNdbasis ,T(E)}T,E)$, where $\aNdbasis \in T(E')$. Hence without loss
of generality one can assume that $T(E')\cap -\roots ^+_E=\emptyset $,
that is $T(E')\subset \roots ^+_E\cap \ndN _0E'$. Moreover,
Equation~(\ref{eq-extWBG'}) gives that $T(\Ndbasis )-\Ndbasis \in \ndZ E'$
for all $\Ndbasis \in E$, and hence $T(E)\subset \roots ^+_E$. This and
Proposition~\ref{s-posroots} yield that
$E\subset \ndN _0T(E)$ and hence $E=T(E)$.
\end{bew}

\begin{satz}\label{s-Wfinite}
Assume that $\extWBG _{\chi ,E}$ is full and that
there exist $(T,E)\in \extWBG _{\chi ,E}$, $\Ndbasis \in E$ and $\aNdbasis \in
\proots $ such that $T(E)=E\setminus \{\Ndbasis \}\cup \{-\aNdbasis \}$.
Set $E':=E\setminus \{\Ndbasis \}$. Let $\chi '$ denote the bicharacter
on $\ndZ ^{d-1}=\ndZ E'\subset \ndZ ^d$
such that $\chi '(\Ndbasis ',\Ndbasis '')=\chi (\Ndbasis ',
\Ndbasis '')$ for all $\Ndbasis ',\Ndbasis ''\in E'$,
and assume that $\extWBG _{\chi ',E'}$ is finite.
Then $\extWBG _{\chi ,E}$ is finite.
\end{satz}

\begin{bew}
It suffices to show that the set $\roots $ is finite. By assumption one has
$(\id ,T(E))\in W_{\chi ,E}$
and hence Proposition~\ref{s-s(roots+)} gives that
$\roots ^+_{T(E)}\cap -\roots ^+_E$ is finite.
Moreover, Proposition~\ref{s-posroots} implies that $\roots =\roots ^+_{T(E)}
\cup -\roots ^+_{T(E)}=\roots ^+_E\cup -\roots ^+_E$ and hence it remains to
show that the set
\begin{align*}
 \roots ^+_{T(E)}\cap \roots ^+_E=\roots \cap (-\ndN _0\aNdbasis +\ndN _0E')
\cap (\ndN _0\Ndbasis +\ndN _0E')=\roots \cap \ndN _0E'
\end{align*}
is finite.
By assumption $\extWBG _{\chi ',E'}$ is finite. Hence
\cite[Proposition\,2]{a-Heck04e} implies that there exists $(\id ,E'')\in
\extWBG _{\chi ',E'}$ such that $E''\subset -\ndN _0E'$. Thus by
Proposition~\ref{s-subgroupoid} there exists $(T',E)\in \extWBG _{\chi ,E}$
such that $T'(E')=E''$ and
$T'(\Ndbasis )\in (\Ndbasis +\ndZ E')\cap \roots ^+_E$.
Since also $(\id ,T'(E))\in \extWBG _{\chi ,E}$,
Proposition~\ref{s-s(roots+)} yields that $\roots ^+_{T'(E)}\cap -\roots ^+_E$
is a finite set. Thus the set $\roots \cap \ndN _0E'$
is finite because of the relations
$\roots \cap -\ndN _0E'=\roots \cap \ndN _0E''\cap -\ndN _0E'\subset
\roots \cap \roots ^+_{T'(E)}\cap -\roots ^+_E$.
\end{bew}

\section{Connected arithmetic root systems of rank two and three}

Recall the notation in Section~\ref{sec-fWg}.
In \cite{a-Heck04e} and \cite{a-Heck05a} arithmetic root systems of rank 2
and rank 3 were classified. For the considerations in rank 4 and higher
the following facts will be needed.

\begin{lemma}\label{l-q12=-1}
Let $(\roots ,\chi ,E)$ be a connected rank 2 arithmetic root system.
If $q_{11}q_{12}q_{21}q_{22}=-1$ then one of the following two systems of
equations holds:
\begin{align*}
&q_{11}+1=q_{12}q_{21}q_{22}-1=0,&
q_{22}+1=q_{11}q_{12}q_{21}-1=0.
\end{align*}
\end{lemma}

\begin{bew}
Since $(\roots ,\chi ,E)$ is connected, one has $q_{12}q_{21}\not=1$. Thus
the claim follows from \cite[Proposition\,9(i)]{a-Heck05a}.
\end{bew}

According to Table~2 in \cite{a-Heck05a} one obtains the following lemma.

\begin{lemma}\label{l-specrank3}
Let $(\roots ,\chi ,E)$ be a connected rank 3 arithmetic root system.
Then the following assertions hold.

(i) If $q_{13}q_{31}=1$ and $q_{11},q_{22},q_{33}\not=-1$ then
either $(\chi ,E)$ is of Cartan type or relations
$q_{ii}\in R_3$, $q_{22},q_{jj}\in R_6\cup R_9$, 
$q_{jj}q_{j2}q_{2j}=q_{22}q_{2i}q_{i2}=1$,
$(q_{j2}q_{2j}q_{22}-1)(q_{j2}q_{2j}q_{22}^2-1)=0$ hold
for an $i\in \{1,3\}$ and $j=4-i$.

(ii) If $q_{ij}q_{ji}\not=1$ for all $i\not=j$ then
$\prod _{i<j}q_{ij}q_{ji}=1$ and $\prod _{i=1}^3(q_{ii}+1)=0$.
Moreover, if also relations $q_{11}=-1$ and $q_{22},q_{33}\not=-1$
hold then $q_{12}^2q_{21}^2=q_{13}^2q_{31}^2\in R_3$ and $q_{12}q_{21}q_{22}
=q_{13}q_{31}q_{33}=1$.

(iii) If $q_{13}q_{31}=1$ and $q_{22}=-1$, $q_{11}q_{12}q_{21}=1$,
$q_{11},q_{33}\not=-1$, then either $q_{23}q_{32}q_{33}=1$ or
$q_{23}q_{32}q_{33}^2=1$, $(q_{11}q_{33}^2-1)(q_{11}q_{33}^3+1)=0$,
or $q_{33}=-q_{11}\in R_3$, $q_{23}q_{32}\in \{-1,-q_{33}\}$.

(iv) If $q_{13}q_{31}=1$, $q_{33}=-1$, $q_{11},q_{22}\not=-1$, and
$q_{11}q_{12}q_{21}= q_{12}q_{21}q_{22}=1$ then either equation
$q_{22}q_{23}q_{32}=1$ or relations
$q_{22}^2q_{23}q_{32}=1$, $q_{11}\in R_3\cup R_4\cup R_6$ hold.

(v) If $q_{13}q_{31}=1$, $q_{33}=-1$, and
$q_{12}q_{21}q_{22}=q_{22}q_{23}q_{32}=1$ then either
$(q_{11}-q_{22})(q_{11}^2-q_{22})=0$ or $q_{11}=-1$ or
$q_{11}\in R_3$, $q_{11}q_{22}=-1$.

(vi) If one has $q_{13}q_{31}=1$, $q_{11}=q_{22}=-1$, and
$q_{33}\not=-1$,
then either relations
$q_{12}q_{21}=-1$, $q_{33}\in R_3$, $q_{23}^2q_{32}^2q_{33}=1$,
or equations $q_{23}q_{32}q_{33}=1$,
$(q_{12}q_{21}+1)(q_{12}q_{21}q_{23}q_{32}+1)
(q_{12}^2q_{21}^2q_{23}q_{32}-1)(q_{12}^3q_{21}^3q_{23}q_{32}-1)=0$,
or equations
$q_{12}q_{21}q_{23}q_{32}=1$,
$(q_{23}q_{32}q_{33}-1)(q_{23}q_{32}q_{33}^2-1)
(q_{23}q_{32}+q_{33})=0$ hold.

(vii) If $q_{13}q_{31}=1$ and $-q_{22},q_{12}q_{21}q_{22},
q_{22}q_{23}q_{32}\not=1$, then $q_{12}q_{21}q_{22}q_{23}q_{32}=-1$
and one has $q_{22}\in R_3\cup R_6$, $q_{ii}=-1$, $q_{22}^2q_{2i}q_{i2}=1$,
$q_{2j}q_{j2}=-q_{22}$, $q_{jj}\in \{-1,-q_{22}^{-1}\}$
for some $i\in \{1,3\}$ and $j=4-i$.

(viii) If $q_{13}q_{31}=1$, $q_{11}=q_{33}=-1$,
$q_{12}q_{21}q_{22}=1$,
and $q_{22}\not=-1$, then either $q_{22}q_{23}q_{32}=1$ or
$q_{23}q_{32}=-1$, $q_{22}\in R_3\cup R_4\cup R_6$,
or $q_{23}^2q_{32}^2=q_{22}^2\in R_3$.

(ix) If $q_{ij}q_{ji}\not=1$ for all $i\not=j$ then for all $i$
there exists $j\not=i$ such that
$(q_{ii}+1)(q_{ii}q_{ij}q_{ji}-1)=0$.
\end{lemma}

\section{Connected arithmetic root systems of rank four}
\label{sec-classrank4}

Let $E$ be a basis of $\ndZ ^d$, where $d\ge 2$, and let $\chi $ be
a bicharacter on $\ndZ ^d$. The following proposition is one of the key
tools in the classification of arithmetic root systems. In most of the
cases it will be used without referring to it.

\begin{satz}{\cite[Prop.\,4,Lemma\,5]{a-Heck05a}}
\label{s-rootsubs}
Let $r\in \ndN $ with $r<d$. Assume that
$F=\{\aNdbasis _1,\ldots ,\aNdbasis _r\}\subset \roots ^+_E$
is a set of linearly independent elements of $\ndZ ^d$
such that for all $j\le r$ and all $m_1,\ldots ,m_{j-1}\in \ndN _0$
one has
\begin{align}\label{eq-subscheck}
 \aNdbasis _j-\sum _{i=1}^{j-1}m_i\aNdbasis _i\notin
 \roots \setminus \{\aNdbasis _j\}.
\end{align}
Then $\roots (\chi ;\aNdbasis _1,\ldots \aNdbasis _r)=(\roots \cap \ndZ F,
\chi \restr{\ndZ F\times \ndZ F},F)$ is an arithmetic root system.
\end{satz}

In the following also the terminology
``$\roots (\chi ;\aNdbasis _1,\ldots \aNdbasis _r)$ is finite''
will be used in order to emphasize the finiteness of the set
$\bigcup \{F'\,|\,(\id ,F')\in W_{\chi ',F}\}=\roots \cap \ndZ F$,
where $\chi '=\chi \restr{\ndZ F\times \ndZ F}$.

Note that relation (\ref{eq-subscheck}) holds in particular if
\begin{align}\label{eq-subscheck2}
 \aNdbasis _j-\sum _{i=1}^{j-1}m_i\aNdbasis _i\notin
 (\ndN _0E\setminus \{\aNdbasis _j\})\cup -\ndN _0E.
\end{align}

In the remaining part of this section let
$(\roots ,\chi ,E)$ be a connected rank 4 arithmetic root system.

\begin{lemma}\label{l-notetraeder}
One has $\prod _{i\not=j}(q_{ij}q_{ji}-1)=0$.
\end{lemma}

\begin{bew}
Assume that $q_{ij}q_{ji}\not=1$ for all $i,j\in \{1,2,3,4\}$ with $i\not=j$.
Since $\roots (\chi ;\Ndbasis _1,\Ndbasis _2,\Ndbasis _3+\Ndbasis _4)$
is finite by Proposition~\ref{s-rootsubs},
one obtains the equation
$\prod _{i=2}^4q_{1i}q_{i1}\prod _{i=3}^4q_{2i}q_{i2}=1$
from Lemma~\ref{l-specrank3}(ii).
Again by Lemma~\ref{l-specrank3}(ii) and the finiteness
of $\roots (\chi ;\Ndbasis _1,\Ndbasis _2,\Ndbasis _3)$ and
$\roots (\chi ;\Ndbasis _1,\Ndbasis _2,\Ndbasis _4)$
this yields $1=q_{14}q_{41}q_{24}q_{42}=(q_{12}q_{21})^{-1}$
which is a contradiction.
\end{bew}

\begin{lemma}\label{l-nodoubletriangle}
For given $a,b\in \{1,2,3,4\}$ with $a\not=b$ one has the equation
$\prod _{i\not=a}(q_{ai}q_{ia}-1)\prod _{i\not=b}(q_{bi}q_{ib}-1)=0$.
\end{lemma}

\begin{bew}
Without loss of generality assume that $a=1$ and $b=2$ and that the claim
of the lemma is false. By the previous lemma one has $q_{34}q_{43}=1$.
Consider $\roots (\chi ;\Ndbasis _1,\Ndbasis _2+\Ndbasis _3,\Ndbasis _4)$.
Then one gets $\left(\prod _{i=2}^4q_{1i}q_{i1}\right)q_{24}q_{42}=1$
which is a contradiction to $q_{13}q_{31}\not=1$ and the finiteness
of $\roots (\chi ;\Ndbasis _1,\Ndbasis _2,\Ndbasis _4)$.
\end{bew}

\begin{lemma}\label{l-no4cycle}
The graph $\gDd _{\chi ,E}$ is not a labeled cycle graph.
\end{lemma}

\begin{bew}
Assume that $(q_{12}q_{21}-1)(q_{23}q_{32}-1)(q_{34}q_{43}-1)
(q_{14}q_{41}-1)\not=0$ and $q_{ij}q_{ji}=1$ if $|i-j|\ge 2$.
Then the finiteness of $\roots (\chi ;\Ndbasis _1+\Ndbasis _2,\Ndbasis _3,
\Ndbasis _4)$ implies that $q_{23}q_{32}q_{34}q_{43}q_{14}q_{41}=1$.
By symmetry this gives that $q_{ij}q_{ji}=q_{12}q_{21}q_{23}q_{32}
q_{34}q_{43}q_{14}q_{41}$ whenever $|i-j|=1$, and $q_{12}q_{21}\in R_3$.
Moreover, from Lemma~\ref{l-specrank3}(ii) one obtains that
$(q_{11}q_{12}q_{21}q_{22}+1)(q_{33}+1)(q_{44}+1)=0$.
By Lemma~\ref{l-q12=-1}
one gets that $\prod _{i=1}^4(q_{ii}+1)=0$. Without loss of generality
assume that $q_{11}=-1$. Then one can apply $s_{\Ndbasis _1,E}$, and
hence Lemma~\ref{l-nodoubletriangle} implies the claim.
\end{bew}

\begin{lemma}\label{l-rightofway}
If $\gDd _{\chi ,E}$ is
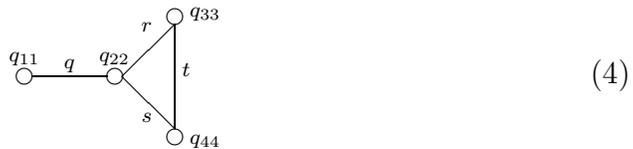
\begin{align}\label{eq-rightofway}
\Drightofway{}{$q_{11}$}{$q$}{$q_{22}$}{$r$}{$s$}{$q_{33}$}{$t$}{$q_{44}$}
\end{align}
then one of the two systems of equations
$q_{22}+1=(qr-1)(qs-1)=0$, $q_{22}q-1=(q-r)(q-s)=0$ is valid.
\end{lemma}

\begin{bew}
The finiteness of $\roots (\chi ;\Ndbasis _1,\Ndbasis _2+\Ndbasis _3,
\Ndbasis _2+\Ndbasis _4)$ implies that equation
$(q_{22}^2-1)(q_{22}^2q^2-1)=0$ holds.
Similarly, considering $\roots (\chi ;\Ndbasis _1+\Ndbasis _2,\Ndbasis _2
+\Ndbasis _3,\Ndbasis _4)$ one obtains that
$(q_{22}^2qr-1)(q_{22}^2qs-1)=0$.
Thus one has two possibilities. In the first case one 
has $q_{22}=-1$ and $(qr-1)(qs-1)=0$, and in the second the equations
$q_{22}^2q^2=1$ and $(q-r)(q-s)=0$ hold. Suppose that $q_{22}q=-1$
and $q_{22}^2qr=1$. Then $q_{22}\not=-1$, and hence Lemma~\ref{l-specrank3}(vii)
applied to
$\roots (\chi ;\Ndbasis _1,\Ndbasis _2,\Ndbasis _3)$ gives that $q_{22}qr=-1$.
This is a contradiction to $r\not=1$.
\end{bew}

\begin{lemma}\label{l-rightofwayhas-1}
If $\gDd _{\chi ,E}$ is
as in Figure (\ref{eq-rightofway}) then $(q_{33}+1)(q_{44}+1)=0$.
In particular, $(\roots ,\chi ,E)$ is Weyl equivalent to an arithmetic
root system $(\roots ',\chi ',E)$ such that $\gDd _{\chi ',E}$
is a labeled path graph.
\end{lemma}

\begin{bew}
Suppose that $(q_{33}+1)(q_{44}+1)\not=0$. By the finiteness of
$\roots (\chi ;\Ndbasis _2,\Ndbasis _3,\Ndbasis _4)$ one gets
$q_{22}=-1$, $q_{33}r=1$, $q_{44}s=1$, $r^2=s^2$, and $r\in R_3\cup R_6$.
Further, the finiteness of
$\roots (\chi ;\Ndbasis _1+\Ndbasis _2,\Ndbasis _3,\Ndbasis _4)$ gives that
$q_{11}q=1$. By Lemma~\ref{l-rightofway} one can assume that
$qr=1$. Since $\chi \chi \op (\Ndbasis _2,\Ndbasis _1+\Ndbasis _2+\Ndbasis _3
+\Ndbasis _4)=s\not=1$, one has $\Ndbasis _1+2\Ndbasis _2+\Ndbasis _3
+\Ndbasis _4\in \roots $. Then one obtains a contradiction to the finiteness
of $\roots (\chi ;\Ndbasis _1+2\Ndbasis _2+\Ndbasis _3+\Ndbasis _4,
\Ndbasis _3+\Ndbasis _4)$ \Dchaintwo{}{$q$}{$q^4$}{$q^4$} and
relation $q\in R_3\cup R_6$.

To obtain the last assertion of the lemma assume that $q_{33}=-1$ and apply
$s_{\Ndbasis _3,E}$.
\end{bew}

\begin{lemma}\label{l-3forkcond1}
If $\gDd _{\chi ,E}$ is
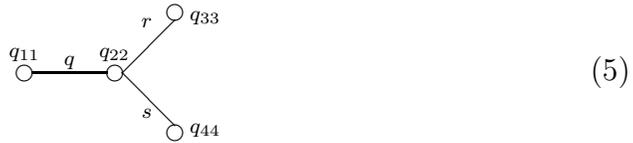
\begin{align}\label{eq-3fork}
\Dthreefork{}{$q_{11}$}{$q$}{$q_{22}$}{$r$}{$s$}{$q_{33}$}{$q_{44}$}
\end{align}
then one has either $q_{22}=-1$, and at least two of $qr,qs$ and $rs$ are
equal to 1, or at least two of $q_{22}q,q_{22}r$ and $q_{22}s$ are equal to 1.
\end{lemma}

\begin{bew}
Assume first that $q_{22}=-1$. Then one can apply $s_{\Ndbasis _2,E}$
and hence Lemma~\ref{l-nodoubletriangle} gives the claim.

Consider now the case $q_{22}\not=-1$.
Assume first that $q_{11}qq_{22}=-1$. By Lemma~\ref{l-q12=-1}
this implies that $q_{11}=-1$ and $q_{22}q=1$.
Then one can apply $s_{\Ndbasis _1,E}$, and the previous case
implies that $(q^{-1}r-1)(q^{-1}s-1)=0$.

Finally, by twist equivalence one can assume that $q_{11}qq_{22},
q_{22}rq_{33}$ and $q_{22}sq_{44}$ are different from $-1$. Then the finiteness
of $\roots (\chi ;\Ndbasis _1+\Ndbasis _2,\Ndbasis _2+\Ndbasis _3,
\Ndbasis _2+\Ndbasis _4)$ and Lemma~\ref{l-specrank3}(ii)
imply that at least one of
$q_{22}^2qr,q_{22}^2qs$ and $q_{22}^2rs$ is equal to 1. Without loss of
generality suppose that $q_{22}^2qr=1$. Then relation $q_{22}\not=-1$
gives that $q_{22}qr\not=-1$, and hence Lemma~\ref{l-specrank3}(vii) applied to
$\roots (\chi ;\Ndbasis _1,\Ndbasis _2,\Ndbasis _3)$ implies that
$q_{22}q=q_{22}r=1$.
\end{bew}

\begin{satz}\label{s-Weqpath}
If $(\roots ,\chi ,E)$ is not of Cartan type then it is Weyl equivalent to
an arithmetic root system $(\roots ',\chi ',E)$ such that $\gDd _{\chi ',E}$
is a labeled path graph.
\end{satz}

\begin{bew}
By Lemmata \ref{l-notetraeder}, \ref{l-nodoubletriangle}, \ref{l-no4cycle}
and \ref{l-rightofwayhas-1} it suffices to show the claim under the assumption
that $\gDd _{\chi ,E}$ is of the form (\ref{eq-3fork}).

If $q_{22}=-1$ then by Lemma~\ref{l-3forkcond1} one can assume that
$qr=1$ and $qs=1$. If $q=-1$ then by \cite[Lemma\,16]{a-Heck05a}
one can suppose additionally
that $q_{33}=q_{44}=-1$. Then the finiteness of $\roots (\chi;\Ndbasis _1,\Ndbasis _2,
\Ndbasis _3)$ implies that $(\chi ,E)$ is of Cartan type which is a contradiction.
On the other hand, if $q\not=-1$ then one can apply $s_{\Ndbasis _2,E}$ and
Lemma~\ref{l-rightofwayhas-1} gives the claim.

Assume now that $q_{22}\not=-1$.
Further, by Lemma~\ref{l-3forkcond1} one can suppose
additionally that $q_{22}r=q_{22}s=1$. Then
Lemma~\ref{l-specrank3}(ii) and the finiteness of
$\roots (\chi ;\Ndbasis _3,\Ndbasis _1+\Ndbasis _2,
\Ndbasis _2+\Ndbasis _4)$
give that equation
\begin{align}\label{eq-3forkcond2}
(q_{22}q-1)(q_{22}-q)=0
\end{align}
holds. Moreover, since $(\chi ,E)$ is not of Cartan type,
\cite[Cor.~13]{a-Heck05a} implies that one of
$\roots (\chi ;\Ndbasis _1,\Ndbasis _2)$,
$\roots (\chi ;\Ndbasis _2,\Ndbasis _3)$ and
$\roots (\chi ;\Ndbasis _2,\Ndbasis _4)$
has a generalized Dynkin diagram appearing in rows~3, 5, 6, or 7 of
\cite[Table\,1]{a-Heck05a}.
One has the following possibilities.

(i) $q_{33}=-1$. Apply $s_{\Ndbasis _3,E}$.
The claim follows from the second paragraph
of the proof of this lemma.

(ii) $q_{11}=-1$, $q=q_{22}^{-2}$. Apply $s_{\Ndbasis _1,E}$ and use
Lemma~\ref{l-3forkcond1} to obtain that
$-q_{22}^{-1}r=1$, that is $q_{22}\in R_4$.
This is a contradiction to Equation~(\ref{eq-3forkcond2}).

(iii) $q_{22}\in R_3$, $q=-q_{22}$, $q_{11}=-1$.
This is again a contradiction to (\ref{eq-3forkcond2}).

(iv) $q_{22}\in R_3$, $q_{11}q=1$, $q^3\not=1$.
Recall again (\ref{eq-3forkcond2}).

(v) $q_{33}\in R_3$, $q_{22}^3\not=1$.
Since $q_{33}\not=-1$ and $q_{33}r\not=1$,
Lemma~\ref{l-specrank3}(ix) and the finiteness of
$\roots (\chi ;\Ndbasis _3,\Ndbasis _1+\Ndbasis _2,
\Ndbasis _2+\Ndbasis _4)$ imply that $q_{22}q=1$, that is
one has $q=r=s=q_{22}^{-1}$.
Thus Lemma~\ref{l-specrank3}(vii) applied to
$\roots (\chi ;\Ndbasis _3,\Ndbasis _1+\Ndbasis _2,
\Ndbasis _2+\Ndbasis _4)$ gives that equation
$q_{33}q^2=-1$ is fulfilled. Apply now the transformation
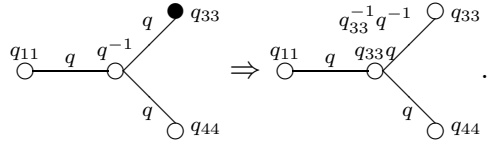
\begin{align*}
\Dthreefork{t}{$q_{11}$}{$q$}{$q^{-1}$}{$q$}{$q$}{$q_{33}$}{$q_{44}$}
\Rightarrow
\Dthreefork{}{$q_{11}$}{$q$}{$q_{33}q$}{$q_{33}^{-1}q^{-1}$}{$q$}%
{$q_{33}$}{$q_{44}$}.
\end{align*}
If $q_{33}q=-1$ then the claim follows from the second paragraph of the proof
of this lemma. Otherwise Lemma~\ref{l-3forkcond1} gives that $q_{33}q^2=1$
which is a contradiction to equation $q_{33}q^2=-1$.
\end{bew}

\begin{thm}\label{t-classrank4}
Let $k$ be a field of characteristic 0.
Then twist equivalence classes of connected arithmetic root systems
of rank 4 are in one-to-one correspondence to generalized Dynkin diagrams
appearing in Table~\ref{a-r4}.
Moreover, two such arithmetic root systems are
Weyl equivalent if and only if their generalized Dynkin diagrams appear
in the same row of Table~\ref{a-r4} and can be presented with the same set of
fixed parameters.
\end{thm}

\begin{bew}
One can check case by case that each row of Table~\ref{a-r4} contains the
generalized Dynkin diagrams of the representants of a unique Weyl
equivalence class $(\chi ,E)$, where $\chi $ is a bicharacter on $\ndZ ^4$
and $E$ is a fixed basis of $\ndZ ^4$. In order to prove that these diagrams
correspond to arithmetic root systems one has to show that $\extWBG _{\chi ,E}$
is finite. If $(\chi ,E)$ is of Cartan type then this follows
from \cite[Theorem\,1]{a-Heck04c}. In all other cases
Proposition~\ref{s-Wfinite} will be applied.
One checks first that $\extWBG _{\chi ,E}$ is full and using
\cite[Theorem\,12]{a-Heck05a} one recognizes that
$\extWBG _{\chi ,E'\subset E}$ is finite for all subsets
$E'\subset E$ of 3 elements.
It remains to find an element $(T,E)$ with respect to one
fixed representant of the Weyl equivalence class
which has the property
$T(E)\subset E\setminus \{\Ndbasis \}\cup \{-\aNdbasis \}$ for some
$\Ndbasis \in E$ and $\aNdbasis \in \roots ^+_E$. These elements
$(T,E)$ are listed in Appendix~\ref{a-(T,E)}.

It remains to prove that all arithmetic root systems
have a generalized
Dynkin diagram appearing in Table~\ref{a-r4}. By
\cite[Theorem\,1]{a-Heck04c}
and Proposition~\ref{s-Weqpath}
it suffices to consider pairs $(\chi ,E)$
such that $\gDd _{\chi ,E}$ is a labeled path graph
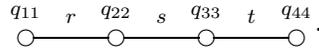
\begin{align*}
\Dchainfour{}{$q_{11}$}{$r$}{$q_{22}$}{$s$}{$q_{33}$}{$t$}{$q_{44}$}~.
\end{align*}
In what follows several cases will be considered according to the
set $I=\{i\in \{1,2,3,4\}\,|\,q_{ii}=-1\}$.

\textit{Step~1.} $q_{ii}^2\not=1$ for all $i\in \{1,2,3,4\}$.
By Lemma~\ref{l-specrank3}(i) for
$\roots (\chi ;\Ndbasis _1,\Ndbasis _2,\Ndbasis _3)$ and
$\roots (\chi ;\Ndbasis _2,\Ndbasis _3,\Ndbasis _4)$
one obtains that there exist $m_1,m_2,m_3,m_4\in \{1,2\}$
such that $q_{22}^{m_1}r=q_{22}^{m_2}s=q_{33}^{m_3}s=q_{33}^{m_4}t=1$.
More precisely, either $(\chi ,E)$
is of Cartan type, in which case $\gDd _{\chi ,E}$ appears in rows
1--5 of Table~\ref{a-r4},
or one has without loss of generality $q_{11}\in R_3$, $q_{22}r=1$,
$q_{33}s=1$, $(q_{22}s-1)(q_{22}^2s-1)=0$, and $q_{22}^3,q_{33}^3\not=1$.
If $q_{22}^2s=1$ then Lemma~\ref{l-specrank3}(i) for
$\roots (\chi ;\Ndbasis _2,\Ndbasis _3,\Ndbasis _4)$ gives that
$q_{33}t=q_{44}t=1$. In this case the finiteness of $\roots (\chi ;
2\Ndbasis _1+\Ndbasis _2,2\Ndbasis _2+\Ndbasis _3,
\Ndbasis _3+\Ndbasis _4)$~
\Dtriangle{}{$q_{11}r$}{$r^{-2}$}{$r^{-2}$}{$r^2$}{$r^2$}{$r^2$}
\rule{0pt}{15\unitlength}~,
relation $q_{33}^2=r^{-6}\not=1$ and Lemma~\ref{l-specrank3}(ii)
give a contradiction.

On the other hand, if $q_{22}s=1$ then Lemma~\ref{l-specrank3}(i) for
$\roots (\chi ;\Ndbasis _1+\Ndbasis _2,\Ndbasis _3,\Ndbasis _4)$
\Dchainthree{}{$q_{11}$}{$r$}{$r^{-1}$}{$t$}{$q_{44}$}~
implies that $q_{44}t=1$, and either $r=t$, $r^3\not=1$, or $r^2=t$,
$r^6\not=1$. Consider first the case
$r^2=t$. Applying Lemma~\ref{l-specrank3}(i) to $\roots (\chi ;
\Ndbasis _2+\Ndbasis _3,2\Ndbasis _1+\Ndbasis _2,2\Ndbasis _3+\Ndbasis _4)$
\Dchainthree{}{$r^{-1}$}{$r$}{$q_{11}r$}{$r^2$}{$r^{-2}$}~
one obtains that
$q_{11}r=-1$ or $q_{11}r^2=1$ or $q_{11}r^3-1=q_{11}^2r^3-1=0$
which is a contradiction to 
$q_{11}\in R_3$ and $r^6\not=1$. In the last case one has $r=t$.
Apply now Lemma~\ref{l-specrank3}(i) to $\roots (\chi ;
\Ndbasis _2+\Ndbasis _3,2\Ndbasis _1+\Ndbasis _2,\Ndbasis _3+\Ndbasis _4)$
\Dchainthree{}{$r^{-1}$}{$r$}{$q_{11}r$}{$r$}{$r^{-1}$}~ and conclude
that $(q_{11}r+1)(q_{11}r^2-1)=0$.
If $q_{11}r=-1$ then $\gDd _{\chi ,E}$
appears in row~17 of Table~\ref{a-r4}.
Otherwise relation $r^3\not=1$ implies that
$r=-q_{11}$, and hence $\gDd _{\chi ,E}$ appears
in row~15 of Table~\ref{a-r4}.

\textit{Step~2.} $q_{11}=-1$, $q_{ii}\not=-1$
for all $i\in \{2,3,4\}$.
By Lemma~\ref{l-specrank3}(ii) for
$\roots (\chi ;\Ndbasis _1+\Ndbasis _2,
\Ndbasis _2+\Ndbasis _3,\Ndbasis _3+\Ndbasis _4)$~~
\Dtriangle{}{$-q_{22}r$}{$q_{22}sq_{33}$}%
{$q_{33}tq_{44}$}{$rq_{22}^2s$}{$sq_{33}^2t$}{$s$}
\rule{0pt}{15\unitlength}~~~~~
and Lemma~\ref{l-q12=-1} one obtains that
either equation $(rq_{22}^2s-1)(sq_{33}^2t-1)=0$ is fulfilled, or
relations $q_{22}sq_{33},q_{33}tq_{44}\not=-1$,
$q_{22}r=1$ and $(rq_{22}^2s-s)(rq_{22}^2s+s)=0$ hold.
The last case is a contradiction to $q_{22}^2\not=1$.

Step~2.1.
Assume first that $rq_{22}^2s=1$. Since $q_{22}\not=-1$,
Lemma~\ref{l-specrank3}(vii) for
$\roots (\chi ;\Ndbasis _1,\Ndbasis _2,\Ndbasis _3)$
gives that $q_{22}r=q_{22}s=1$.
Further, Lemma~\ref{l-specrank3}(i) for 
$\roots (\chi ;\Ndbasis _2,\Ndbasis _3,\Ndbasis _4)$
implies that $(q_{33}s-1)(q_{33}^2s-1)=0$.

If $q_{22}r=q_{22}s=q_{33}^2s=1$ then again by
Lemma~\ref{l-specrank3}(i) for 
$\roots (\chi ;\Ndbasis _2,\Ndbasis _3,\Ndbasis _4)$
one obtains that $q_{33}t=1$, and either $q_{44}\in R_3$, $t^6\not=1$,
or $q_{44}t=1$. Note that also relation $t^4=q_{22}^{-2}\not=1$ holds.
Therefore Lemma~\ref{l-specrank3}(iii) for 
$\roots (\chi ;\Ndbasis _2+2\Ndbasis _3,\Ndbasis _1+\Ndbasis _2,
\Ndbasis _3+\Ndbasis _4)$
\Dchainthree{}{$t^{-2}$}{$t^2$}{$-1$}{$t^2$}{$q_{44}$}~ implies that
$q_{44}=t^{-1}\in R_{10}$.
This is a contradiction to the finiteness of
$\roots (\chi ;\Ndbasis _2+\Ndbasis _3,\Ndbasis _3+\Ndbasis _4,
\Ndbasis _1+\Ndbasis _2)$
\Dchainthree{}{$t^{-1}$}{$t$}{$q_{44}$}{$t^2$}{$-1$}
and Lemma~\ref{l-specrank3}(iv).

On the other hand, if $q_{22}r=q_{22}s=q_{33}s=1$ then by
Lemma~\ref{l-specrank3}(i) for 
$\roots (\chi ;\Ndbasis _2,\Ndbasis _3,\Ndbasis _4)$
one obtains that either equations $q_{33}^2t=1$, $q_{44}t=1$ hold,
or equation $q_{33}t=1$ is valid. In the first case $\gDd _{\chi ,E}$
appears in row~8 of Table~\ref{a-r4}.
In the second case the finiteness of
$\roots (\chi ;\Ndbasis _1,\Ndbasis _2,\Ndbasis _3+\Ndbasis _4)$
\Dchainthree{}{$-1$}{$t$}{$t^{-1}$}{$t$}{$q_{44}$}~
and Lemma~\ref{l-specrank3}(v) imply
that $q_{44}t=1$ or $q_{44}^2t=1$ or $q_{44}=-t\in R_3$.
Then $\gDd _{\chi ,E}$ appears in rows~6, 7, and 16
of Table~\ref{a-r4}, respectively.

Step~2.2.
Now suppose that $rq_{22}^2s\not=1$ and $sq_{33}^2t=1$. Since
relation $q_{33}\not=-1$ holds,
Lemma~\ref{l-specrank3}(vii) for
$\roots (\chi ;\Ndbasis _2,\Ndbasis _3,\Ndbasis _4)$
gives $q_{33}s=q_{33}t=1$.
Further, Lemma~\ref{l-specrank3}(vii) for 
$\roots (\chi ;\Ndbasis _1,\Ndbasis _2,\Ndbasis _3)$
yields $(q_{22}r-1)(q_{22}s-1)(q_{22}rs+1)=0$.
If $(q_{22}r-1)(q_{22}s-1)\not=0$ then one can change the representant
of the Weyl equivalence class via the transformation
\begin{align*}
\Dchainfour{1}{$-1$}{$r$}{$q_{22}$}{$s$}{$q_{33}$}{$t$}{$q_{44}$}
\quad \Rightarrow \quad 
\Dchainfour{}{$-1$}{$r^{-1}$}{$-rq_{22}$}{$s$}{$q_{33}$}{$t$}{$q_{44}$}
\end{align*}
and therefore one can assume that $(q_{22}r-1)(q_{22}s-1)=0$.

In the first case one has $q_{22}r=1$, and hence $r\not=s$.
Lemma~\ref{l-specrank3}(i) for
$\roots (\chi ;\Ndbasis _2,\Ndbasis _3,\Ndbasis _4)$
implies that $q_{33}=q_{44}$ and either $q_{22}\in R_3$
or $q_{22}^2s=1$. If $q_{22}^2s=1$ then $\gDd _{\chi ,E}$ appears in
row~9 or 18 of Table~\ref{a-r4}.
Otherwise one has the relations $r\in R_3$, $s^3\not=1$, and hence
Lemma~\ref{l-specrank3}(ii) for
$\roots (\chi ;\Ndbasis _1+\Ndbasis _2,2\Ndbasis _2+\Ndbasis _3,
\Ndbasis _3+\Ndbasis _4)$~
\mbox{\Dtriangle{}{$-1$}{$r^{-1}s$}{$s^{-1}$}{$rs$}{$s$}{$s$}
\rule{0pt}{15\unitlength}}~
implies that $rs^3=1$, and either $r^{-1}s=-1$ or $r^2s^2=s^2$.
This is impossible since $r\in R_3$.

In the second case equation $q_{22}s=1$ holds, and hence relation
$rq_{22}^2s\not=1$ implies that $r\not=s$. Therefore
Lemma~\ref{l-specrank3}(iv) for
$\roots (\chi ;\Ndbasis _1,\Ndbasis _2,\Ndbasis _3)$
tells that $r=s^2$ and $s\in R_3\cup R_4\cup R_6$.
If $s\in R_4$ then Lemma~\ref{l-specrank3}(i) for
$\roots (\chi ;\Ndbasis _2,\Ndbasis _3,\Ndbasis _4)$
implies that $(\chi ,E)$ is of Cartan type
and hence $\gDd _{\chi ,E}$ appears in row~3 of Table~\ref{a-r4}.
On the other hand, if $s\in R_3\cup R_6$ then
Lemma~\ref{l-specrank3}(i) for
$\roots (\chi ;\Ndbasis _2+\Ndbasis _3,\Ndbasis _1+\Ndbasis _2,
\Ndbasis _3+\Ndbasis _4)$~
\Dchainthree{}{$s^{-1}$}{$s$}{$-s$}{$s$}{$q_{44}$}
gives a contradiction.

\textit{Step~3.} $q_{22}=-1$, $q_{ii}\not=-1$
for all $i\in \{1,3,4\}$.

Step~3.1. Assume first that $q_{11}r\not=1$. Then by
\cite[Lemma\,16]{a-Heck05a} for
$\roots (\chi ;\Ndbasis _1,\Ndbasis _2,\Ndbasis _3)$ and
$\roots (\chi ;\Ndbasis _1,\Ndbasis _2+\Ndbasis _3,\Ndbasis _4)$
one has $q_{33}s=q_{44}t=1$. Further, \cite[Lemma\,16]{a-Heck05a} for
$\roots (\chi ;\Ndbasis _1,\Ndbasis _2+\Ndbasis _3,
\Ndbasis _3+\Ndbasis _4)$~
\Dchainthree{}{$q_{11}$}{$r$}{$-1$}{$s^{-1}t$}{$s^{-1}$}
implies that $(s-t)(s^2-t)=0$.

If $s^2=t$ then relation $s^4\not=1$ and Lemma~\ref{l-specrank3}(iii) for
$\roots (\chi ;\Ndbasis _1,\Ndbasis _2,\Ndbasis _3)$ and
$\roots (\chi ;\Ndbasis _1,\Ndbasis _2+\Ndbasis _3,\Ndbasis _4)$
imply that $q_{11}^2r=1$ and $(s^{-1}q_{11}^2-1)(s^{-1}q_{11}^3+1)=
(s^{-2}q_{11}^2-1)(s^{-2}q_{11}^3+1)=0$. The last equations contradict
to the assumption $s^4\not=1$ and $q_{11}^2\not=1$.

If $s=t$ then consider first two special cases. If $q_{11}\in R_3$,
$r^3\not=1$, and $q_{11}r^2\not=1$, then the transformation
\begin{align*}
\Dchainfour{1}{$q_{11}$}{$r$}{$-1$}{$s$}{$s^{-1}$}{$s$}{$s^{-1}$}
\quad \Rightarrow \quad
\Dchainfour{}{$q_{11}~~~$}{$q_{11}^2r^{-1}~$}{$~~~~~-q_{11}r^2$}%
{$~~~~s$}{$~~s^{-1}$}{$s$}{$s^{-1}$}
\end{align*}
shows that $(\roots ,\chi ,E)$
is Weyl equivalent to an arithmetic root
system appearing in Step~1.

On the other hand, if $rs=1$ then the transformations
\begin{align*}
&\Dchainfour{2}{$q_{11}$}{$s^{-1}$}{$-1$}{$s$}%
{$s^{-1}$}{$s$}{$s^{-1}$}
&&\Rightarrow &
&\Dchainfour{3}{$-q_{11}s^{-1}$}{$s$}{$-1$}%
{$s^{-1}$}{$-1$}{$s$}{$s^{-1}$}
\quad \Rightarrow \\
&\Dchainfour{4}{$-q_{11}s^{-1}$}{$s$}{$s^{-1}$}%
{$s$}{$-1$}{$s^{-1}$}{$-1$}
&&\Rightarrow &
&\Dchainfour{}{$-q_{11}s^{-1}$}{$s$}{$s^{-1}$}%
{$s$}{$s^{-1}$}{$s$}{$-1$}
\end{align*}
and relation $q_{11}s^{-1}=q_{11}r\not=1$
show that $(\roots ,\chi ,E)$
is Weyl equivalent to an arithmetic root system
appearing in Step~2.

Finally, apply Lemma~\ref{l-specrank3}(iii) to
$\roots (\chi ;\Ndbasis _3,\Ndbasis _2,\Ndbasis _1)$
and conclude that either 
$q_{11}^2r=1$, $(s^{-1}q_{11}^2-1)(s^{-1}q_{11}^3+1)=0$
or $q_{11}=-s^{-1}\in R_3$, $r\in \{-1,-q_{11}\}$.
If $q_{11}^2r=q_{11}^2s^{-1}=1$ then one gets $rs=1$,
and if $q_{11}^2r=-q_{11}^3s^{-1}=1$ then
$\roots (\chi ;\Ndbasis _1+\Ndbasis _2,\Ndbasis _3,\Ndbasis _4)$
\Dchainthree{}{$-q_{11}^{-1}$}{$-q_{11}^3$}%
{$-q_{11}^{-3}$}{$-q_{11}^3$}{$-q_{11}^{-3}$}~
is of infinite Cartan type. It remains to check
the case $q_{11}=-s^{-1}\in R_3$, $r^3\not=1$. As noted above one can assume
also that $q_{11}r^2=1$ and $rs\not=1$. However relations
$q_{11}=-s^{-1}\in R_3$ and $q_{11}r^2=1$ imply that $(rs-1)(rs+1)=0$ which is
a contradiction to relations $-s\in R_3$, $r^3\not=1$ and $rs\not=1$.

Step~3.2. Suppose now that $q_{11}r=1$. Then $r^2=q_{11}^{-2}\not=1$.

Step~3.2.1. Consider the case when $rs=1$. Then the transformations
\begin{align*}
\Dchainfour{2}{$r^{-1}$}{$r$}{$-1$}{$r^{-1}$}{$q_{33}$}{$t$}{$q_{44}$}
\Rightarrow
\Dchainfour{1}{$-1$}{$r^{-1}$}{$-1$}{$r~~~$}{$-q_{33}r^{-1}$}{$~~~t$}{$q_{44}$}
\Rightarrow
\Dchainfour{}{$-1$}{$r$}{$r^{-1}$}{$r~~~$}{$-q_{33}r^{-1}$}{$~~~t$}{$q_{44}$}
\end{align*}
show that if $q_{33}\not=r$ then
$(\roots ,\chi ,E)$ is Weyl equivalent
to an arithmetic root system appearing in Step~2.
Therefore one can assume
that $q_{33}=r$. Now use that
$\roots (\chi ;\Ndbasis _2,\Ndbasis _3,\Ndbasis _4)$
is finite and $q_{33},q_{44}\not=-1$ and apply
\cite[Theorem\,12]{a-Heck05a}.
If $r=t^{-1}=q_{44}$ then $\gDd _{\chi ,E}$ appears in row~10 of
Table~\ref{a-r4}.
If $q_{44}^2t=q_{33}t=1$ then $\gDd _{\chi ,E}$ can be found in row~11
of Table~\ref{a-r4}.
If $q_{44}t=q_{33}^2t=1$ then $\gDd _{\chi ,E}$ appears in row~12 of
Table~\ref{a-r4}.
If $q_{44}t=q_{33}^3t=1$ then Lemma~\ref{l-specrank3}(iii) applied to
$\roots (\chi ;\Ndbasis _1,\Ndbasis _2,2\Ndbasis _3+\Ndbasis _4)$
\Dchainthree{}{$r^{-1}$}{$r$}{$-1$}{$r^{-2}$}{$r$}~
implies that $r\in R_4$ and then one can find $\gDd _{\chi ,E}$
in row~22 of Table~\ref{a-r4}. If $q_{44}=-t=-r^{-1}\in R_3$ then
row~19 of Table~\ref{a-r4} contains
$\gDd _{\chi ,E}$. Finally, if $r=-t^{-1}=-q_{44}\in R_3$
then Lemma~\ref{l-specrank3}(iii) tells that
$\roots (\chi ;\Ndbasis _1,\Ndbasis _2+\Ndbasis _3,\Ndbasis _3+\Ndbasis _4)$
\Dchainthree{}{$r^{-1}$}{$r$}{$-1$}{$-1$}{$r$}
is not finite.

Step~3.2.2. Assume now that $rs\not=1$. Note that
Lemma~\ref{l-q12=-1} applied to
$\roots (\chi ;\Ndbasis _3,\Ndbasis _4)$ implies also relation
$q_{33}tq_{44}\not=-1$. Then Lemma~\ref{l-specrank3}(ii) for
$\roots (\chi ;\Ndbasis _1+\Ndbasis _2,\Ndbasis _2+\Ndbasis _3,
\Ndbasis _3+\Ndbasis _4)$
\makebox{
\Dtriangle{}{$-1$}{$-q_{33}s$}{$q_{33}tq_{44}$}{$rs$}{$q_{33}^2st$}{$s$}
\rule{0pt}{15\unitlength}}~~~~~
gives that one of equations $q_{33}^2st=1$, $q_{33}s=1$ holds.
Therefore by
Lemma~\ref{l-specrank3}(vii)
for $\roots (\chi ;\Ndbasis _2,\Ndbasis _3,\Ndbasis _4)$
one obtains that $q_{33}s=1$. Further, Lemma~\ref{l-rightofway}
and the transformations
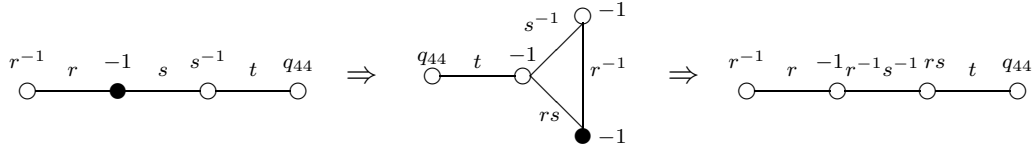
\begin{align*}
\Dchainfour{2}{$r^{-1}$}{$r$}{$-1$}{$s$}{$s^{-1}$}{$t$}{$q_{44}$}
\quad \Rightarrow \quad
\Drightofway{b}{$q_{44}$}{$t$}{$-1$}{$s^{-1}$}{$rs$}{$-1$}{$r^{-1}$}{$-1$}
\quad \Rightarrow \quad
\Dchainfour{}{$r^{-1}$}{$r$}{$-1~~$}{$r^{-1}s^{-1}$}{$~~rs$}{$t$}{$q_{44}$}
\end{align*}
give that $(rst-1)(s^{-1}t-1)=0$, and hence either $rs=-1$ or
(by Weyl equivalence) one can assume that $q_{33}t=1$.
Note that relations
$rs=-1$ and $s\not=t$ imply that $rst=1$ and hence $t=-1$
which is a contradiction
to $q_{33},q_{44}\not=-1$ and \cite[Lemma\,10(i)]{a-Heck05a}
for $\roots (\chi ;\Ndbasis _3,\Ndbasis _4)$. Therefore
one can suppose that $s=t$.

By relation $rs\not=1$ and Lemma~\ref{l-specrank3}(vi) for
$\roots (\chi ;\Ndbasis _1+\Ndbasis _2,\Ndbasis _2+\Ndbasis _3,\Ndbasis _4)$
\Dchainthree{}{$-1$}{$rs$}{$-1$}{$s$}{$q_{44}$}~
one obtains that $(rs+1)(q_{44}s-1)(rs^2-1)=0$. On the other hand,
Lemma~\ref{l-specrank3}(iii) for
$\roots (\chi ;\Ndbasis _1,\Ndbasis _2+\Ndbasis _3,\Ndbasis _4)$
\Dchainthree{}{$r^{-1}$}{$r$}{$-1$}{$s$}{$q_{44}$}~
together with relations $s^2\not=1$, $rs\not=1$ implies that
either $q_{44}s=1$ or $q_{44}^2s=1$, $q_{44}^3=-r$.
Since $q_{44}^2\not=1$, these two conditions give that $q_{44}s=1$.

If $rs=-1$ or $rs^2=1$ or $r^2s^3=1$ then $\gDd _{\chi ,E}$ appears in row~14,
13 and 9 of Table~\ref{a-r4}, respectively. Otherwise the transformations
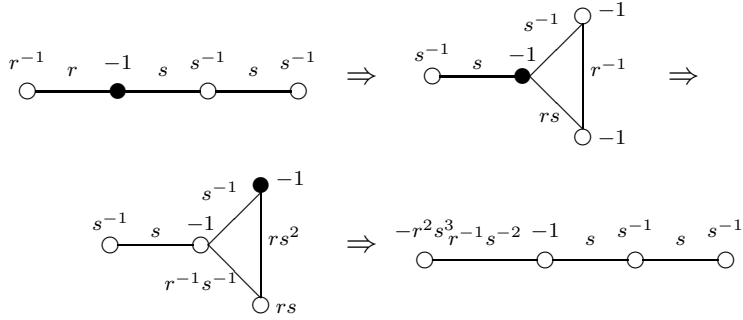
\begin{align*}
\Dchainfour{2}{$r^{-1}$}{$r$}{$-1$}{$s$}{$s^{-1}$}{$s$}{$s^{-1}$}
\quad &\Rightarrow \quad
\Drightofway{m}{$s^{-1}$}{$s$}{$-1$}{$s^{-1}$}{$rs$}{$-1$}{$r^{-1}$}{$-1$}
\quad \Rightarrow \quad \\
\Drightofway{t}{$s^{-1}$}{$s$}{$-1$}{$s^{-1}$}{$r^{-1}s^{-1}$}%
{$-1$}{$rs^2$}{$rs$}
\quad & \Rightarrow \quad
\bigDchainfour{}{$-r^2s^3$}{$r^{-1}s^{-2}$}{$-1$}{$s$}{$s^{-1}$}{$s$}{$s^{-1}$}
\end{align*}
show that $(\roots ,\chi ,E)$ is Weyl equivalent to an
arithmetic root system appearing in Step~3.1.

\textit{Step~4.} $q_{11}=q_{22}=-1$, $q_{33},q_{44}\not=-1$.
One can assume that $r=-1$, because otherwise the transformation
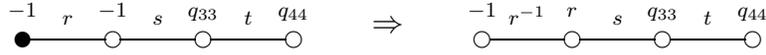
\begin{align*}
\Dchainfour{1}{$-1$}{$r$}{$-1$}{$s$}{$q_{33}$}{$t$}{$q_{44}$}
\qquad \Rightarrow \qquad
\Dchainfour{}{$-1$}{$r^{-1}$}{$r$}{$s$}{$q_{33}$}{$t$}{$q_{44}$}
\end{align*}
gives an arithmetic root system from Step~2.
Further, Lemma~\ref{l-specrank3}(vii) for
$\roots (\chi ;\Ndbasis _2,\Ndbasis _3,\Ndbasis _4)$
implies that equation $(q_{33}s-1)(q_{33}t-1)(q_{33}st+1)=0$ holds.

Step~4.1. Consider the case $q_{33}s=1$. Since $q_{33}^2\not=1$, one obtains
that $s^2\not=1$. Then using Lemma~\ref{l-rightofway} and the first of the
transformations
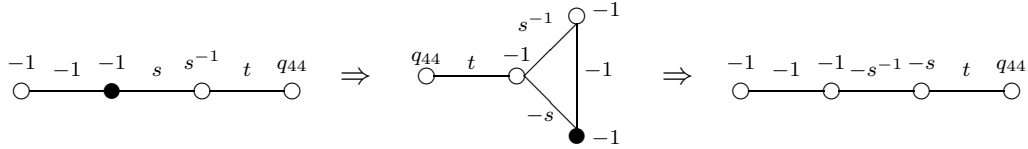
\begin{align*}
\Dchainfour{2}{$-1$}{$-1$}{$-1$}{$s$}{$s^{-1}$}{$t$}{$q_{44}$}
\quad \Rightarrow \quad
\Drightofway{b}{$q_{44}$}{$t$}{$-1$}{$s^{-1}$}{$-s$}{$-1$}{$-1$}{$-1$}
\quad \Rightarrow \quad
\Dchainfour{}{$-1$}{$-1$}{$-1$}{$-s^{-1}$}{$-s$}{$t$}{$q_{44}$}
\end{align*}
one obtains that $(s-t)(st+1)=0$ and hence by the second of the above
transformations one can assume that $s=t$.
Further, Lemma~\ref{l-specrank3}(i) applied to
$\roots (\chi ;\Ndbasis _3,\Ndbasis _4,
\Ndbasis _1+2\Ndbasis _2+\Ndbasis _3)$
\Dchainthree{}{$s^{-1}$}{$s$}{$q_{44}$}{$s$}{$-s$}
gives relations $q_{44}s=1$ and $s\in R_3\cup R_4\cup R_6$.
Then Lemma~\ref{l-specrank3}(vi) for
$\roots (\chi ;\Ndbasis _1+\Ndbasis _2,\Ndbasis _2+\Ndbasis _3,\Ndbasis _4)$
\Dchainthree{}{$-1$}{$-s$}{$-1$}{$s$}{$s^{-1}$}
gives that $s\in R_3\cup R_4$. In this case $\gDd _{\chi ,E}$ appears in
row~9 and 13 of Table~\ref{a-r4}, respectively.

Step~4.2. Assume that $q_{33}t=1$ and $q_{33}s\not=1$.
Then Lemma~\ref{l-specrank3}(ii) for
$\roots (\chi ;\Ndbasis _1+\Ndbasis _2,\Ndbasis _2+\Ndbasis _3,
\Ndbasis _3+\Ndbasis _4)$
\makebox{
\Dtriangle{}{$-1$}{$-st^{-1}$}{$q_{44}$}{$-s$}{$st^{-1}$}{$s$}
\rule{0pt}{15\unitlength}}~
implies that either $s=-1$ or $-s^3t^{-1}=s^2t^{-1}=1$. Therefore
it remains to consider the case $s=-1$. By
Lemma~\ref{l-specrank3}(vii) for
$\roots (\chi ;\Ndbasis _1,\Ndbasis _2+\Ndbasis _3,\Ndbasis _3+\Ndbasis _4)$
\Dchainthree{}{$-1$}{$-1$}{$t^{-1}$}{$-t^{-1}$}{$q_{44}$}
and since $t^2\not=1$ one gets $t\in R_4$.
If $q_{44}t\not=1$ then \cite[Lemma\,10]{a-Heck05a}(i) gives that
$\roots (\chi ;\Ndbasis _1,\Ndbasis _2+\Ndbasis _3,
\Ndbasis _3+2\Ndbasis _4)$
\Dchainthree{}{$-1$}{$-1$}{$t^{-1}$}{$-1$}{$tq_{44}^4$}
~is of infinite Cartan type, a contradiction. Hence one has
$q_{44}=t^{-1}\in R_4$.
In this case $\gDd _{\chi ,E}$ appears in row~4 of Table~\ref{a-r4}.

Step~4.3. Suppose now that $q_{33}st=-1$ and $q_{33}s\not=1\not=q_{33}t$.
Then $s\not=-1$, $t\not=-1$, $st\not=1$,
and by Lemma~\ref{l-q12=-1} for
$\roots (\Ndbasis _3,\Ndbasis _4)$
also relation $s^{-1}q_{44}\not=1$ holds. Hence by
Lemma~\ref{l-specrank3}(ii) for
$\roots (\chi ;\Ndbasis _1+\Ndbasis _2,\Ndbasis _2+\Ndbasis _3,
\Ndbasis _3+\Ndbasis _4)$
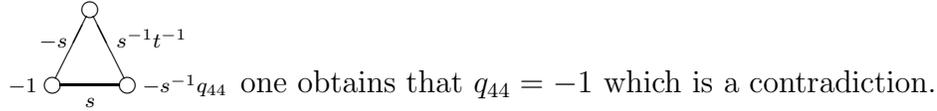
\makebox{
\Dtriangle{}{$-1$}{$t^{-1}$}{$-s^{-1}q_{44}$}{$-s$}{$s^{-1}t^{-1}$}{$s$}
\rule{0pt}{15\unitlength}}~~~~~~
one obtains that $q_{44}=-1$ which is a contradiction.

\textit{Step~5.} Consider the setting $q_{11}=q_{33}=-1$,
$q_{22},q_{44}\not=-1$.
Lemma~\ref{l-specrank3}(vii) applied to
$\roots (\chi ;\Ndbasis _1,\Ndbasis _2,\Ndbasis _3)$ gives that
$(q_{22}r-1)(q_{22}s-1)(q_{22}rs+1)=0$.

Step~5.1. The case $q_{22}r=1$. If $q_{22}s=1$ then the transformations
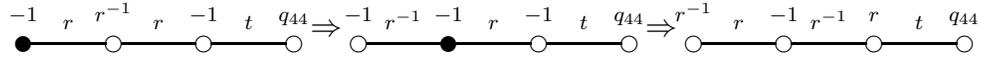
\begin{align*}
\Dchainfour{1}{$-1$}{$r$}{$r^{-1}$}{$r$}{$-1$}{$t$}{$q_{44}$} \Rightarrow
\Dchainfour{2}{$-1$}{$r^{-1}$}{$-1$}{$r$}{$-1$}{$t$}{$q_{44}$} \Rightarrow
\Dchainfour{}{$r^{-1}$}{$r$}{$-1$}{$r^{-1}$}{$r$}{$t$}{$q_{44}$}
\end{align*}
show that $(\roots, \chi ,E)$ is Weyl equivalent to an arithmetic
root system in Step~3. Therefore assume that $q_{22}s\not=1$. In this case
\cite[Lemma\,16]{a-Heck05a} for $\roots (\chi ;\Ndbasis _2,\Ndbasis _3,
\Ndbasis _4)$ implies that $q_{44}t=1$.
Moreover, the transformations
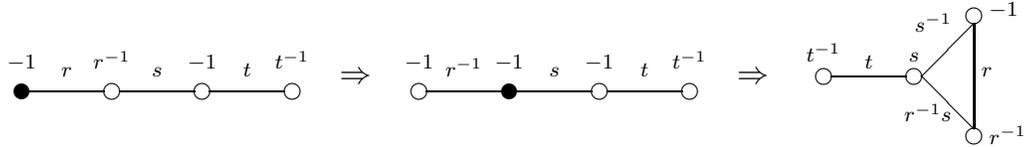
\begin{align*}
\Dchainfour{1}{$-1$}{$r$}{$r^{-1}$}{$s$}{$-1$}{$t$}{$t^{-1}$}
\quad \Rightarrow \quad
\Dchainfour{2}{$-1$}{$r^{-1}$}{$-1$}{$s$}{$-1$}{$t$}{$t^{-1}$}
\quad \Rightarrow \quad
\Drightofway{}{$t^{-1}$}{$t$}{$s$}{$s^{-1}$}{$r^{-1}s$}{$-1$}{$r$}{$r^{-1}$}
\end{align*}
and Lemma~\ref{l-rightofway} show that either equations
$s=-1$, $-r^{-1}t=1$, or equation $st=1$ holds.
If $s=-1$ and $t=-r$ then Lemma~\ref{l-specrank3}(vii) for
$\roots (\chi ;\Ndbasis _1+\Ndbasis _2,
\Ndbasis _2+\Ndbasis _3,\Ndbasis _4)$
\Dchainthree{}{$-1$}{$-r^{-1}$}{$r^{-1}$}{$-r$}{$-r^{-1}$}
and relation $r^2=q_{22}^{-2}\not=1$ imply that $r\in R_4$, and hence
$\gDd _{\chi ,E}$ appears in row~22 of Table~\ref{a-r4}.
Otherwise one has $st=1$
and then Lemma~\ref{l-specrank3}(vii) for
$\roots (\chi ;\Ndbasis _1,\Ndbasis _2+\Ndbasis _3,\Ndbasis _4)$
\Dchainthree{}{$-1$}{$r$}{$-r^{-1}s$}{$~~s^{-1}$}{$s$}
tells that relations $s\in R_3\cup R_6$ and $r=s^2$ hold. Moreover
Lemma~\ref{l-specrank3}(vii) for
$\roots (\chi ;\Ndbasis _1+\Ndbasis _2,\Ndbasis _2+\Ndbasis _3,\Ndbasis _4)$
\Dchainthree{}{$-1~~$}{$r^{-1}s~~~$}{$-r^{-1}s$}{$~~s^{-1}$}{$s$}
implies that $s\in R_3$, and then $\gDd _{\chi ,E}$ can be found in
row~21 of Table \ref{a-r4}.

Step~5.2. The case $q_{22}r\not=1$, $q_{22}s=1$, $s^2,q_{44}^2\not=1$.
By \cite[Lemma\,16]{a-Heck05a} for
$\roots (\chi ;\Ndbasis _1+\Ndbasis _2,\Ndbasis _3,\Ndbasis _4)$
\Dchainthree{}{$-rs^{-1}$}{$s$}{$-1$}{$t$}{$q_{44}$}
and for
$\roots (\chi ;\Ndbasis _1+\Ndbasis _2,\Ndbasis _2+\Ndbasis _3,\Ndbasis _4)$
\Dchainthree{}{$-rs^{-1}$}{$rs^{-1}$}{$-1$}{$t$}{$q_{44}$}
one obtains that either $q_{44}t=1$ or $r=-1$, $s\in R_4$. If $q_{44}t\not=1$
then Lemma~\ref{l-specrank3}(ii) for
$\roots (\chi ;\Ndbasis _1+\Ndbasis _2,\Ndbasis _2+\Ndbasis _3,
\Ndbasis _3+\Ndbasis _4)$
\makebox{
\Dtriangle{}{$-s$}{$-1$}{$-q_{44}t$}{$s$}{$st$}{$s$}
\rule{0pt}{15\unitlength}}~~~
implies that $st=1$.
Therefore Lemma~\ref{l-specrank3}(vi) for
$\roots (\chi ;\Ndbasis _1,\Ndbasis _2+\Ndbasis _3,\Ndbasis _4)$
\Dchainthree{}{$-1$}{$-1$}{$-1$}{$s^{-1}$}{$q_{44}$}
gives a contradiction.

Thus it remains to consider the case $q_{22}s=q_{44}t=1$,
$q_{22}r,s^2,t^2\not=1$. By Lemma~\ref{l-specrank3}(ii) for
$\roots (\chi ;\Ndbasis _1+\Ndbasis _2,\Ndbasis _2+\Ndbasis _3,
\Ndbasis _3+\Ndbasis _4)$
\makebox{
\Dtriangle{}{$-rs^{-1}$}{$-1$}{$-1$}{$rs^{-1}$}{$st$}{$s$}
\rule{0pt}{15\unitlength}}~~~
one gets equation $(st-1)(rst-1)=0$.

Step~5.2.1. Assume that $q_{22}s=q_{44}t=st=1$,
$r\not=s$, $s^2\not=1$. By Lemma~\ref{l-specrank3}(viii) for
$\roots (\chi ;\Ndbasis _1,\Ndbasis _2,\Ndbasis _3)$
one gets that either $r=-1$, $s\in R_3\cup R_4\cup R_6$,
or $r^2=s^{-2}\in R_3$. If $s\in R_6$ and $r^2=s^{-2}$
then Lemma~\ref{l-specrank3}(vi) for
$\roots (\chi ;\Ndbasis _1,\Ndbasis _2+\Ndbasis _3,\Ndbasis _4)$
\Dchainthree{}{$-1$}{$r$}{$-1$}{$s^{-1}$}{$s$}
gives a contradiction. Otherwise one has either $s\in R_3$, $r^2=s$,
or $r=-1$, $s\in R_3\cup R_4\cup R_6$.
Then Lemma~\ref{l-specrank3}(iii) for
$\roots (\chi ;\Ndbasis _4,\Ndbasis _2+\Ndbasis _3,\Ndbasis _1+\Ndbasis _2)$
\Dchainthree{}{$s$}{$s^{-1}$}{$-1$}{$rs^{-1}~~$}{$~~-rs^{-1}$}~~
implies that $r=-1$ and $s\in R_4$. Then $\gDd _{\chi ,E}$
appears in row~12 of Table~\ref{a-r4}.

Step~5.2.2. Assume that $q_{22}s=q_{44}t=rst=1$, and
$r\not=s$, $s^2\not=1$, $r^2s^2\not=1$. If $r\not=-1$ then the transformations
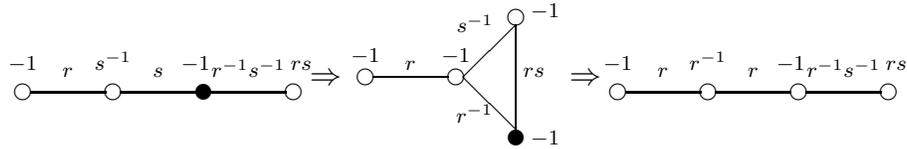
\begin{align*}
\Dchainfour{3}{$-1$}{$r$}{$s^{-1}$}{$s$}{$-1~~$}{$r^{-1}s^{-1}$}{$~~rs$}
\Rightarrow
\Drightofway{b}{$-1$}{$r$}{$-1$}{$s^{-1}$}{$r^{-1}$}{$-1$}{$rs$}{$-1$}
\Rightarrow
\Dchainfour{}{$-1$}{$r$}{$r^{-1}$}{$r$}{$-1~~$}{$r^{-1}s^{-1}$}{$~~rs$}
\end{align*}
imply that $(\roots ,\chi ,E)$ is Weyl equivalent to an arithmetic root
system from Step~5.1. Therefore one can assume that $r=-1$.
By Lemma~\ref{l-specrank3}(viii) for
$\roots (\chi ;\Ndbasis _3,\Ndbasis _2,\Ndbasis _1)$
\Dchainthree{}{$-1$}{$s$}{$s^{-1}$}{$-1$}{$-1$}
one gets relation $s\in R_3\cup R_4\cup R_6$. If $s\in R_4$
then $\gDd _{\chi ,E}$ can be found in row~9 of Table~\ref{a-r4}. Otherwise
Lemma~\ref{l-specrank3}(iii) for
$\roots (\chi ;\Ndbasis _4,\Ndbasis _2+\Ndbasis _3,\Ndbasis _1+\Ndbasis _2)$
\Dchainthree{}{$-s$}{$-s^{-1}$}{$-1$}{$-s^{-1}$}{$s^{-1}$}
gives a contradiction.

Step~5.3. Assume now that $q_{22}r\not=1\not=q_{22}s$ and $q_{22}rs=-1$.
Then the transformation
\begin{align*}
\Dchainfour{1}{$-1$}{$r$}{$q_{22}$}{$s$}{$-1$}{$t$}{$q_{44}$}
\qquad \Rightarrow \qquad
\Dchainfour{}{$-1$}{$r^{-1}$}{$-q_{22}r$}{$s$}{$-1$}{$t$}{$q_{44}$}
\end{align*}
shows that $(\roots, \chi ,E)$ is Weyl equivalent to an arithmetic
root system in Step~5.2.

\textit{Step~6.} $q_{11}=q_{44}=-1$, $q_{22},q_{33}\not=-1$.
By Lemma~\ref{l-specrank3}(vii) for
$\roots (\Ndbasis _1,\Ndbasis _2,\Ndbasis _3)$
one gets $(q_{22}r-1)(q_{22}s-1)(q_{22}rs+1)=0$. The transformation
\begin{align*}
\Dchainfour{1}{$-1$}{$r$}{$q_{22}$}{$s$}{$q_{33}$}{$t$}{$-1$}
\qquad \Rightarrow \qquad
\Dchainfour{}{$-1$}{$r^{-1}~$}{$-q_{22}r$}{$~s$}{$q_{33}$}{$t$}{$-1$}
\end{align*}
shows that one can assume that $(q_{22}r-1)(q_{22}s-1)=0$. By symmetry
one can suppose that equation $(q_{33}s-1)(q_{33}t-1)=0$ holds, too.

Step~6.1. Suppose first that $q_{22}r=q_{33}t=1$ and $r^2,t^2\not=1$.
If $q_{22}s=q_{33}s=1$ then $\gDd _{\chi ,E}$ appears in row~10 of
Table~\ref{a-r4}.
If $q_{22}s=1$ and $q_{33}s\not=1$ then the transformation
\begin{align*}
\Dchainfour{1}{$-1$}{$r$}{$r^{-1}$}{$r$}{$t^{-1}$}{$t$}{$-1$} \Rightarrow
\Dchainfour{2}{$-1$}{$r^{-1}$}{$-1$}{$r$}{$t^{-1}$}{$t$}{$-1$} \Rightarrow
\Dchainfour{}{$r^{-1}$}{$r$}{$-1$}{$r^{-1}$}{$-rt^{-1}$}{$t$}{$-1$}
\end{align*}
shows that $(\roots ,\chi ,E)$ is Weyl equivalent to an arithmetic root system
in Step~3. Finally, if $q_{22}s\not=1$ and $q_{33}s\not=1$ then one obtains
a contradiction to \cite[Corollary 13]{a-Heck05a} for
$\roots (\chi ;\Ndbasis _2,\Ndbasis _3)$.

Step~6.2. Assume now that $q_{22}r=q_{33}s=1$ and
$s\not=t$, $r^2,s^2\not=1$.
Then Lemma~\ref{l-specrank3}(ii) for
$\roots (\chi ;
\Ndbasis _1+\Ndbasis _2,\Ndbasis _2+\Ndbasis _3,
\Ndbasis _3+\Ndbasis _4)$
\makebox{
\Dtriangle{}{$-1$}{$r^{-1}$}{$-s^{-1}t$}{$r^{-1}s$}{$s^{-1}t$}{$s$}
\rule{0pt}{15\unitlength}}~~~
implies that $r=s$.
Therefore if $t\not=-1$ then by the transformations
\begin{align*}
\Dchainfour{1}{$-1$}{$r$}{$r^{-1}$}{$r$}{$r^{-1}$}{$t$}{$-1$} &
\qquad \Rightarrow \qquad
\Dchainfour{2}{$-1$}{$r^{-1}$}{$-1$}{$r$}{$r^{-1}$}{$t$}{$-1$}\\
\Rightarrow \qquad
\Dchainfour{4}{$r^{-1}$}{$r$}{$-1$}{$r^{-1}$}{$-1$}{$t$}{$-1$} &
\qquad \Rightarrow \qquad
\Dchainfour{}{$r^{-1}$}{$r$}{$-1$}{$r^{-1}$}{$t$}{$t^{-1}$}{$-1$}
\end{align*}
one can see that $(\roots ,\chi ,E)$
is Weyl equivalent to an arithmetic
root system from Step~5. Thus suppose that $t=-1$. If $r\in R_4$ then
$\gDd _{\chi ,E}$ appears in row~8 of Table~\ref{a-r4}.
Otherwise the transformations
\begin{align*}
\Dchainfour{3}{$r^{-1}$}{$r$}{$-1$}{$r^{-1}$}{$-1$}{$-1$}{$-1$}
\quad \Rightarrow \quad
\Drightofway{b}{$r^{-1}$}{$r$}{$r^{-1}$}{$r$}{$-r^{-1}$}{$-1$}{$-1$}{$-1$}
\quad \Rightarrow \quad
\Dchainfour{}{$r^{-1}$}{$r$}{$r^{-2}$}{$-r$}{$-1$}{$-1$}{$-1$}
\end{align*}
show that $(\roots ,\chi ,E)$ is Weyl equivalent to an arithmetic
root system from Step~4.

Step~6.3. Suppose now that $q_{22}s=q_{33}s=1$ and
$q_{22}r,q_{33}t\not=1$,
$s^2\not=1$. Then Lemma~\ref{l-specrank3}(ii) for
$\roots (\chi ;
\Ndbasis _1+\Ndbasis _2,\Ndbasis _2+\Ndbasis _3,\Ndbasis _3+\Ndbasis _4)$~
\makebox{
\Dtriangle{}{$-rs^{-1}$}{$s^{-1}$}{$-s^{-1}t$}{$rs^{-1}$}{$s^{-1}t$}{$s$}
\rule{0pt}{15\unitlength}}~~~
gives a contradiction.

\textit{Step~7.} $q_{22}=q_{33}=-1$, $q_{11},q_{44}\not=-1$.
If $s\not=-1$ then the transformation
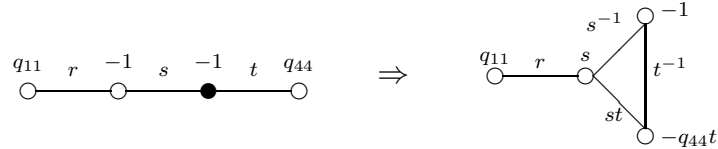
\begin{align*}
\Dchainfour{3}{$q_{11}$}{$r$}{$-1$}{$s$}{$-1$}{$t$}{$q_{44}$}
\qquad \Rightarrow \qquad
\Drightofway{}{$q_{11}$}{$r$}{$s$}{$s^{-1}$}{$st$}{$-1$}{$t^{-1}$}{$-q_{44}t$}
\end{align*}
and Lemma~\ref{l-rightofway} give first that either $st=1$ or $rs=1$.
By twist equivalence one can assume that $st=1$. Then the above transformation
shows also that $(\roots ,\chi ,E)$ is Weyl equivalent
to an arithmetic root system from Step~3 or Step~4.

Assume now that $s=-1$. By \cite[Lemma\,16]{a-Heck05a} for
$\roots (\chi ;\Ndbasis _1,\Ndbasis _2+\Ndbasis _3,\Ndbasis _4)$
and twist equivalence one can assume that $q_{44}t=1$, and hence
$t^2\not=1$. Therefore the transformation
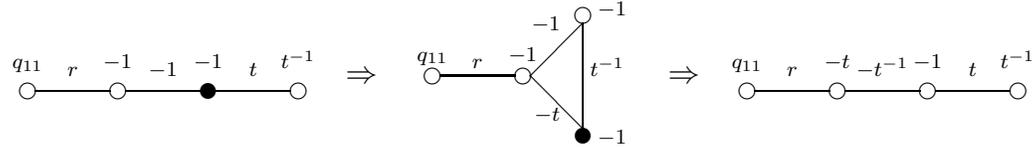
\begin{align*}
\Dchainfour{3}{$q_{11}$}{$r$}{$-1$}{$-1$}{$-1$}{$t$}{$t^{-1}$}
\quad \Rightarrow \quad
\Drightofway{b}{$q_{11}$}{$r$}{$-1$}{$-1$}{$-t$}{$-1$}{$t^{-1}$}{$-1$}
\quad \Rightarrow \quad
\Dchainfour{}{$q_{11}$}{$r$}{$-t$}{$-t^{-1}$}{$-1$}{$t$}{$t^{-1}$}
\end{align*}
tells that $(\roots ,\chi ,E)$ is Weyl equivalent to an arithmetic
root system in Step~3.

\textit{Step~8.} $q_{11}=q_{22}=q_{33}=-1$, $q_{44}\not=-1$.
Like at the beginning of Step~4 one can conclude that $r=-1$.

If $s\not=-1$ then Lemma~\ref{l-rightofway} and the transformation
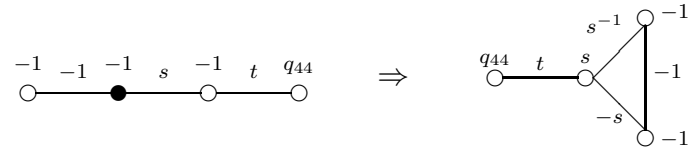
\begin{align*}
\Dchainfour{2}{$-1$}{$-1$}{$-1$}{$s$}{$-1$}{$t$}{$q_{44}$}
\qquad \Rightarrow \qquad
\Drightofway{}{$q_{44}$}{$t$}{$s$}{$s^{-1}$}{$-s$}{$-1$}{$-1$}{$-1$}
\end{align*}
show that $st=1$. In this case, if $q_{44}t=1$ or $q_{44}t\not=1$
then the transformations
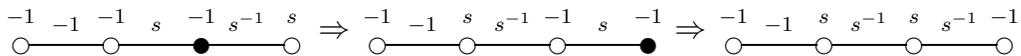
\begin{align*}
\Dchainfour{3}{$-1$}{$-1$}{$-1$}{$s$}{$-1$}{$s^{-1}$}{$s$}
\,\, \Rightarrow \,\,
\Dchainfour{4}{$-1$}{$-1$}{$s$}{$s^{-1}$}{$-1$}{$s$}{$-1$}
\,\, \Rightarrow \,\,
\Dchainfour{}{$-1$}{$-1$}{$s$}{$s^{-1}$}{$s$}{$s^{-1}$}{$-1$}
\end{align*}
and
\begin{align*}
\Dchainfour{3}{$-1$}{$-1$}{$-1$}{$s$}{$-1$}{$s^{-1}$}{$q_{44}$}
\qquad \Rightarrow \qquad
\Dchainfour{}{$-1$}{$-1$}{$s$}{$s^{-1}$}{$-1$}{$s$}{$-s^{-1}q_{44}$},
\end{align*}
give an arithmetic root system from Step~6 and Step~5, respectively.

Consider now the case when $s=-1$. If $q_{44}t=1$ then $t^2\not=1$ and
hence the transformations
\begin{align*}
\Dchainfour{3}{$-1$}{$-1$}{$-1$}{$-1$}{$-1$}{$t$}{$t^{-1}$}
\quad \Rightarrow \quad
\Drightofway{b}{$-1$}{$-1$}{$-1$}{$-1$}{$-t$}{$-1$}{$t^{-1}$}{$-1$}
\quad \Rightarrow \quad
\Dchainfour{}{$-1$}{$-1$}{$-t$}{$-t^{-1}$}{$-1$}{$t$}{$t^{-1}$}
\end{align*}
show that $(\roots ,\chi ,E)$ is Weyl equivalent to an arithmetic root system
in Step~5. On the other hand, if $q_{44}t\not=1$ then
Lemma~\ref{l-specrank3}(vi) for
$\roots (\chi ;\Ndbasis _2,\Ndbasis _3,\Ndbasis _4)$
\Dchainthree{}{$-1$}{$-1$}{$-1$}{$t$}{$q_{44}$}
implies that either $q_{44}\in R_3$, $t^2q_{44}=1$, or $t=-1$, $q_{44}\in R_4$.
In the second case $\gDd _{\chi ,E}$ appears in row~2 of Table~\ref{a-r4},
and in the
first case, where $t\in R_3\cup R_6$, Lemma~\ref{l-specrank3}(vi) for
$\roots (\chi ;\Ndbasis _1,\Ndbasis _2,\Ndbasis _3+\Ndbasis _4)$
\Dchainthree{}{$-1$}{$-1$}{$-1$}{$-1$}{$-t^{-1}$}
gives a contradiction.

\textit{Step~9.} $q_{11}=q_{22}=q_{44}=-1$, $q_{33}\not=-1$.
By the same arguments as in Steps~4 and 5.3 one can assume that
equations $r=-1$, $(q_{33}s-1)(q_{33}t-1)=0$ hold.

Suppose that $r=-1$, $q_{33}s=1$, and $s^2\not=1$.
As in Step~4.1 it suffices to consider the case $t=s$.
Then the transformations
\begin{align*}
\Dchainfour{4}{$-1$}{$-1$}{$-1$}{$s$}{$s^{-1}$}{$s$}{$-1$}
\,\, \Rightarrow \,\, 
\Dchainfour{3}{$-1$}{$-1$}{$-1$}{$s$}{$-1$}{$s^{-1}$}{$-1$}
\,\,  \Rightarrow \,\, 
\Dchainfour{}{$-1$}{$-1$}{$s$}{$s^{-1}$}{$-1$}{$s$}{$s^{-1}$}
\end{align*}
give an arithmetic root system from Step~5.

Assume that $r=-1$, $q_{33}t=1$, and $s\not=t$, $t^2\not=1$.
Then Lemma~\ref{l-rightofway} and the transformations
\begin{align*}
\Dchainfour{4}{$-1$}{$-1$}{$-1$}{$s$}{$t^{-1}$}{$t$}{$-1$}
\quad \Rightarrow \quad
\Dchainfour{3}{$-1$}{$-1$}{$-1$}{$s$}{$-1$}{$t^{-1}$}{$-1$}
\quad \Rightarrow \quad
\Drightofway{}{$-1$}{$-1$}{$s$}{$s^{-1}$}{$st^{-1}$}{$-1$}{$t$}{$t^{-1}$}
\end{align*}
give that $s=-1$. Moreover, by Lemma~\ref{l-specrank3}(vi) for
the arithmetic root system
$\roots (\chi ;\Ndbasis _1,\Ndbasis _2,\Ndbasis _3)$
\Dchainthree{}{$-1$}{$-1$}{$-1$}{$-1$}{$t^{-1}$}
one has $t\in R_4$, and hence $\gDd _{\chi ,E}$ appears in row~9
of Table~\ref{a-r4}.

\textit{Step~10.} $q_{11}=q_{22}=q_{33}=q_{44}=-1$.
Following the argument at the beginning of Step~4 one can
assume that $r=t=-1$. Then Lemma~\ref{l-specrank3}(vii) for
$\roots (\chi ;\Ndbasis _1,\Ndbasis _2+\Ndbasis _3,\Ndbasis _4)$
\Dchainthree{}{$-1$}{$-1$}{$s$}{$-1$}{$-1$}
tells that $s=-1$ and hence $\gDd _{\chi ,E}$ appears in row~1
of Table~\ref{a-r4}.
\end{bew}

\section{Simple chains}
\label{sec-simplechains}

The results of the previous section show that labeled path graphs
play an important role in the study of generalized Dynkin diagrams.
As a preparation for Theorem~\ref{t-classrank>4} here
a special class of such graphs will be studied.

Let $E$ be a basis of $\ndZ ^d$, where $d\ge 2$,
and let $\chi $ be a bicharacter on $\ndZ ^d$.

\begin{defin}
Assume that $(\roots ,\chi ,E)$ is an arithmetic root system of rank $d\ge 2$
and $\gDd _{\chi ,E}$ is a labeled path graph. Call this graph
a \textit{simple chain} (of length $d$) if
\begin{gather*}
(q_{11}q_{12}q_{21}-1)(q_{11}+1)=0,\quad
(q_{dd}q_{d,d-1}q_{d-1,d}-1)(q_{dd}+1)=0,\\
q_{ii}^2q_{i-1,i}q_{i,i-1}q_{i,i+1}q_{i+1,i}=1\text{ for $1<i<d$.}
\end{gather*}
By Lemma~\ref{l-specrank3}(vii) the latter equations hold
if and only if
\begin{align}\label{eq-simplechaincond}
q_{ii}+1=q_{i-1,i}q_{i,i-1}q_{i,i+1}q_{i+1,i}-1=0 \text{ or }
q_{ii}q_{i-1,i}q_{i,i-1}=q_{ii}q_{i,i+1}q_{i+1,i}=1.
\end{align}
If $(\roots _i,\chi _i,E_i)$, where $i\in \{1,2\}$, are arithmetic root
systems, the generalized Dynkin diagrams
$\gDd _{\chi _i,E_i}$ are labeled path graphs,
and $E_1=\{\Ndbasis _1,\ldots ,\Ndbasis _n\}\subset E_2
=\{\Ndbasis _1,\ldots ,\Ndbasis _d\}$, where $2\le n<d$,
then $\gDd _{\chi _2,E_2}$ is called a
\textit{prolongation of $\gDd _{\chi _1,E_1}$
to the right by a simple chain (of length $d-n$)}
if
\begin{align*}
(q_{dd}q_{d,d-1}q_{d-1,d}-1)(q_{dd}+1)=&0\text{ and }\\
q_{jj}^2q_{j-1,j}q_{j,j-1}q_{j,j+1}q_{j+1,j}=&1\text{ for $n\le j<d$.}
\end{align*}
Similarly one defines prolongation of a labeled path graph
to the left by a simple chain. Because no other prolongations
will be considered, usually the attribute ``by a simple chain'' will be
omitted.
\end{defin}

\begin{lemma}
If $\gDd _{\chi ,E}$ is a simple chain and $\gDd _{\chi ',E}$
is Weyl equivalent to $\gDd _{\chi ,E}$ then $\gDd _{\chi ',E}$ is a simple
chain.
\end{lemma}

\begin{bew}
If $i\in \{1,\ldots ,d\}$ and $q_{ii}\not=-1$ then
Equation~(\ref{eq-simplechaincond})
for $\gDd _{\chi ,E}$ implies that
$\gDd _{\chi ,E}=\gDd _{\chi ',E}$, where $\chi '=\chi \circ (s_{\Ndbasis _i,E}
\times s_{\Ndbasis _i,E})$. Otherwise application of $s_{\Ndbasis _i,E}$
gives
\begin{gather*}
\rule[-4\unitlength]{0pt}{5\unitlength}
\begin{picture}(70,10)(0,3)
\put(1,1){\framebox(11,8){}}
\put(12,3){\line(1,0){10}}
\put(23,3){\circle{2}}
\put(24,3){\line(1,0){10}}
\put(35,3){\circle{2}}
\put(36,3){\line(1,0){10}}
\put(47,3){\circle{2}}
\put(48,3){\line(1,0){10}}
\put(58,1){\framebox(11,8){}}
\put(17,5){\makebox[0pt]{\scriptsize $p$}}
\put(22,6){\makebox[0pt]{\scriptsize $r$}}
\put(29,5){\makebox[0pt]{\scriptsize $s$}}
\put(36,6){\makebox[0pt]{\scriptsize $-1$}}
\put(41,5){\makebox[0pt]{\scriptsize $s^{-1}$}}
\put(47,6){\makebox[0pt]{\scriptsize $t$}}
\put(53,5){\makebox[0pt]{\scriptsize $u$}}
\end{picture}\quad \Rightarrow \\
\rule[-4\unitlength]{0pt}{5\unitlength}
\begin{picture}(70,10)(0,3)
\put(1,1){\framebox(11,8){}}
\put(12,3){\line(1,0){10}}
\put(23,3){\circle{2}}
\put(24,3){\line(1,0){10}}
\put(35,3){\circle{2}}
\put(36,3){\line(1,0){10}}
\put(47,3){\circle{2}}
\put(48,3){\line(1,0){10}}
\put(58,1){\framebox(11,8){}}
\put(17,5){\makebox[0pt]{\scriptsize $p$}}
\put(22,6){\makebox[0pt]{\scriptsize $-rs$}}
\put(29,5){\makebox[0pt]{\scriptsize $s^{-1}$}}
\put(36,6){\makebox[0pt]{\scriptsize $-1$}}
\put(41,5){\makebox[0pt]{\scriptsize $s$}}
\put(47,6){\makebox[0pt]{\scriptsize $-s^{-1}t$}}
\put(53,5){\makebox[0pt]{\scriptsize $u$}}
\end{picture}.
\end{gather*}
Thus if $(rs-1)(r+1)=0$ then $(-rs+1)((-rs)s^{-1}-1)=0$, and if
$(pr-1)(r+1)=0$ and $s=1/pr^2$ then
$$(p(-rs)-1)(-rs+1)=(-pr/pr^2-1)(-r/pr^2+1)=0.$$
Thus the claim follows by symmetry.
\end{bew}

Let $\gDd _{\chi ,E}$ be a simple chain and set $q:=q_{dd}^2q_{d-1,d}
q_{d,d-1}$. By (\ref{eq-simplechaincond}) this means that
if $q_{dd}\not=-1$ then $q=q_{dd}$ and if $q_{dd}=-1$ then
$q=q_{d-1,d}q_{d,d-1}$.
Equations~(\ref{eq-simplechaincond})
imply that for all $i\in \{1,\ldots ,d-1\}$
one has $q_{i,i+1}q_{i+1,1}\in \{q,q^{-1}\}$ and $q_{11}^2q_{12}q_{21}\in
\{q,q^{-1}\}$. Moreover, the knowledge of $q$ and all indices $i$
with $1\le i\le d$ and $q_{i-1,i}q_{i,i-1}=q$, where $q_{01}q_{10}:=
1/q_{11}^2q_{12}q_{21}$, determines $\gDd _{\chi ,E}$ uniquely.
Therefore the symbol
\begin{align}\label{eq-simplechaingraph}
\begin{picture}(26,8)
\put(13,3){\oval(26,6)}
\put(0,0){\makebox(26,6){\scriptsize $C(d,q;i_1,\ldots ,i_j)$}}
\end{picture}
\end{align}
will be used for the simple chain of length $d$ for which
$q=q_{d-1,d}q_{d,d-1}q_{dd}^2$ holds, and for which equation
$q_{i-1,i}q_{i,i-1}=q$, where $1\le i\le d$,
is valid if and only if $i\in \{i_1,i_2,\ldots ,i_j\}$.
Additionally the convention $1\le i_1<i_2<\ldots <i_j\le d$ is fixed, where
$0\le j\le d$.
For brevity, in the running text also the notation $C(d,q;i_1,\ldots ,i_j)$
will be used.

\begin{bsp}
For $q\in k^\ast \setminus \{1\}$ the generalized Dynkin diagrams
\begin{align*}
&\Dchainfour{}{$-1$}{$q^{-1}$}{$-1$}{$q$}{$q^{-1}$}{$q$}{$-1$}
&
%
%
&
\rule[-4\unitlength]{0pt}{5\unitlength}
\begin{picture}(50,4)(0,3)
\put(1,1){\circle{2}}
\put(2,1){\line(1,0){10}}
\put(13,1){\circle{2}}
\put(14,1){\line(1,0){10}}
\put(25,1){\circle{2}}
\put(26,1){\line(1,0){10}}
\put(37,1){\circle{2}}
\put(38,1){\line(1,0){10}}
\put(49,1){\circle{2}}
\put(1,4){\makebox[0pt][r]{\scriptsize $-1$}}
\put(7,3){\makebox[0pt]{\scriptsize $q^{-1}$}}
\put(13,4){\makebox[0pt]{\scriptsize $-1$}}
\put(19,3){\makebox[0pt]{\scriptsize $q$}}
\put(25,4){\makebox[0pt]{\scriptsize $q^{-1}$}}
\put(31,3){\makebox[0pt]{\scriptsize $q$}}
\put(37,4){\makebox[0pt]{\scriptsize $-1$}}
\put(43,3){\makebox[0pt]{\scriptsize $q^{-1}$}}
\put(49,4){\makebox[0pt]{\scriptsize $q$}}
\end{picture}
\end{align*}
are simple chains of length 4 and 5, respectively.
One has $q_{i-1,i}q_{i,i-1}=q$ if and only if
$i\in \{1,3,4\}$, and hence for these simple chains the symbols
\begin{align*}
&\begin{picture}(26,8)
\put(13,3){\oval(26,6)}
\put(0,0){\makebox(26,6){\scriptsize $C(4,q;1,3,4)$}}
\end{picture}
&
&\begin{picture}(26,8)
\put(13,3){\oval(26,6)}
\put(0,0){\makebox(26,6){\scriptsize $C(5,q;1,3,4)$}}
\end{picture}
\end{align*}
will be used.
\end{bsp}

\begin{satz}\label{s-Weqsimplechains}
Two simple chains $C(d,q;i_1,\ldots ,i_j)$, where $0\le j\le d$,
and $C(d',q';i'_1,i'_2,\ldots ,i'_{j'})$,  where $0\le j'\le d'$,
are Weyl equivalent if and only if $d=d'$ and one of the following conditions
hold:
\begin{align*}
&(i) \quad q=q',\quad j=j', &
&(ii) \quad qq'=1,\quad j+j'=d+1.
\end{align*}
\end{satz}

\begin{bew}
If $q_{i-1,i}q_{i,i-1}\not=q_{i,i+1}q_{i+1,i}$ for some $1\le i<d$ then
Equation~(\ref{eq-simplechaincond}) implies that
$q_{ii}=-1$ and exactly
one of $q_{i-1,i}q_{i,i-1}$, $q_{i,i+1}q_{i+1,i}$ is equal to $q$.
Hence one has either $i\in \{i_1,\ldots ,i_j\}$ or
$i+1\in \{i_1,\ldots ,i_j\}$, but both relations can not be fulfilled.
In this case applying the reflection $s_{\Ndbasis _i,E}$ is equivalent to
replacing the index $i$ respectively $i+1$ in the argument of $C$
by $i+1$ respectively $i$.
Therefore $(i)$ implies that the given simple chains are Weyl equivalent.
Further, if $i_j=d$ then $q_{d-1,d}q_{d,d-1}=q$ and $q_{dd}=-1$.
Application of $s_{\Ndbasis _d,E}$ transforms $q$ to $q^{-1}$
and leaves all $q_{i-1,i}q_{i,i-1}$ with $1\le i\le d-1$ invariant.
Thus the simple chain
$C(d,q;i_1,\ldots ,i_{j-1},d)$ is Weyl equivalent to
$C(d,q^{-1};i'_1,\ldots ,i'_{d-j},d)$,
where relation $\{i_1,\ldots ,i_{j-1}\}
\cup \{i'_1,\ldots ,i'_{d-j}\}=\{1,2,\ldots ,d-1\}$ holds.
Hence using the first part of the proof one concludes that
if condition (ii) holds then the two chains
$C(d,q;i_1,\ldots ,i_j)$
and $C(d,q';i'_1,i'_2,\ldots ,i'_{j'})$ are Weyl equivalent.

On the other hand, if $q_{i-1,i}q_{i,i-1}=q_{i,i+1}q_{i+1,i}$
for some $1\le i<d$ then
the simple chain is invariant under the reflection $s_{\Ndbasis _i,E}$.
Similarly, if $i_j\not=d$ then $q_{d-1,d}q_{d,d-1}q_{dd}=1$ and
the simple chain is invariant under the reflection $s_{\Ndbasis _d,E}$.
Therefore the claim of the proposition follows
from the first part of the proof.
\end{bew}

\section{Connected arithmetic root systems
of rank higher than four}
\label{sec-classrank>4}

Having once all arithmetic root systems of rank 4,
it is relatively easy to
classify those of higher rank too.
The methods used are the same as in
Section~\ref{sec-classrank4}. The aim of this section
is to prove the following theorem.

\begin{thm}\label{t-classrank>4}
Let $k$ be a field of characteristic 0.
Then twist equivalence classes of connected arithmetic root systems
of rank bigger than 4
are in one-to-one correspondence to generalized Dynkin diagrams
appearing in Table~\ref{a-rb4}.
Moreover, two such arithmetic root systems are
Weyl equivalent if and only if their generalized Dynkin diagrams appear
in the same row of Table~\ref{a-rb4} and can be presented with the same
fixed parameter.
\end{thm}

Again first several special cases
will be analyzed before the theorem can be proven.

Let $E=\{\Ndbasis _1,\ldots ,\Ndbasis _d\}$
be a basis of $\ndZ ^d$, where $d\ge 5$, and let $\chi $ be a
bicharacter on $\ndZ ^d$ such that $(\roots ,\chi ,E)$ is a connected
arithmetic root system.

Lemmata \ref{l-notetraeder}, \ref{l-nodoubletriangle}, and \ref{l-no4cycle}
can be generalized as follows.

\begin{lemma}\label{l-nobigcycle}
Let $\cG $ be a labeled cycle graph which can be obtained
from $\gDd _{\chi ,E}$ by omitting some vertices and some edges.
Then $\cG $ has three vertices.
\end{lemma}

\begin{bew}
Assume to the contrary that there exists $r\ge 4$ such that
$\prod _{i=1}^{r-1}(q_{i,i+1}q_{i+1,i}-1)(q_{1r}q_{r1}-1)\not=0$.
By \cite[Lemma\,2]{a-Heck05a} and Proposition~\ref{s-rootsubs}
one gets $\sum _{i=3}^{r-1}\Ndbasis _i\in \roots ^+_E$. Thus
by Proposition~\ref{s-rootsubs} the triple $\roots (\chi ;
\Ndbasis _1,\Ndbasis _2,\sum _{i=3}^{r-1}\Ndbasis _i,\Ndbasis _r)$ is
an arithmetic root system. By equations
$\chi \chi \op (\Ndbasis _2,\sum _{i=3}^{r-1}\Ndbasis _i)=q_{23}q_{32}$ and
$\chi \chi \op (\sum _{i=3}^{r-1}\Ndbasis _i,\Ndbasis _r)=q_{r-1,r}q_{r,r-1}$
the generalized Dynkin diagram of the arithmetic root system
$\roots (\chi ;\Ndbasis _1,\Ndbasis _2,\sum _{i=3}^{r-1}\Ndbasis _i,
\Ndbasis _r)$ contains a labeled cycle with four vertices which is a
contradiction to 
Lemmata \ref{l-notetraeder}, \ref{l-nodoubletriangle}, and \ref{l-no4cycle}.
\end{bew}

The following lemmata will be helpful to find $(T,E)\in \extWBG _{\chi ,E}$
such that the generalized Dynkin diagram
$\gDd _{\chi \circ (T\times T),E}$ of the Weyl equivalent arithmetic
root system $(T^{-1}(\roots ),\chi \circ (T\times T), E)$
is of a particularly simple shape.

\begin{lemma}\label{l-tent}
Assume that $\gDd _{\chi ,E}$ is of the form
 \begin{align*}
\rule[-4\unitlength]{0pt}{5\unitlength}
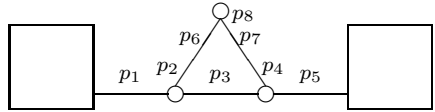
\begin{picture}(58,13)(0,3)
\put(1,1){\framebox(11,11){}}
\put(12,3){\line(1,0){10}}
\put(23,3){\circle{2}}
\put(23,4){\line(2,3){6}}
\put(24,3){\line(1,0){10}}
\put(35,3){\circle{2}}
\put(35,4){\line(-2,3){6}}
\put(29,14){\circle{2}}
\put(36,3){\line(1,0){10}}
\put(46,1){\framebox(11,11){}}
\put(17,5){\makebox[0pt]{\scriptsize $p_1$}}
\put(22,6){\makebox[0pt]{\scriptsize $p_2$}}
\put(29,5){\makebox[0pt]{\scriptsize $p_3$}}
\put(36,6){\makebox[0pt]{\scriptsize $p_4$}}
\put(41,5){\makebox[0pt]{\scriptsize $p_5$}}
\put(25,10){\makebox[0pt]{\scriptsize $p_6$}}
\put(33,10){\makebox[0pt]{\scriptsize $p_7$}}
\put(32,13){\makebox[0pt]{\scriptsize $p_8$}}
\end{picture}
 \end{align*}
where both boxes may contain an arbitrary finite labeled graph.
Then equation $p_8=-1$ holds,
and $(\roots ,\chi ,E)$ is Weyl equivalent to an arithmetic root system
with generalized Dynkin diagram
 \begin{align*}
\rule[-4\unitlength]{0pt}{5\unitlength}
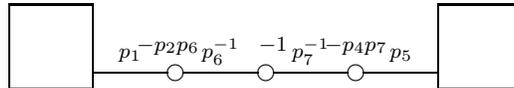
\begin{picture}(70,13)(0,3)
\put(1,1){\framebox(11,11){}}
\put(12,3){\line(1,0){10}}
\put(23,3){\circle{2}}
\put(24,3){\line(1,0){10}}
\put(35,3){\circle{2}}
\put(36,3){\line(1,0){10}}
\put(47,3){\circle{2}}
\put(48,3){\line(1,0){10}}
\put(58,1){\framebox(11,11){}}
\put(17,5){\makebox[0pt]{\scriptsize $p_1$}}
\put(22,6){\makebox[0pt]{\scriptsize $-p_2p_6$}}
\put(29,5){\makebox[0pt]{\scriptsize $p_6^{-1}$}}
\put(36,6){\makebox[0pt]{\scriptsize $-1$}}
\put(41,5){\makebox[0pt]{\scriptsize $p_7^{-1}$}}
\put(47,6){\makebox[0pt]{\scriptsize $-p_4p_7$}}
\put(53,5){\makebox[0pt]{\scriptsize $p_5$}}
\end{picture}
 \end{align*}
\end{lemma}

\begin{bew}
Without loss of generality assume that the vertex $i\in \{1,2,\ldots ,5\}$
of the subgraph 
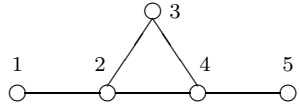
\begin{align*}
\rule[-4\unitlength]{0pt}{5\unitlength}
\begin{picture}(38,11)(0,3)
\put(1,1){\circle{2}}
\put(2,1){\line(1,0){10}}
\put(13,1){\circle{2}}
\put(13,2){\line(2,3){6}}
\put(14,1){\line(1,0){10}}
\put(25,1){\circle{2}}
\put(25,2){\line(-2,3){6}}
\put(19,12){\circle{2}}
\put(26,1){\line(1,0){10}}
\put(37,1){\circle{2}}
\put(1,4){\makebox[0pt]{\scriptsize $1$}}
\put(12,4){\makebox[0pt]{\scriptsize $2$}}
\put(26,4){\makebox[0pt]{\scriptsize $4$}}
\put(37,4){\makebox[0pt]{\scriptsize $5$}}
\put(22,11){\makebox[0pt]{\scriptsize $3$}}
\end{picture}
 \end{align*}
of $\gDd _{\chi ,E}$ corresponds to the basis element
$\Ndbasis _i\in E$. Then
$\roots (\chi ;\Ndbasis _1+\Ndbasis _2,\Ndbasis _2+\Ndbasis _3,
\Ndbasis _3+\Ndbasis _4,\Ndbasis _4+\Ndbasis _5)$
is an arithmetic root system
and by Lemma~\ref{l-specrank3}(ii)
it has the generalized Dynkin diagram
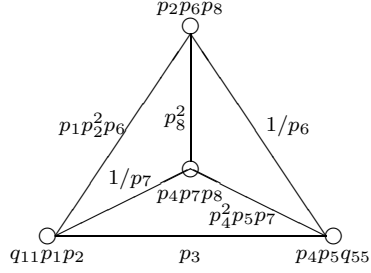
\begin{align*}
\rule[-4\unitlength]{0pt}{5\unitlength}
\begin{picture}(48,33)(0,3)
\put(1,3){\circle{2}}
\put(2,3){\line(1,0){36}}
\put(39,3){\circle{2}}
\put(2,3){\line(2,3){18}}
\put(38,3){\line(-2,3){18}}
\put(20,31){\circle{2}}
\put(20,12){\circle{2}}
\put(2,3){\line(2,1){18}}
\put(38,3){\line(-2,1){18}}
\put(20,13){\line(0,1){17}}
\put(1,0){\makebox[0pt]{\scriptsize $q_{11}p_1p_2$}}
\put(20,0){\makebox[0pt]{\scriptsize $p_3$}}
\put(39,0){\makebox[0pt]{\scriptsize $p_4p_5q_{55}$}}
\put(20,33){\makebox[0pt]{\scriptsize $p_2p_6p_8$}}
\put(7,17){\makebox[0pt]{\scriptsize $p_1p_2^2p_6$}}
\put(12,10){\makebox[0pt]{\scriptsize $1/p_7$}}
\put(18,18){\makebox[0pt]{\scriptsize $p_8^2$}}
\put(20,8){\makebox[0pt]{\scriptsize $p_4p_7p_8$}}
\put(27,5){\makebox[0pt]{\scriptsize $p_4^2p_5p_7$}}
\put(33,17){\makebox[0pt]{\scriptsize $1/p_6$}}
\end{picture}
 \end{align*}
Therefore one has $p_8=-1$. Indeed, since
$p_3,p_6,p_7\not=1$, relation $p_8^2\not=1$ would contradict
to Lemma~\ref{l-nobigcycle}.
Now set $T:=s_{\Ndbasis _3,E}$ and observe that
the arithmetic root system
$(T^{-1}(\roots ),\chi \circ (T\times T),E)$, which
is Weyl equivalent to $(\roots ,\chi ,E)$,
satisfies the claim of the lemma.
\end{bew}

\begin{lemma}\label{l-nooctopus}
For all $d\ge 5$
the graph $\gDd _{\chi ,E}$ does not contain a subgraph of the form
 \begin{align*}
\rule[-4\unitlength]{0pt}{5\unitlength}
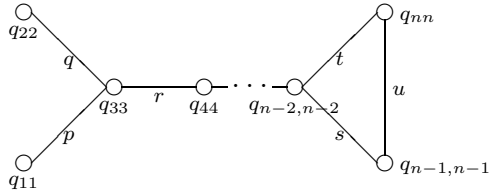
\begin{picture}(53,25)(0,3)
\put(1,3){\circle{2}}
\put(2,3){\line(1,1){10}}
\put(13,13){\circle{2}}
\put(1,23){\circle{2}}
\put(2,23){\line(1,-1){10}}
\put(14,13){\line(1,0){10}}
\put(25,13){\circle{2}}
\put(26,13){\line(1,0){2}}
\put(31,13){\makebox[0pt]{$\ldots $}}
\put(34,13){\line(1,0){2}}
\put(37,13){\circle{2}}
\put(38,13){\line(1,1){10}}
\put(38,13){\line(1,-1){10}}
\put(49,23){\circle{2}}
\put(49,3){\circle{2}}
\put(49,4){\line(0,1){18}}
\put(1,0){\makebox[0pt]{\scriptsize $q_{11}$}}
\put(1,20){\makebox[0pt]{\scriptsize $q_{22}$}}
\put(7,6){\makebox[0pt]{\scriptsize $p$}}
\put(7,16){\makebox[0pt]{\scriptsize $q$}}
\put(13,10){\makebox[0pt]{\scriptsize $q_{33}$}}
\put(19,11){\makebox[0pt]{\scriptsize $r$}}
\put(25,10){\makebox[0pt]{\scriptsize $q_{44}$}}
\put(37,10){\makebox[0pt]{\scriptsize $q_{n-2,n-2}$}}
\put(43,6){\makebox[0pt]{\scriptsize $s$}}
\put(43,16){\makebox[0pt]{\scriptsize $t$}}
\put(51,2){\makebox{\scriptsize $q_{n-1,n-1}$}}
\put(50,12){\makebox{\scriptsize $u$}}
\put(51,22){\makebox{\scriptsize $q_{nn}$}}
\end{picture}
 \end{align*}
where $p,q,r,s,t\not=1$ (but $u$ is allowed to be equal to 1), and
$5\le n\le d$.
\end{lemma}

\begin{bew}
Assume to the contrary that the claim of the lemma is false.
Then $\roots (\chi ;\Ndbasis _1,\Ndbasis _2,\sum _{i=3}^{n-2}\Ndbasis _i,
\Ndbasis _{n-1},\Ndbasis _n)$ is an arithmetic root system with
generalized Dynkin diagram of the same form but only with five vertices,
that is $\Ndbasis _3=\Ndbasis _{n-2}$. Therefore it is sufficient
to prove the lemma for $n=5$.

Set $\aNdbasis _1:=\Ndbasis _1+\Ndbasis _3$, $\aNdbasis _2:=\Ndbasis _2
+\Ndbasis _3$, $\aNdbasis _3:=\Ndbasis _4$ and $\aNdbasis _4:=\Ndbasis _5$.
On the one hand
$\roots (\chi ;\aNdbasis _1,\aNdbasis _2,\aNdbasis _3,\aNdbasis _4)$
is finite. On the other hand one has
$\chi \chi \op (\aNdbasis _1,\aNdbasis _3)=s$,
$\chi \chi \op (\aNdbasis _3,\aNdbasis _2)=s$,
$\chi \chi \op (\aNdbasis _2,\aNdbasis _4)=t$ and
$\chi \chi \op (\aNdbasis _4,\aNdbasis _1)=t$, which is in view of $s\not=1$,
$t\not=1$ a contradiction to Lemmata \ref{l-notetraeder},
\ref{l-nodoubletriangle} and \ref{l-no4cycle}.
\end{bew}

\begin{lemma}\label{l-noufo}
The graph $\gDd _{\chi ,E}$ does not contain a subgraph of the form
 \begin{align*}
\rule[-4\unitlength]{0pt}{5\unitlength}
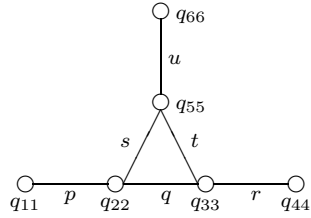
\begin{picture}(38,24)(0,3)
\put(1,3){\circle{2}}
\put(2,3){\line(1,0){10}}
\put(13,3){\circle{2}}
\put(14,3){\line(1,0){10}}
\put(25,3){\circle{2}}
\put(26,3){\line(1,0){10}}
\put(37,3){\circle{2}}
\put(14,3){\line(1,2){5}}
\put(24,3){\line(-1,2){5}}
\put(19,14){\circle{2}}
\put(19,15){\line(0,1){10}}
\put(19,26){\circle{2}}
\put(1,0){\makebox[0pt]{\scriptsize $q_{11}$}}
\put(13,0){\makebox[0pt]{\scriptsize $q_{22}$}}
\put(25,0){\makebox[0pt]{\scriptsize $q_{33}$}}
\put(37,0){\makebox[0pt]{\scriptsize $q_{44}$}}
\put(21,13){{\scriptsize $q_{55}$}}
\put(21,25){{\scriptsize $q_{66}$}}
\put(7,1){\makebox[0pt]{\scriptsize $p$}}
\put(19,1){\makebox{\scriptsize $q$}}
\put(31,1){\makebox{\scriptsize $r$}}
\put(15,8){\makebox[0pt][r]{\scriptsize $s$}}
\put(23,8){\makebox[0pt][l]{\scriptsize $t$}}
\put(20,19){\makebox[0pt][l]{\scriptsize $u$}}
\end{picture}
 \end{align*}
where $p,q,r,s,t,u\not=1$.
\end{lemma}

\begin{bew}
Again one can perform an indirect proof. If $\gDd _{\chi ,E}$ would contain
a subgraph as above then $\roots (\chi ;\Ndbasis _1,\Ndbasis _2,
\Ndbasis _3+\Ndbasis _4,\Ndbasis _3+\Ndbasis _5,
\Ndbasis _5+\Ndbasis _6)$
would be an arithmetic root system. Note that
$\chi \chi \op (\Ndbasis _1,\Ndbasis _3+\Ndbasis _4)=
\chi \chi \op (\Ndbasis _1,\Ndbasis _5+\Ndbasis _6)=
\chi \chi \op (\Ndbasis _1,\Ndbasis _3+\Ndbasis _5)=1$
and
$\chi \chi \op (\Ndbasis _2,\Ndbasis _1)$,
$\chi \chi \op (\Ndbasis _2,\Ndbasis _3+\Ndbasis _4)$,
$\chi \chi \op (\Ndbasis _2,\Ndbasis _5+\Ndbasis _6)$ and
$\chi \chi \op (\Ndbasis _2,\Ndbasis _3+\Ndbasis _5)$ are all
different from 1, the latter because of Lemma~\ref{l-specrank3}(ii).
Depending on the other values of $\chi \chi \op $
this gives a contradiction
to one of Lemmata~\ref{l-notetraeder}, \ref{l-nodoubletriangle},
or \ref{l-nooctopus}.
\end{bew}

\begin{lemma}\label{l-Dngraph}
Assume that $\gDd _{\chi ,E}$ contains a subgraph of the form
 \begin{align}\label{eq-gDdDn}
\rule[-4\unitlength]{0pt}{5\unitlength}
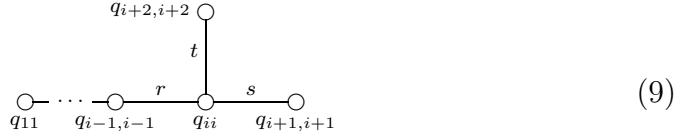
\begin{picture}(38,14)(0,3)
\put(1,3){\circle{2}}
\put(2,3){\line(1,0){2}}
\put(7,3){\scriptsize \makebox[0pt]{$\ldots $}}
\put(10,3){\line(1,0){2}}
\put(13,3){\circle{2}}
\put(14,3){\line(1,0){10}}
\put(25,3){\circle{2}}
\put(26,3){\line(1,0){10}}
\put(37,3){\circle{2}}
\put(25,4){\line(0,1){10}}
\put(25,15){\circle{2}}
\put(1,0){\makebox[0pt]{\scriptsize $q_{11}$}}
\put(13,0){\makebox[0pt]{\scriptsize $q_{i-1,i-1}$}}
\put(19,4){\makebox[0pt]{\scriptsize $r$}}
\put(25,0){\makebox[0pt]{\scriptsize $q_{ii}$}}
\put(31,4){\makebox[0pt]{\scriptsize $s$}}
\put(37,0){\makebox[0pt]{\scriptsize $q_{i+1,i+1}$}}
\put(24,9){\makebox[0pt][r]{\scriptsize $t$}}
\put(23,15){\makebox[0pt][r]{\scriptsize $q_{i+2,i+2}$}}
\end{picture}
 \end{align}
where $i\ge 2$.
Then for all $j\le i+2$ one has $q_{jj}\in \{t,t^{-1},-1\}$,
and for all $j\le i$ relation $q_{j,j+1}q_{j+1,j}\in \{t,t^{-1}\}$
holds. Moreover, if $q_{jj}=-1$ for some $j<i$ then
$(q_{j+1,j+1}+1)(q_{j,j+1}q_{j+1,j}q_{j+1,j+1}-1)=0$.
\end{lemma}

\begin{bew}
Perform induction on $i$. If $i=2$ then the claim can be obtained
immediately from Theorem~\ref{t-classrank4}. In the general case
consider the arithmetic root systems
$\roots (\chi ;\Ndbasis _2,\Ndbasis _3,\ldots ,\Ndbasis _{i+2})$ and
$\roots (\chi ;\Ndbasis _1,\ldots ,\Ndbasis _{i-2},
\Ndbasis _{i-1}+\Ndbasis _i,\Ndbasis _{i+1},\Ndbasis _{i+2})$.
To both of them one can apply the induction hypothesis and hence the first
part of the claim follows. If $q_{jj}=-1$ then consider
the arithmetic root system $(s_{\Ndbasis _j,E}^{-1}(\roots ),
\chi \circ (s_{\Ndbasis _j,E}\times s_{\Ndbasis _j,E}),E)$.
By Lemma~\ref{l-nooctopus} and the first part of the claim one obtains
that $q_{j,j+1}q_{j+1,j}=q_{j+1,j+1}\not=-1$ is not possible.
\end{bew}

\begin{satz}\label{s-Weqpathd}
Suppose that $(\chi ,E)$ is not of Cartan type. Then $(\roots ,\chi ,E)$
is Weyl equivalent to an arithmetic root system $(\roots ',\chi ',E)$
such that $\gDd _{\chi ',E}$ is a labeled path graph.
\end{satz}

\begin{bew}
Note that the proposition holds for arithmetic root systems of rank at most
four according to Proposition~\ref{s-Weqpath} and the classification
results in rank at most 3. The general case will be proven by induction
over the rank $d$.

Assume that $d\ge 5$ and $(\chi ,E)$ is not of Cartan type. Then there exists
a subset $E'\subset E$ with $d-1$ elements
such that $(\chi \restr{\ndZ E'\times \ndZ E'},E')$
is not of Cartan type and $\roots (\chi ;E')$ is a connected arithmetic
root system. By induction hypothesis it is Weyl equivalent to an arithmetic
root system $(\roots ',\chi ',E')$ such that $\gDd _{\chi ',E'}$ is a
labeled path graph. Hence, since $\extWBG _{\chi ,E}$ is full,
$(\roots ,\chi ,E)$ is Weyl equivalent to an arithmetic root system
with
unlabeled generalized Dynkin diagram
 \begin{align}\label{eq-gDdrankd}
\rule[-4\unitlength]{0pt}{5\unitlength}
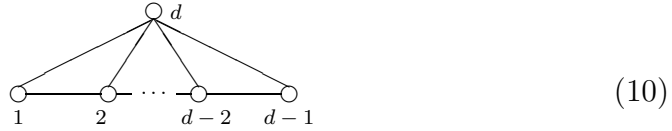
\begin{picture}(38,14)(0,3)
\put(1,4){\circle{2}}
\put(1,5){\line(2,1){18}}
\put(2,4){\line(1,0){10}}
\put(13,4){\circle{2}}
\put(13,5){\line(2,3){6}}
\put(14,4){\line(1,0){2}}
\put(19,4){\scriptsize \makebox[0pt]{$\ldots $}}
\put(22,4){\line(1,0){2}}
\put(25,4){\circle{2}}
\put(25,5){\line(-2,3){6}}
\put(19,15){\circle{2}}
\put(26,4){\line(1,0){10}}
\put(37,5){\line(-2,1){18}}
\put(37,4){\circle{2}}
\put(1,0){\makebox[0pt]{\scriptsize $1$}}
\put(12,0){\makebox[0pt]{\scriptsize $2$}}
\put(26,0){\makebox[0pt]{\scriptsize $d-2$}}
\put(37,0){\makebox[0pt]{\scriptsize $d-1$}}
\put(22,14){\makebox[0pt]{\scriptsize $d$}}
\end{picture}
 \end{align}
where it is allowed that some (but not all) nonhorizontal
edges in this diagram are omitted. Moreover, according to Lemmata
\ref{l-nodoubletriangle}, \ref{l-no4cycle} and \ref{l-nobigcycle}
all nonhorizontal edges in (\ref{eq-gDdrankd})
up to either two neighboring
or one single edge have to be omitted (that is labeled by 1).
If two neighboring nonhorizontal edges are labeled by numbers different
from 1,
 \begin{align}\label{eq-triangle+1or2legs}
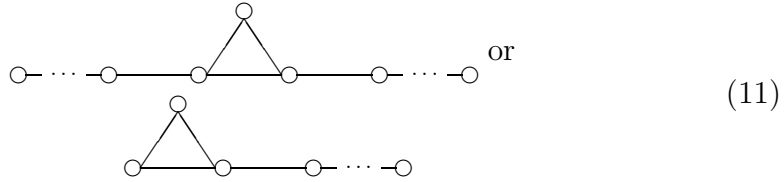
\begin{gathered}
\rule[-4\unitlength]{0pt}{5\unitlength}
\begin{picture}(62,7)(0,3)
\put(1,1){\circle{2}}
\put(2,1){\line(1,0){2}}
\put(7,1){\scriptsize \makebox[0pt]{$\ldots $}}
\put(10,1){\line(1,0){2}}
\put(13,1){\circle{2}}
\put(14,1){\line(1,0){10}}
\put(25,1){\circle{2}}
\put(26,1){\line(1,0){10}}
\put(37,1){\circle{2}}
\put(38,1){\line(1,0){10}}
\put(49,1){\circle{2}}
\put(50,1){\line(1,0){2}}
\put(55,1){\scriptsize \makebox[0pt]{$\ldots $}}
\put(58,1){\line(1,0){2}}
\put(61,1){\circle{2}}
\put(26,1){\line(2,3){5}}
\put(36,1){\line(-2,3){5}}
\put(31,9.5){\circle{2}}
\end{picture}
\text{ or }\\
\rule[-4\unitlength]{0pt}{5\unitlength}
\begin{picture}(38,7)(0,3)
\put(1,1){\circle{2}}
\put(2,1){\line(1,0){10}}
\put(13,1){\circle{2}}
\put(14,1){\line(1,0){10}}
\put(25,1){\circle{2}}
\put(26,1){\line(1,0){2}}
\put(31,1){\scriptsize \makebox[0pt]{$\ldots $}}
\put(34,1){\line(1,0){2}}
\put(37,1){\circle{2}}
\put(2,1){\line(2,3){5}}
\put(12,1){\line(-2,3){5}}
\put(7,9.5){\circle{2}}
\end{picture}
\end{gathered}
 \end{align}
then Lemmata \ref{l-rightofwayhas-1} and \ref{l-tent} imply the
claim of this proposition.

Thus without loss of generality
one can assume that $\gDd _{\chi ,E}$ is a connected labeled graph of the
form
 \begin{align}\label{eq-gDdEd}
\rule[-4\unitlength]{0pt}{5\unitlength}
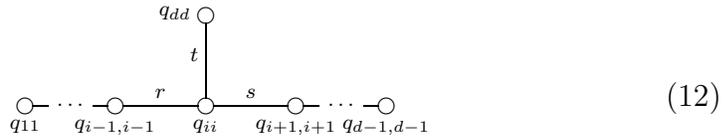
\begin{picture}(50,14)(0,3)
\put(1,3){\circle{2}}
\put(2,3){\line(1,0){2}}
\put(7,3){\scriptsize \makebox[0pt]{$\ldots $}}
\put(10,3){\line(1,0){2}}
\put(13,3){\circle{2}}
\put(14,3){\line(1,0){10}}
\put(25,3){\circle{2}}
\put(26,3){\line(1,0){10}}
\put(37,3){\circle{2}}
\put(38,3){\line(1,0){2}}
\put(43,3){\scriptsize \makebox[0pt]{$\ldots $}}
\put(46,3){\line(1,0){2}}
\put(49,3){\circle{2}}
\put(25,4){\line(0,1){10}}
\put(25,15){\circle{2}}
\put(1,0){\makebox[0pt]{\scriptsize $q_{11}$}}
\put(13,0){\makebox[0pt]{\scriptsize $q_{i-1,i-1}$}}
\put(19,4){\makebox[0pt]{\scriptsize $r$}}
\put(25,0){\makebox[0pt]{\scriptsize $q_{ii}$}}
\put(31,4){\makebox[0pt]{\scriptsize $s$}}
\put(37,0){\makebox[0pt]{\scriptsize $q_{i+1,i+1}$}}
\put(49,0){\makebox[0pt]{\scriptsize $q_{d-1,d-1}$}}
\put(24,9){\makebox[0pt][r]{\scriptsize $t$}}
\put(23,15){\makebox[0pt][r]{\scriptsize $q_{dd}$}}
\end{picture}
 \end{align}
where $1\le i\le d-1$ and $r,s,t\not=1$, and
$\roots (\chi ;\Ndbasis _1,\ldots ,\Ndbasis _i)$ is not
of Cartan type. If $i=1$ or $i=d-1$ then
$\gDd _{\chi ,E}$ is a labeled path graph and we are done.
Otherwise the finiteness of $\roots (\chi ;\Ndbasis _{i-1},
\Ndbasis _i,\Ndbasis _{i+1},\Ndbasis _d)$ and Theorem~\ref{t-classrank4}
imply that one of the following is true:
\begin{enumerate}
\item $r=s=t=q_{ii}=q_{i-1,i-1}=q_{i+1,i+1}=q_{dd}=-1$.
\item $t\not=-1$, $q_{ii}=-1$, $r,s\in \{t,t^{-1}\}$.
\item $t\not=-1$, $q_{ii}\not=-1$, $r,s\in \{t,t^{-1}\}$, and
there exist $m_{i,i-1},m_{i,i+1},m_{i,d}
\in \{1,2\}$ such that
$q_{ii}^{m_{i,i-1}}r=q_{ii}^{m_{i,i+1}}s=q_{ii}^{m_{id}}t=1$.
\end{enumerate}
Consider the first case.
By the same argument one obtains for $1<j<i$ by induction on
$i-j$ and by considering
$\roots (\chi ;\Ndbasis _1,\ldots ,\Ndbasis _{j-1},
\sum _{l=j}^i\Ndbasis _l,\Ndbasis _{i+1},\Ndbasis _d)$
that $q_{j-1,j-1}=q_{j-1,j}q_{j,j-1}=-1$.
By symmetry this gives that
$(\chi ,E)$ is of Cartan type which is a contradiction to the hypothesis
of the proposition.

In the second case for $T=s_{\Ndbasis _i,E}$ the arithmetic root system
$(T^{-1}(\roots ),\chi \circ (T\times T),E)$ 
is Weyl equivalent to $(\roots ,\chi ,E)$.
Because of Lemmata~\ref{l-notetraeder}, \ref{l-nodoubletriangle}
and \ref{l-noufo}
its unlabeled generalized Dynkin diagram is of the form as in
(\ref{eq-triangle+1or2legs})
and hence the claim of the proposition holds following the
arguments above.

It remains to prove the proposition in the third case
under the assumption
that $\roots (\chi ;\Ndbasis _1,\ldots ,\Ndbasis _i)$
is not of
Cartan type. This means that there exists $j<i$
such that $q_{jj}^mq_{j,j+1}q_{j+1,j}\not=1$ for all $m\in \ndZ $
or
$q_{jj}^mq_{j,j-1}q_{j-1,j}\not=1$ for all $m\in \ndZ $. Obviously, then
$q_{jj}\in R_{m+1}$ for some $m\in \ndN $. Suppose
that $j$ is maximal with this
property. Then Lemma~\ref{l-Dngraph} implies that $q_{jj}=-1$. Indeed,
otherwise relations
$q_{j,j-1}q_{j-1,j},q_{j,j+1}q_{j+1,j}\in \{q_{jj},q_{jj}^{-1}\}$
and $q_{jj}\in R_{m+1}$
give a contradiction to the choice of $j$.
Moreover, $q_{jj}=-1$ and Lemma~\ref{l-Dngraph} also imply that
$q_{j,j+1}q_{j+1,j}q_{j+1,j+1}=1$. In this case
$(s_{\Ndbasis _j,E}(\roots ),
\chi \circ (s_{\Ndbasis _j,E}\times s_{\Ndbasis _j,E}),E)$
is an arithmetic root system which in view of Lemma~\ref{l-nooctopus}
has a generalized Dynkin diagram as in (\ref{eq-gDdEd}). Further, the
properties of the third case apply, but
$\chi '(\Ndbasis _{j+1},\Ndbasis _{j+1})=-1$ where $\chi '=
\chi \circ (s_{\Ndbasis _j,E}\times s_{\Ndbasis _j,E})$.
Therefore by induction on $i-j$ one can show that $(\roots ,\chi ,E)$
is Weyl equivalent to an arithmetic root system as in case 2 and hence
the claim of the proposition is proven.
\end{bew}

\begin{lemma}\label{l-5chain}
If $\gDd _{\chi ,E}$ contains a subdiagram of the form
\begin{align*}
\Dchainfive{$q_{11}$}{$p$}{$q_{22}$}{$r$}{$q_{33}$}%
{$s$}{$q_{44}$}{$t$}{$q_{55}$}
\end{align*}
where $p,r,s,t\not=1$,
then $(q_{22}^2pr-1)(q_{44}^2st-1)=0$.
In particular, if for all $i,j\in \{1,2,\ldots ,d\}$
with $|i-j|\ge 2$
one has $q_{ij}q_{ji}=1$, and $I$
is the set of numbers $i\in \{2,3,\ldots ,d-1\}$ such that $q_{ii}^2
q_{i,i-1}q_{i-1,i}q_{i,i+1}q_{i+1,i}\not=1$, then $|I|\le 2$ and
$|i-j|\le 1$ for all $i,j\in I$.
\end{lemma}

\begin{bew}
Consider $\roots (\chi ;\Ndbasis _1+\Ndbasis _2,
\Ndbasis _2+\Ndbasis _3,\Ndbasis _3+\Ndbasis _4,
\Ndbasis _4+\Ndbasis _5)$. Its generalized Dynkin
diagram is
 \begin{align*}
\rule[-4\unitlength]{0pt}{5\unitlength}
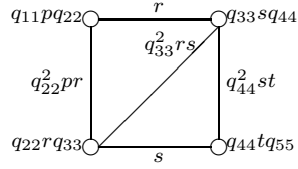
\begin{picture}(20,20)(0,3)
\put(1,3){\circle{2}}
\put(2,3){\line(1,0){15}}
\put(18,3){\circle{2}}
\put(2,3){\line(1,1){16}}
\put(18,20){\circle{2}}
\put(18,4){\line(0,1){15}}
\put(1,4){\line(0,1){15}}
\put(1,20){\circle{2}}
\put(2,20){\line(1,0){15}}
\put(0,20){\makebox[0pt][r]{\scriptsize $q_{11}pq_{22}$}}
\put(0,3){\makebox[0pt][r]{\scriptsize $q_{22}rq_{33}$}}
\put(19,20){\makebox[0pt][l]{\scriptsize $q_{33}sq_{44}$}}
\put(19,3){\makebox[0pt][l]{\scriptsize $q_{44}tq_{55}$}}
\put(10,1){\makebox[0pt]{\scriptsize $s$}}
\put(10,21){\makebox[0pt]{\scriptsize $r$}}
\put(0,11){\makebox[0pt][r]{\scriptsize $q_{22}^2pr$}}
\put(19,11){\makebox[0pt][l]{\scriptsize $q_{44}^2st$}}
\put(15,16){\makebox[0pt][r]{\scriptsize $q_{33}^2rs$}}
 \end{picture}
 \end{align*}
and hence Lemmata \ref{l-nodoubletriangle} and \ref{l-no4cycle} imply
the first part of the claim. For the second part, suppose that
$|I|\ge 2$ and $i,j\in I$ with $j-i\ge 2$. Then the finiteness
of $\roots (\chi ;
\Ndbasis _{i-1},\Ndbasis _i,\sum _{l=i+1}^{j-1}\Ndbasis _l,
\Ndbasis _j,\Ndbasis _{j+1})$ contradicts the first part of
the proof.
\end{bew}

\begin{lemma}\label{l-5chainmid}
If $\gDd _{\chi ,E}$ contains a subdiagram of the form
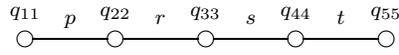
\begin{align*}
\Dchainfive{$q_{11}$}{$p$}{$q_{22}$}{$r$}{$q_{33}$}%
{$s$}{$q_{44}$}{$t$}{$q_{55}$}
\end{align*}
such that $q_{22}^2pr=q_{44}^2st=1$,
$q_{33}^2rs\not=1$, and $p,r,s,t\not=1$ then the following relations hold:
\begin{gather}\label{eq-5chainmiderel}
(p-s)(ps-1)=(r-t)(rt-1)=0,\qquad q_{33}^2rs\in \{r,r^{-1}\},\\
\label{eq-5chainmidvrel}
q_{22}rq_{33},q_{33}sq_{44}\in \{-1,1/q_{33}^2rs\}.
\end{gather}
Moreover, if $q_{33}sq_{44}=t$ then
equation $q_{55}t=1$ holds.
\end{lemma}

\begin{bew}
The arithmetic root system $\roots (\chi ;\Ndbasis _1+\Ndbasis _2,
\Ndbasis _2+\Ndbasis _3,\Ndbasis _3+\Ndbasis _4,\Ndbasis _5)$
has generalized Dynkin diagram
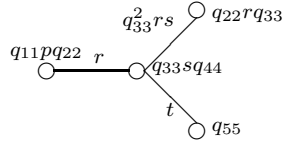
\begin{align*}
\rule[-9\unitlength]{0pt}{12\unitlength}
\begin{picture}(28,12)(0,9)
\put(2,10){\circle{2}}
\put(3,10){\line(1,0){10}}
\put(14,10){\circle{2}}
\put(15,10){\line(1,1){7}}
\put(15,10){\line(1,-1){7}}
\put(22,18){\circle{2}}
\put(22,2){\circle{2}}
\put(2,12){\makebox[0pt]{\scriptsize $q_{11}pq_{22}$}}
\put(9,11){\makebox[0pt]{\scriptsize $r$}}
\put(16,10){\makebox[0pt][l]{\scriptsize $q_{33}sq_{44}$}}
\put(19,16){\makebox[0pt][r]{\scriptsize $q_{33}^2rs$}}
\put(19,4){\makebox[0pt][r]{\scriptsize $t$}}
\put(24,17){\makebox[0pt][l]{\scriptsize $q_{22}rq_{33}$}}
\put(24,2){\makebox[0pt][l]{\scriptsize $q_{55}$}}
\end{picture}
\end{align*}
and hence by Theorem~\ref{t-classrank4} one gets relations $(r-t)(rt-1)=1$,
$q_{33}^2rs\in \{r,r^{-1}\}$, and $q_{22}rq_{33}\in \{-1,1/q_{33}^2rs\}$.
By the same reason one gets $q_{55}t=1$ whenever $q_{33}sq_{44}=t$.
The remaining relations follow by symmetry.
\end{bew}

\begin{bew}[ of Theorem~\ref{t-classrank>4}]
As in the proof of Theorem~\ref{t-classrank4} one shows that if
$\gDd _{\chi ,E}$ appears in Table~\ref{a-rb4} then
$\extWBG _{\chi ,E}$ is finite,
and two arithmetic root systems corresponding to generalized Dynkin
diagrams in Table~\ref{a-rb4} are Weyl equivalent if and only if these diagrams
appear in the same row of Table~\ref{a-rb4} and can be presented with the same
fixed parameter. For the list of elements $(T,E)\in W_{\chi ,E}$
mentioned in Proposition \ref{s-Wfinite} confer
Appendix~\ref{a-rb4}. In order to determine
the Weyl equivalence classes
Proposition~\ref{s-Weqsimplechains} is helpful.

It remains to prove that if $\extWBG _{\chi ,E}$ is full and finite
then $\gDd _{\chi ,E}$ appears in Table~\ref{a-rb4}.
This will be done by induction over $d$.
Further, Proposition~\ref{s-Weqpathd} is used to reduce to the case
where $\gDd _{\chi ,E}$ is a labeled path graph. More precisely,
assume that $q_{ij}q_{ji}=1$ whenever $|i-j|\ge 2$. Moreover,
\cite[Theorem\,1]{a-Heck04c} tells that if $(\chi ,E)$
is of Cartan
type then $\extWBG _{\chi ,E}$ is full and finite if and only if
$(\chi ,E)$ is of finite Cartan type.
These examples appear in rows~1, 3, 7, 8, 16, 20, and 22 of
Table~\ref{a-rb4}. Thus it
suffices to consider pairs $(\chi ,E)$
such that $\extWBG _{\chi ,E}$
is full and finite but $(\chi ,E)$ is not of Cartan type.

Since $\roots (\chi ;\Ndbasis _i,\Ndbasis _{i+1},\Ndbasis _{i+2},
\Ndbasis _{i+3})$ is finite for all $i\le d-3$,
by Theorem~\ref{t-classrank4}
any connected
subgraph of $\gDd _{\chi ,E}$ with four verices has to appear in
Table~\ref{a-r4}.
Using this fact it is sufficient to consider prolongations
of the labeled path graphs in Table~\ref{a-r4}
by (not necessarily simple) chains. 
This will be done in several steps.

\textit{Step~1.} Assume that
$(q_{dd}+1)(q_{dd}q_{d,d-1}q_{d-1,d}-1)\not=1$.
Then, since
$\roots (\chi ;\sum _{j=1}^i\Ndbasis _j,\Ndbasis _{i+1},
\sum _{j=i+2}^{d-1}\Ndbasis _j,\Ndbasis _d)$
is finite for all $i\le d-3$,
Theorem~\ref{t-classrank4} implies that
the generalized Dynkin diagram of
$\roots (\chi ;\Ndbasis _1,\Ndbasis _2,\ldots ,\Ndbasis _{d-2})$
is a simple chain. Moreover, since
$\roots (\chi ;\Ndbasis _{d-3},\Ndbasis _{d-2},
\Ndbasis _{d-1},\Ndbasis _d)$
is finite, either $\gDd _{\chi ,E}$
appears in one of rows~4, 5, and 6 of
Table~\ref{a-rb4},
or it is a prolongation of the first or the second graph
in row~17 of Table~\ref{a-r4} to the left by a simple chain.
Using the reflection $s_{\Ndbasis _d,E}$ the two latter
cases are Weyl equivalent,
and hence it is sufficient to consider the
second diagram in row~17 of Table~\ref{a-r4}.

Assume that $d=5$. Since
$\roots (\chi ;\Ndbasis _1,\Ndbasis _2+\Ndbasis _3,
\Ndbasis _4,\Ndbasis _5)$ is finite,
$\gDd _{\chi ,E}$ is of the form
\begin{align*}
\Dchainfive{$-\zeta $}{$-\zeta ^{-1}$}{$-\zeta $}{$-\zeta ^{-1}$}%
{$-\zeta $}{$-\zeta ^{-1}$}{$-1$}{$-1$}{$\zeta $}
\end{align*}
where $\zeta \in R_3$. Then $\roots (\chi ;\Ndbasis _1,
\Ndbasis _2+\Ndbasis _3,\Ndbasis _3+\Ndbasis _4,
\Ndbasis _4+\Ndbasis _5)$
has generalized Dynkin diagram
\begin{gather*}
%
%
\Dchainfour{}{$-\zeta $}{$-\zeta ^{-1}$}{$-\zeta $}{$-\zeta ^{-1}$}%
{$\zeta $}{$\zeta ^{-1}$}{$-1$}
\end{gather*}
which gives a contradiction to Theorem~\ref{t-classrank4}.

\textit{Step~2.} Assume that equations
$(q_{11}+1)(q_{11}q_{12}q_{21}-1)=0$ and
$(q_{dd}+1)(q_{dd}q_{d,d-1}q_{d-1,d}-1)=0$ hold,
and there exist at least
two integers $i$ with $1<i<d$ such that
$q_{ii}^2q_{i,i-1}q_{i-1,i}q_{i,i+1}q_{i+1,i}\not=1$. Then by
Lemma~\ref{l-5chain} there exist exactly two such numbers and they differ by 1.
Thus by Theorem~\ref{t-classrank4} $\gDd _{\chi ,E}$ is a prolongation of
one of the following diagrams in Table~\ref{a-r4}
to the left and/or to the right
by simple chains: diagram 6 in row~9,
diagram 4 in row~14, diagram 5 in row~17, diagram 6 in row~22,
diagram 7 in row~22.

\textit{Step~2.1. Prolongations of diagram 6 in row~9 of Table~\ref{a-r4}.}
One has $q\in k^\ast $ and $q^2,q^3\not=1$.
To the left:
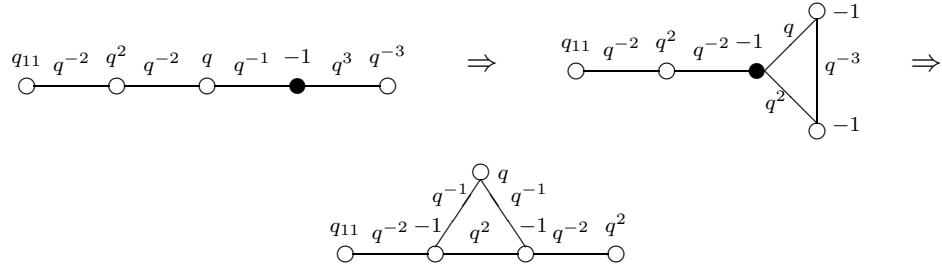
\begin{gather*}
%
%
\rule[-4\unitlength]{0pt}{5\unitlength}
\begin{picture}(50,4)(0,3)
\put(1,1){\circle{2}}
\put(2,1){\line(1,0){10}}
\put(13,1){\circle{2}}
\put(14,1){\line(1,0){10}}
\put(25,1){\circle{2}}
\put(26,1){\line(1,0){10}}
\put(37,1){\circle*{2}}
\put(38,1){\line(1,0){10}}
\put(49,1){\circle{2}}
\put(1,4){\makebox[0pt]{\scriptsize $q_{11}$}}
\put(7,3){\makebox[0pt]{\scriptsize $q^{-2}$}}
\put(13,4){\makebox[0pt]{\scriptsize $q^2$}}
\put(19,3){\makebox[0pt]{\scriptsize $q^{-2}$}}
\put(25,4){\makebox[0pt]{\scriptsize $q$}}
\put(31,3){\makebox[0pt]{\scriptsize $q^{-1}$}}
\put(37,4){\makebox[0pt]{\scriptsize $-1$}}
\put(43,3){\makebox[0pt]{\scriptsize $q^3$}}
\put(49,4){\makebox[0pt]{\scriptsize $q^{-3}$}}
\end{picture}
\qquad \Rightarrow \qquad
\rule[-10\unitlength]{0pt}{12\unitlength}
\begin{picture}(36,10)(0,9)
\put(1,9){\circle{2}}
\put(2,9){\line(1,0){10}}
\put(13,9){\circle{2}}
\put(14,9){\line(1,0){10}}
\put(25,9){\circle*{2}}
\put(26,9){\line(1,1){7}}
\put(26,9){\line(1,-1){7}}
\put(33,17){\circle{2}}
\put(33,1){\circle{2}}
\put(33,2){\line(0,1){14}}
\put(1,12){\makebox[0pt]{\scriptsize $q_{11}$}}
\put(7,11){\makebox[0pt]{\scriptsize $q^{-2}$}}
\put(13,12){\makebox[0pt]{\scriptsize $q^2$}}
\put(19,11){\makebox[0pt]{\scriptsize $q^{-2}$}}
\put(24,12){\makebox[0pt]{\scriptsize $-1$}}
\put(30,14){\makebox[0pt][r]{\scriptsize $q$}}
\put(29,4){\makebox[0pt][r]{\scriptsize $q^2$}}
\put(35,16){\makebox[0pt][l]{\scriptsize $-1$}}
\put(34,9){\makebox[0pt][l]{\scriptsize $q^{-3}$}}
\put(35,1){\makebox[0pt][l]{\scriptsize $-1$}}
\end{picture}
\qquad \Rightarrow \\
\rule[-4\unitlength]{0pt}{5\unitlength}
\begin{picture}(38,11)(0,3)
\put(1,1){\circle{2}}
\put(2,1){\line(1,0){10}}
\put(13,1){\circle{2}}
\put(13,2){\line(2,3){6}}
\put(14,1){\line(1,0){10}}
\put(25,1){\circle{2}}
\put(25,2){\line(-2,3){6}}
\put(19,12){\circle{2}}
\put(26,1){\line(1,0){10}}
\put(37,1){\circle{2}}
\put(1,4){\makebox[0pt]{\scriptsize $q_{11}$}}
\put(7,3){\makebox[0pt]{\scriptsize $q^{-2}$}}
\put(12,4){\makebox[0pt]{\scriptsize $-1$}}
\put(19,3){\makebox[0pt]{\scriptsize $q^2$}}
\put(26,4){\makebox[0pt]{\scriptsize $-1$}}
\put(31,3){\makebox[0pt]{\scriptsize $q^{-2}$}}
\put(37,4){\makebox[0pt]{\scriptsize $q^2$}}
\put(15,8){\makebox[0pt]{\scriptsize $q^{-1}$}}
\put(23,8){\makebox[0pt][l]{\scriptsize $q^{-1}$}}
\put(22,11){\makebox[0pt]{\scriptsize $q$}}
\end{picture}
\end{gather*}
This is a contradiction to $q\not=-1$ and Lemma~\ref{l-tent}.

Prolongations to the right:
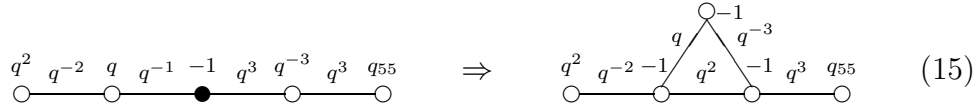
\begin{gather}\label{eq-9-6-l}
%
%
\rule[-4\unitlength]{0pt}{5\unitlength}
\begin{picture}(50,4)(0,3)
\put(1,1){\circle{2}}
\put(2,1){\line(1,0){10}}
\put(13,1){\circle{2}}
\put(14,1){\line(1,0){10}}
\put(25,1){\circle*{2}}
\put(26,1){\line(1,0){10}}
\put(37,1){\circle{2}}
\put(38,1){\line(1,0){10}}
\put(49,1){\circle{2}}
\put(1,4){\makebox[0pt]{\scriptsize $q^2$}}
\put(7,3){\makebox[0pt]{\scriptsize $q^{-2}$}}
\put(13,4){\makebox[0pt]{\scriptsize $q$}}
\put(19,3){\makebox[0pt]{\scriptsize $q^{-1}$}}
\put(25,4){\makebox[0pt]{\scriptsize $-1$}}
\put(31,3){\makebox[0pt]{\scriptsize $q^3$}}
\put(37,4){\makebox[0pt]{\scriptsize $q^{-3}$}}
\put(43,3){\makebox[0pt]{\scriptsize $q^3$}}
\put(49,4){\makebox[0pt]{\scriptsize $q_{55}$}}
\end{picture}
\qquad \Rightarrow \qquad
\rule[-4\unitlength]{0pt}{5\unitlength}
\begin{picture}(38,11)(0,3)
\put(1,1){\circle{2}}
\put(2,1){\line(1,0){10}}
\put(13,1){\circle{2}}
\put(13,2){\line(2,3){6}}
\put(14,1){\line(1,0){10}}
\put(25,1){\circle{2}}
\put(25,2){\line(-2,3){6}}
\put(19,12){\circle{2}}
\put(26,1){\line(1,0){10}}
\put(37,1){\circle{2}}
\put(1,4){\makebox[0pt]{\scriptsize $q^2$}}
\put(7,3){\makebox[0pt]{\scriptsize $q^{-2}$}}
\put(12,4){\makebox[0pt]{\scriptsize $-1$}}
\put(19,3){\makebox[0pt]{\scriptsize $q^2$}}
\put(26,4){\makebox[0pt]{\scriptsize $-1$}}
\put(31,3){\makebox[0pt]{\scriptsize $q^3$}}
\put(37,4){\makebox[0pt]{\scriptsize $q_{55}$}}
\put(15,8){\makebox[0pt]{\scriptsize $q$}}
\put(23,8){\makebox[0pt][l]{\scriptsize $q^{-3}$}}
\put(22,11){\makebox[0pt]{\scriptsize $-1$}}
\end{picture}
\end{gather}
Let $\Ndbasis _1,\Ndbasis _2,\Ndbasis _3,\Ndbasis _4$
correspond to the
four vertices of the last diagram lying on the same straight line,
and let $\Ndbasis _5$ correspond to the fifth vertex.
The generalized Dynkin diagram of
$\roots (\chi ;\Ndbasis _3+\Ndbasis _5,
\Ndbasis _2+\Ndbasis _5,\Ndbasis _3+\Ndbasis _4,
\Ndbasis _1+\Ndbasis _2)$ is
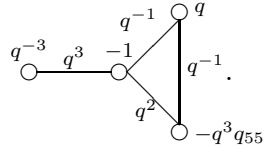
\begin{gather*}
\Drightofway{}{$q^{-3}$}{$q^3$}{$-1$}{$q^{-1}$}{$q^2$}{$q$}{$q^{-1}$}%
{$-q^3q_{55}$}.
\end{gather*}
By Lemma~\ref{l-rightofwayhas-1} one obtains that $q_{55}=q^{-3}$,
and Lemma~\ref{l-rightofway} gives that $q\in R_5$.
Then the first diagram in (\ref{eq-9-6-l}) coincides with the first
diagram in row~15 of Table~\ref{a-rb4}.

In a further prolongation to the right by a simple chain of length 1
one has $q\in R_5$.  Moreover,
since $\roots (\chi ;\Ndbasis _1,\Ndbasis _2,\Ndbasis _3,\Ndbasis _4+
\Ndbasis _5,\Ndbasis _6)$ is finite, one can only have
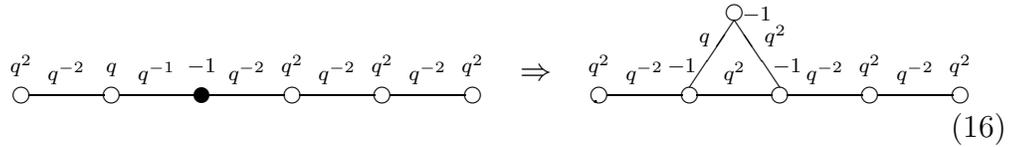
\begin{gather}\label{eq-9-6-le}
%
%
\rule[-4\unitlength]{0pt}{5\unitlength}
\begin{picture}(62,4)(0,3)
\put(1,1){\circle{2}}
\put(2,1){\line(1,0){10}}
\put(13,1){\circle{2}}
\put(14,1){\line(1,0){10}}
\put(25,1){\circle*{2}}
\put(26,1){\line(1,0){10}}
\put(37,1){\circle{2}}
\put(38,1){\line(1,0){10}}
\put(49,1){\circle{2}}
\put(50,1){\line(1,0){10}}
\put(61,1){\circle{2}}
\put(1,4){\makebox[0pt]{\scriptsize $q^2$}}
\put(7,3){\makebox[0pt]{\scriptsize $q^{-2}$}}
\put(13,4){\makebox[0pt]{\scriptsize $q$}}
\put(19,3){\makebox[0pt]{\scriptsize $q^{-1}$}}
\put(25,4){\makebox[0pt]{\scriptsize $-1$}}
\put(31,3){\makebox[0pt]{\scriptsize $q^{-2}$}}
\put(37,4){\makebox[0pt]{\scriptsize $q^2$}}
\put(43,3){\makebox[0pt]{\scriptsize $q^{-2}$}}
\put(49,4){\makebox[0pt]{\scriptsize $q^2$}}
\put(55,3){\makebox[0pt]{\scriptsize $q^{-2}$}}
\put(61,4){\makebox[0pt]{\scriptsize $q^2$}}
\end{picture}
\quad \Rightarrow \quad
\rule[-4\unitlength]{0pt}{5\unitlength}
\begin{picture}(50,11)(0,3)
\put(1,1){\circle{2}}
\put(2,1){\line(1,0){10}}
\put(13,1){\circle{2}}
\put(13,2){\line(2,3){6}}
\put(14,1){\line(1,0){10}}
\put(25,1){\circle{2}}
\put(25,2){\line(-2,3){6}}
\put(19,12){\circle{2}}
\put(26,1){\line(1,0){10}}
\put(37,1){\circle{2}}
\put(38,1){\line(1,0){10}}
\put(49,1){\circle{2}}
\put(1,4){\makebox[0pt]{\scriptsize $q^2$}}
\put(7,3){\makebox[0pt]{\scriptsize $q^{-2}$}}
\put(12,4){\makebox[0pt]{\scriptsize $-1$}}
\put(19,3){\makebox[0pt]{\scriptsize $q^2$}}
\put(26,4){\makebox[0pt]{\scriptsize $-1$}}
\put(31,3){\makebox[0pt]{\scriptsize $q^{-2}$}}
\put(37,4){\makebox[0pt]{\scriptsize $q^2$}}
\put(43,3){\makebox[0pt]{\scriptsize $q^{-2}$}}
\put(49,4){\makebox[0pt]{\scriptsize $q^2$}}
\put(15,8){\makebox[0pt]{\scriptsize $q$}}
\put(23,8){\makebox[0pt][l]{\scriptsize $q^2$}}
\put(22,11){\makebox[0pt]{\scriptsize $-1$}}.
\end{picture}
\end{gather}
However for the last diagram, where $\Ndbasis _1,\ldots ,\Ndbasis _5$
correspond to the vertices lying on the same straight line and $\Ndbasis _6$
corresponds to the remaining vertex,
the generalized Dynkin diagram of $\roots (\chi ;\Ndbasis _1+\Ndbasis _2,
\Ndbasis _2+\Ndbasis _6,\Ndbasis _3+\Ndbasis _6,\Ndbasis _3+\Ndbasis _4,
\Ndbasis _5)$ has the form
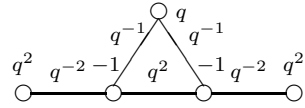
\begin{gather*}
\rule[-4\unitlength]{0pt}{5\unitlength}
\begin{picture}(38,11)(0,3)
\put(1,1){\circle{2}}
\put(2,1){\line(1,0){10}}
\put(13,1){\circle{2}}
\put(13,2){\line(2,3){6}}
\put(14,1){\line(1,0){10}}
\put(25,1){\circle{2}}
\put(25,2){\line(-2,3){6}}
\put(19,12){\circle{2}}
\put(26,1){\line(1,0){10}}
\put(37,1){\circle{2}}
\put(1,4){\makebox[0pt]{\scriptsize $q^2$}}
\put(7,3){\makebox[0pt]{\scriptsize $q^{-2}$}}
\put(12,4){\makebox[0pt]{\scriptsize $-1$}}
\put(19,3){\makebox[0pt]{\scriptsize $q^2$}}
\put(26,4){\makebox[0pt]{\scriptsize $-1$}}
\put(31,3){\makebox[0pt]{\scriptsize $q^{-2}$}}
\put(37,4){\makebox[0pt]{\scriptsize $q^2$}}
\put(15,8){\makebox[0pt]{\scriptsize $q^{-1}$}}
\put(23,8){\makebox[0pt][l]{\scriptsize $q^{-1}$}}
\put(22,11){\makebox[0pt]{\scriptsize $q$}}
\end{picture}
\end{gather*}
which is a contradiction to Lemma~\ref{l-tent} and $q\in R_5$.

\textit{Step~2.2. Prolongations of diagram 4 in row~14 of Table~\ref{a-r4}.}
By symmetry it is enough to prolong to the left or to both
directions by simple chains.
Assume that $q^2\not=1$ and $(q_{11}+1)(q_{11}-q)=0$.
Then Weyl equivalence gives
\begin{gather}\label{eq-14-4-l}
%
%
\rule[-4\unitlength]{0pt}{5\unitlength}
\begin{picture}(50,4)(0,3)
\put(1,1){\circle{2}}
\put(2,1){\line(1,0){10}}
\put(13,1){\circle{2}}
\put(14,1){\line(1,0){10}}
\put(25,1){\circle*{2}}
\put(26,1){\line(1,0){10}}
\put(37,1){\circle{2}}
\put(38,1){\line(1,0){10}}
\put(49,1){\circle{2}}
\put(1,4){\makebox[0pt]{\scriptsize $q_{11}$}}
\put(7,3){\makebox[0pt]{\scriptsize $q^{-1}$}}
\put(13,4){\makebox[0pt]{\scriptsize $q$}}
\put(19,3){\makebox[0pt]{\scriptsize $q^{-1}$}}
\put(25,4){\makebox[0pt]{\scriptsize $-1$}}
\put(31,3){\makebox[0pt]{\scriptsize $-1$}}
\put(37,4){\makebox[0pt]{\scriptsize $-1$}}
\put(43,3){\makebox[0pt]{\scriptsize $-q$}}
\put(49,4){\makebox[0pt]{\scriptsize $-q^{-1}$}}
\end{picture}
\qquad \Rightarrow \qquad
\rule[-4\unitlength]{0pt}{5\unitlength}
\begin{picture}(38,11)(0,3)
\put(1,1){\circle{2}}
\put(2,1){\line(1,0){10}}
\put(13,1){\circle{2}}
\put(13,2){\line(2,3){6}}
\put(14,1){\line(1,0){10}}
\put(25,1){\circle{2}}
\put(25,2){\line(-2,3){6}}
\put(19,12){\circle{2}}
\put(26,1){\line(1,0){10}}
\put(37,1){\circle{2}}
\put(1,4){\makebox[0pt]{\scriptsize $q_{11}$}}
\put(7,3){\makebox[0pt]{\scriptsize $q^{-1}$}}
\put(12,4){\makebox[0pt]{\scriptsize $-1$}}
\put(19,3){\makebox[0pt]{\scriptsize $-q^{-1}$}}
\put(26,4){\makebox[0pt]{\scriptsize $-1$}}
\put(31,3){\makebox[0pt]{\scriptsize $-q$}}
\put(37,4){\makebox[0pt]{\scriptsize $-q^{-1}$}}
\put(15,8){\makebox[0pt]{\scriptsize $q$}}
\put(23,8){\makebox[0pt][l]{\scriptsize $-1$}}
\put(22,11){\makebox[0pt]{\scriptsize $-1$}}
\end{picture}
\end{gather}
With the same notation as below (\ref{eq-9-6-l}) the finiteness of
$\roots (\chi ;\Ndbasis _2+\Ndbasis _5,\Ndbasis _3+\Ndbasis _5,
\Ndbasis _3+\Ndbasis _4,\Ndbasis _1+\Ndbasis _2)$,
which has the generalized
Dynkin diagram
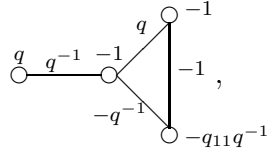
\begin{align*}
\Drightofway{}{$q$}{$q^{-1}$}{$-1$}{$q$}{$-q^{-1}$}{$-1$}{$-1$}%
{$-q_{11}q^{-1}$},
\end{align*}
and Theorem~\ref{t-classrank4} imply that $q_{11}\not=-1$. Therefore it remains
to check the case when $q_{11}=q$. Starting again with the original diagram,
Weyl equivalence gives also
\begin{gather}\label{eq-14-4-le}
%
%
\rule[-4\unitlength]{0pt}{5\unitlength}
\begin{picture}(50,4)(0,3)
\put(1,1){\circle{2}}
\put(2,1){\line(1,0){10}}
\put(13,1){\circle{2}}
\put(14,1){\line(1,0){10}}
\put(25,1){\circle{2}}
\put(26,1){\line(1,0){10}}
\put(37,1){\circle*{2}}
\put(38,1){\line(1,0){10}}
\put(49,1){\circle{2}}
\put(1,4){\makebox[0pt]{\scriptsize $q$}}
\put(7,3){\makebox[0pt]{\scriptsize $q^{-1}$}}
\put(13,4){\makebox[0pt]{\scriptsize $q$}}
\put(19,3){\makebox[0pt]{\scriptsize $q^{-1}$}}
\put(25,4){\makebox[0pt]{\scriptsize $-1$}}
\put(31,3){\makebox[0pt]{\scriptsize $-1$}}
\put(37,4){\makebox[0pt]{\scriptsize $-1$}}
\put(43,3){\makebox[0pt]{\scriptsize $-q$}}
\put(49,4){\makebox[0pt]{\scriptsize $-q^{-1}$}}
\end{picture}
\qquad \Rightarrow \qquad
\rule[-10\unitlength]{0pt}{12\unitlength}
\begin{picture}(36,10)(0,9)
\put(1,9){\circle{2}}
\put(2,9){\line(1,0){10}}
\put(13,9){\circle{2}}
\put(14,9){\line(1,0){10}}
\put(25,9){\circle{2}}
\put(26,9){\line(1,1){7}}
\put(26,9){\line(1,-1){7}}
\put(33,17){\circle{2}}
\put(33,1){\circle{2}}
\put(33,2){\line(0,1){14}}
\put(1,12){\makebox[0pt]{\scriptsize $q$}}
\put(7,11){\makebox[0pt]{\scriptsize $q^{-1}$}}
\put(13,12){\makebox[0pt]{\scriptsize $q$}}
\put(19,11){\makebox[0pt]{\scriptsize $q^{-1}$}}
\put(24,12){\makebox[0pt]{\scriptsize $-1$}}
\put(30,14){\makebox[0pt][r]{\scriptsize $-1$}}
\put(29,4){\makebox[0pt][r]{\scriptsize $q$}}
\put(35,16){\makebox[0pt][l]{\scriptsize $-1$}}
\put(34,9){\makebox[0pt][l]{\scriptsize $-q^{-1}$}}
\put(35,1){\makebox[0pt][l]{\scriptsize $-1$}}
\end{picture}
.
\end{gather}
Assign the nodes of the right diagram from left to right and
top to bottom to the basis vectors $\Ndbasis _i$, where $1\le i\le 5$.
Applying Lemma~\ref{l-rightofway} to
$\roots (\chi ;\Ndbasis _1+\Ndbasis _2,
\Ndbasis _2+\Ndbasis _3,\Ndbasis _3+\Ndbasis _4,\Ndbasis _5)$,
which has
as generalized Dynkin diagram the labeled graph
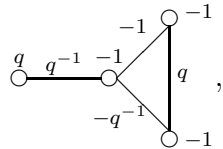
\begin{align*}
\Drightofway{}{$q$}{$q^{-1}$}{$-1$}{$-1$}%
{$-q^{-1}$}{$-1$}{$q$}{$-1$},
\end{align*}
one obtains that $q\in R_4$. In this case (see the left graph in
(\ref{eq-14-4-le}))
$\gDd _{\chi ,E}$ appears in row~14 of Table~\ref{a-rb4}.

Next it will be shown that there is no prolongation
of the left graph
in (\ref{eq-14-4-le}) to the right. To this end note that 
since
$\roots (\chi ;\Ndbasis _2,\Ndbasis _3,\Ndbasis _4,\Ndbasis _5,
\Ndbasis _6)$ is finite, by the previous arguments
any prolongation to the right has to be of the form
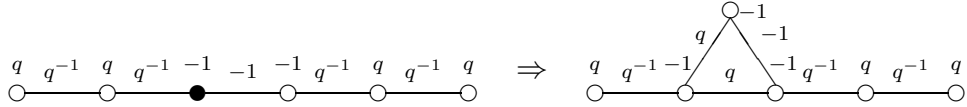
\begin{gather*}
%
%
\rule[-4\unitlength]{0pt}{5\unitlength}
\begin{picture}(62,4)(0,3)
\put(1,1){\circle{2}}
\put(2,1){\line(1,0){10}}
\put(13,1){\circle{2}}
\put(14,1){\line(1,0){10}}
\put(25,1){\circle*{2}}
\put(26,1){\line(1,0){10}}
\put(37,1){\circle{2}}
\put(38,1){\line(1,0){10}}
\put(49,1){\circle{2}}
\put(50,1){\line(1,0){10}}
\put(61,1){\circle{2}}
\put(1,4){\makebox[0pt]{\scriptsize $q$}}
\put(7,3){\makebox[0pt]{\scriptsize $q^{-1}$}}
\put(13,4){\makebox[0pt]{\scriptsize $q$}}
\put(19,3){\makebox[0pt]{\scriptsize $q^{-1}$}}
\put(25,4){\makebox[0pt]{\scriptsize $-1$}}
\put(31,3){\makebox[0pt]{\scriptsize $-1$}}
\put(37,4){\makebox[0pt]{\scriptsize $-1$}}
\put(43,3){\makebox[0pt]{\scriptsize $q^{-1}$}}
\put(49,4){\makebox[0pt]{\scriptsize $q$}}
\put(55,3){\makebox[0pt]{\scriptsize $q^{-1}$}}
\put(61,4){\makebox[0pt]{\scriptsize $q$}}
\end{picture}
\quad \Rightarrow \quad
\rule[-4\unitlength]{0pt}{5\unitlength}
\begin{picture}(50,11)(0,3)
\put(1,1){\circle{2}}
\put(2,1){\line(1,0){10}}
\put(13,1){\circle{2}}
\put(13,2){\line(2,3){6}}
\put(14,1){\line(1,0){10}}
\put(25,1){\circle{2}}
\put(25,2){\line(-2,3){6}}
\put(19,12){\circle{2}}
\put(26,1){\line(1,0){10}}
\put(37,1){\circle{2}}
\put(38,1){\line(1,0){10}}
\put(49,1){\circle{2}}
\put(1,4){\makebox[0pt]{\scriptsize $q$}}
\put(7,3){\makebox[0pt]{\scriptsize $q^{-1}$}}
\put(12,4){\makebox[0pt]{\scriptsize $-1$}}
\put(19,3){\makebox[0pt]{\scriptsize $q$}}
\put(26,4){\makebox[0pt]{\scriptsize $-1$}}
\put(31,3){\makebox[0pt]{\scriptsize $q^{-1}$}}
\put(37,4){\makebox[0pt]{\scriptsize $q$}}
\put(43,3){\makebox[0pt]{\scriptsize $q^{-1}$}}
\put(49,4){\makebox[0pt]{\scriptsize $q$}}
\put(15,8){\makebox[0pt]{\scriptsize $q$}}
\put(23,8){\makebox[0pt][l]{\scriptsize $-1$}}
\put(22,11){\makebox[0pt]{\scriptsize $-1$}}
\end{picture}
\end{gather*}
where $q\in R_4$.
Using the conventions below (\ref{eq-9-6-le}) for the last graph,
the generalized Dynkin diagram of
$\roots (\chi ;\Ndbasis _1+\Ndbasis _2,\Ndbasis _2+\Ndbasis _3,
\Ndbasis _3+\Ndbasis _6,\Ndbasis _4,\Ndbasis _5)$
has the form
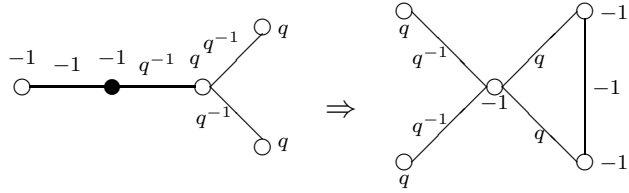
\begin{align*}
\rule[-10\unitlength]{0pt}{12\unitlength}
\begin{picture}(36,14)(0,9)
\put(1,13){\circle{2}}
\put(2,13){\line(1,0){10}}
\put(13,13){\circle*{2}}
\put(14,13){\line(1,0){10}}
\put(25,13){\circle{2}}
\put(26,13){\line(1,1){7}}
\put(26,13){\line(1,-1){7}}
\put(33,21){\circle{2}}
\put(33,5){\circle{2}}
\put(1,16){\makebox[0pt]{\scriptsize $-1$}}
\put(7,15){\makebox[0pt]{\scriptsize $-1$}}
\put(13,16){\makebox[0pt]{\scriptsize $-1$}}
\put(19,15){\makebox[0pt]{\scriptsize $q^{-1}$}}
\put(24,16){\makebox[0pt]{\scriptsize $q$}}
\put(30,18){\makebox[0pt][r]{\scriptsize $q^{-1}$}}
\put(29,8){\makebox[0pt][r]{\scriptsize $q^{-1}$}}
\put(35,20){\makebox[0pt][l]{\scriptsize $q$}}
\put(35,5){\makebox[0pt][l]{\scriptsize $q$}}
\end{picture}
\quad \Rightarrow \quad
\rule[-10\unitlength]{0pt}{5\unitlength}
\begin{picture}(29,17)(0,9)
\put(1,3){\circle{2}}
\put(2,3){\line(1,1){10}}
\put(13,13){\circle{2}}
\put(1,23){\circle{2}}
\put(2,23){\line(1,-1){10}}
\put(14,13){\line(1,1){10}}
\put(14,13){\line(1,-1){10}}
\put(25,23){\circle{2}}
\put(25,3){\circle{2}}
\put(25,4){\line(0,1){18}}
\put(1,0){\makebox[0pt]{\scriptsize $q$}}
\put(1,20){\makebox[0pt]{\scriptsize $q$}}
\put(7,7){\makebox[0pt][r]{\scriptsize $q^{-1}$}}
\put(7,16){\makebox[0pt][r]{\scriptsize $q^{-1}$}}
\put(13,10){\makebox[0pt]{\scriptsize $-1$}}
\put(19,6){\makebox[0pt]{\scriptsize $q$}}
\put(19,16){\makebox[0pt]{\scriptsize $q$}}
\put(27,2){\makebox{\scriptsize $-1$}}
\put(26,12){\makebox{\scriptsize $-1$}}
\put(27,22){\makebox{\scriptsize $-1$}}
\end{picture}
\end{align*}
which is a contradiction to Lemma~\ref{l-nooctopus}.

It remains to consider prolongations of the left graph in (\ref{eq-14-4-le})
to the left, where $q\in R_4$. Since
$\roots (\chi ;E\setminus \{\Ndbasis _2\})$ is finite, one has to have
$q_{11}=q$ and $q_{12}q_{21}=q^{-1}$.
This prolongation of length 1 is the first diagram in row~19 of
Table~\ref{a-rb4}.
The prolongation of length 2 looks like
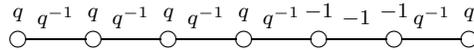
\begin{align*}
\rule[-4\unitlength]{0pt}{5\unitlength}
\begin{picture}(62,4)(0,3)
\put(1,1){\circle{2}}
\put(2,1){\line(1,0){8}}
\put(11,1){\circle{2}}
\put(12,1){\line(1,0){8}}
\put(21,1){\circle{2}}
\put(22,1){\line(1,0){8}}
\put(31,1){\circle{2}}
\put(32,1){\line(1,0){8}}
\put(41,1){\circle{2}}
\put(42,1){\line(1,0){8}}
\put(51,1){\circle{2}}
\put(52,1){\line(1,0){8}}
\put(61,1){\circle{2}}
\put(1,4){\makebox[0pt]{\scriptsize $q$}}
\put(6,3){\makebox[0pt]{\scriptsize $q^{-1}$}}
\put(11,4){\makebox[0pt]{\scriptsize $q$}}
\put(16,3){\makebox[0pt]{\scriptsize $q^{-1}$}}
\put(21,4){\makebox[0pt]{\scriptsize $q$}}
\put(26,3){\makebox[0pt]{\scriptsize $q^{-1}$}}
\put(31,4){\makebox[0pt]{\scriptsize $q$}}
\put(36,3){\makebox[0pt]{\scriptsize $q^{-1}$}}
\put(41,4){\makebox[0pt]{\scriptsize $-1$}}
\put(46,3){\makebox[0pt]{\scriptsize $-1$}}
\put(51,4){\makebox[0pt]{\scriptsize $-1$}}
\put(56,3){\makebox[0pt]{\scriptsize $q^{-1}$}}
\put(61,4){\makebox[0pt]{\scriptsize $q$}}
\end{picture}
\end{align*}
in which case the generalized Dynkin diagram of $\roots (\chi ;
\Ndbasis _1,\Ndbasis _2+\Ndbasis _3,\Ndbasis _3+\Ndbasis _4,
\Ndbasis _4+\Ndbasis _5,\Ndbasis _5+\Ndbasis _6,
\Ndbasis _6+\Ndbasis _7)$ is
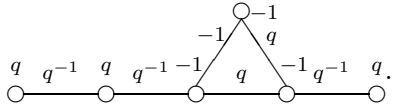
\begin{align*}
\rule[-4\unitlength]{0pt}{5\unitlength}
\begin{picture}(50,11)(0,3)
\put(1,1){\circle{2}}
\put(2,1){\line(1,0){10}}
\put(13,1){\circle{2}}
\put(14,1){\line(1,0){10}}
\put(25,1){\circle{2}}
\put(25,2){\line(2,3){6}}
\put(26,1){\line(1,0){10}}
\put(37,1){\circle{2}}
\put(37,2){\line(-2,3){6}}
\put(31,12){\circle{2}}
\put(38,1){\line(1,0){10}}
\put(49,1){\circle{2}}
\put(1,4){\makebox[0pt]{\scriptsize $q$}}
\put(7,3){\makebox[0pt]{\scriptsize $q^{-1}$}}
\put(13,4){\makebox[0pt]{\scriptsize $q$}}
\put(19,3){\makebox[0pt]{\scriptsize $q^{-1}$}}
\put(24,4){\makebox[0pt]{\scriptsize $-1$}}
\put(31,3){\makebox[0pt]{\scriptsize $q$}}
\put(38,4){\makebox[0pt]{\scriptsize $-1$}}
\put(43,3){\makebox[0pt]{\scriptsize $q^{-1}$}}
\put(49,4){\makebox[0pt]{\scriptsize $q$}}
\put(27,8){\makebox[0pt]{\scriptsize $-1$}}
\put(35,8){\makebox[0pt]{\scriptsize $q$}}
\put(34,11){\makebox[0pt]{\scriptsize $-1$}}
\end{picture}
.
\end{align*}
This is a contradiction to Theorem~\ref{t-classrank>4} for $d=6$.

\textit{Step~2.3. Prolongations of diagram 5 in row~17 of Table~\ref{a-r4}.}
Let $\zeta \in R_3$.

Prolongations to the left:
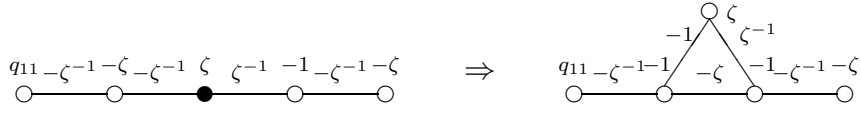
\begin{gather*}\label{eq-17-5-l}
%
%
\rule[-4\unitlength]{0pt}{5\unitlength}
\begin{picture}(50,4)(0,3)
\put(1,1){\circle{2}}
\put(2,1){\line(1,0){10}}
\put(13,1){\circle{2}}
\put(14,1){\line(1,0){10}}
\put(25,1){\circle*{2}}
\put(26,1){\line(1,0){10}}
\put(37,1){\circle{2}}
\put(38,1){\line(1,0){10}}
\put(49,1){\circle{2}}
\put(1,4){\makebox[0pt]{\scriptsize $q_{11}$}}
\put(7,3){\makebox[0pt]{\scriptsize $-\zeta ^{-1}$}}
\put(13,4){\makebox[0pt]{\scriptsize $-\zeta $}}
\put(19,3){\makebox[0pt]{\scriptsize $-\zeta ^{-1}$}}
\put(25,4){\makebox[0pt]{\scriptsize $\zeta $}}
\put(31,3){\makebox[0pt]{\scriptsize $\zeta ^{-1}$}}
\put(37,4){\makebox[0pt]{\scriptsize $-1$}}
\put(43,3){\makebox[0pt]{\scriptsize $-\zeta ^{-1}$}}
\put(49,4){\makebox[0pt]{\scriptsize $-\zeta $}}
\end{picture}
\qquad \Rightarrow \qquad
\rule[-4\unitlength]{0pt}{5\unitlength}
\begin{picture}(38,11)(0,3)
\put(1,1){\circle{2}}
\put(2,1){\line(1,0){10}}
\put(13,1){\circle{2}}
\put(13,2){\line(2,3){6}}
\put(14,1){\line(1,0){10}}
\put(25,1){\circle{2}}
\put(25,2){\line(-2,3){6}}
\put(19,12){\circle{2}}
\put(26,1){\line(1,0){10}}
\put(37,1){\circle{2}}
\put(1,4){\makebox[0pt]{\scriptsize $q_{11}$}}
\put(7,3){\makebox[0pt]{\scriptsize $-\zeta ^{-1}$}}
\put(12,4){\makebox[0pt]{\scriptsize $-1$}}
\put(19,3){\makebox[0pt]{\scriptsize $-\zeta $}}
\put(26,4){\makebox[0pt]{\scriptsize $-1$}}
\put(31,3){\makebox[0pt]{\scriptsize $-\zeta ^{-1}$}}
\put(37,4){\makebox[0pt]{\scriptsize $-\zeta $}}
\put(15,8){\makebox[0pt]{\scriptsize $-1$}}
\put(23,8){\makebox[0pt][l]{\scriptsize $\zeta ^{-1}$}}
\put(22,11){\makebox[0pt]{\scriptsize $\zeta $}}
\end{picture}
\end{gather*}
This is a contradiction to Lemma~\ref{l-tent}.

To the right:
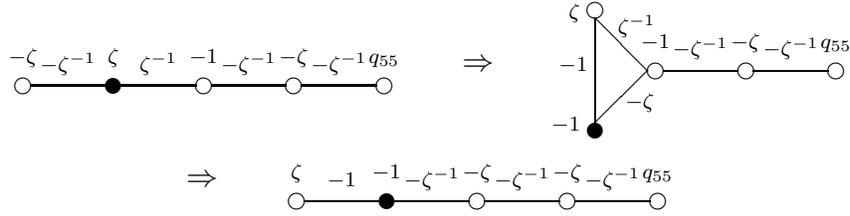
\begin{gather*}\label{eq-17-5-r}
%
%
\rule[-4\unitlength]{0pt}{5\unitlength}
\begin{picture}(50,4)(0,3)
\put(1,1){\circle{2}}
\put(2,1){\line(1,0){10}}
\put(13,1){\circle*{2}}
\put(14,1){\line(1,0){10}}
\put(25,1){\circle{2}}
\put(26,1){\line(1,0){10}}
\put(37,1){\circle{2}}
\put(38,1){\line(1,0){10}}
\put(49,1){\circle{2}}
\put(1,4){\makebox[0pt]{\scriptsize $-\zeta $}}
\put(7,3){\makebox[0pt]{\scriptsize $-\zeta ^{-1}$}}
\put(13,4){\makebox[0pt]{\scriptsize $\zeta $}}
\put(19,3){\makebox[0pt]{\scriptsize $\zeta ^{-1}$}}
\put(25,4){\makebox[0pt]{\scriptsize $-1$}}
\put(31,3){\makebox[0pt]{\scriptsize $-\zeta ^{-1}$}}
\put(37,4){\makebox[0pt]{\scriptsize $-\zeta $}}
\put(43,3){\makebox[0pt]{\scriptsize $-\zeta ^{-1}$}}
\put(49,4){\makebox[0pt]{\scriptsize $q_{55}$}}
\end{picture}
\qquad \Rightarrow \qquad
\rule[-10\unitlength]{0pt}{12\unitlength}
\begin{picture}(36,10)(0,9)
\put(36,9){\circle{2}}
\put(35,9){\line(-1,0){10}}
\put(24,9){\circle{2}}
\put(23,9){\line(-1,0){10}}
\put(12,9){\circle{2}}
\put(11,9){\line(-1,1){7}}
\put(11,9){\line(-1,-1){7}}
\put(4,17){\circle{2}}
\put(4,1){\circle*{2}}
\put(4,2){\line(0,1){14}}
\put(36,12){\makebox[0pt]{\scriptsize $q_{55}$}}
\put(30,11){\makebox[0pt]{\scriptsize $-\zeta ^{-1}$}}
\put(24,12){\makebox[0pt]{\scriptsize $-\zeta $}}
\put(18,11){\makebox[0pt]{\scriptsize $-\zeta ^{-1}$}}
\put(12,12){\makebox[0pt]{\scriptsize $-1$}}
\put(7,14){\makebox[0pt][l]{\scriptsize $\zeta ^{-1}$}}
\put(8,4){\makebox[0pt][l]{\scriptsize $-\zeta $}}
\put(2,16){\makebox[0pt][r]{\scriptsize $\zeta $}}
\put(3,9){\makebox[0pt][r]{\scriptsize $-1$}}
\put(2,1){\makebox[0pt][r]{\scriptsize $-1$}}
\end{picture}\\
\Rightarrow \qquad
%
%
\rule[-4\unitlength]{0pt}{5\unitlength}
\begin{picture}(50,4)(0,3)
\put(1,1){\circle{2}}
\put(2,1){\line(1,0){10}}
\put(13,1){\circle*{2}}
\put(14,1){\line(1,0){10}}
\put(25,1){\circle{2}}
\put(26,1){\line(1,0){10}}
\put(37,1){\circle{2}}
\put(38,1){\line(1,0){10}}
\put(49,1){\circle{2}}
\put(1,4){\makebox[0pt]{\scriptsize $\zeta $}}
\put(7,3){\makebox[0pt]{\scriptsize $-1$}}
\put(13,4){\makebox[0pt]{\scriptsize $-1$}}
\put(19,3){\makebox[0pt]{\scriptsize $-\zeta ^{-1}$}}
\put(25,4){\makebox[0pt]{\scriptsize $-\zeta $}}
\put(31,3){\makebox[0pt]{\scriptsize $-\zeta ^{-1}$}}
\put(37,4){\makebox[0pt]{\scriptsize $-\zeta $}}
\put(43,3){\makebox[0pt]{\scriptsize $-\zeta ^{-1}$}}
\put(49,4){\makebox[0pt]{\scriptsize $q_{55}$}}
\end{picture}
\end{gather*}
The latter diagram does not correspond to an arithmetic root system,
as shown in Step~1.

\textit{Step~2.4. Prolongations of diagram 6 in row~22
of Table~\ref{a-r4}.}
Assume that $\zeta \in R_4$.
For a prolongation to the left by a simple chain of length 1
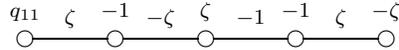
\begin{gather*}
%
%
\rule[-4\unitlength]{0pt}{5\unitlength}
\begin{picture}(50,4)(0,3)
\put(1,1){\circle{2}}
\put(2,1){\line(1,0){10}}
\put(13,1){\circle{2}}
\put(14,1){\line(1,0){10}}
\put(25,1){\circle{2}}
\put(26,1){\line(1,0){10}}
\put(37,1){\circle{2}}
\put(38,1){\line(1,0){10}}
\put(49,1){\circle{2}}
\put(1,4){\makebox[0pt]{\scriptsize $q_{11}$}}
\put(7,3){\makebox[0pt]{\scriptsize $\zeta $}}
\put(13,4){\makebox[0pt]{\scriptsize $-1$}}
\put(19,3){\makebox[0pt]{\scriptsize $-\zeta $}}
\put(25,4){\makebox[0pt]{\scriptsize $\zeta $}}
\put(31,3){\makebox[0pt]{\scriptsize $-1$}}
\put(37,4){\makebox[0pt]{\scriptsize $-1$}}
\put(43,3){\makebox[0pt]{\scriptsize $\zeta $}}
\put(49,4){\makebox[0pt]{\scriptsize $-\zeta $}}
\end{picture}
\end{gather*}
$\roots (\chi ;\Ndbasis _1+\Ndbasis _2,\Ndbasis _2+\Ndbasis _3,
\Ndbasis _3+\Ndbasis _4,\Ndbasis _5)$ has
the generalized Dynkin diagram
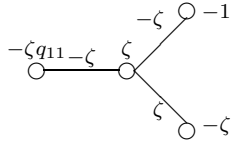
\begin{align*}
\Dthreefork{}{$-\zeta q_{11}$}{$-\zeta $}%
{$\zeta $}{$-\zeta $}{$\zeta $}{$-1$}{$-\zeta $}
\end{align*}
which is a contradiction to Theorem~\ref{t-classrank4}.

For a prolongation to the right by a simple chain of length 1
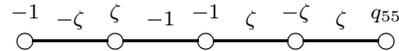
\begin{gather*}
%
%
\rule[-4\unitlength]{0pt}{5\unitlength}
\begin{picture}(50,4)(0,3)
\put(1,1){\circle{2}}
\put(2,1){\line(1,0){10}}
\put(13,1){\circle{2}}
\put(14,1){\line(1,0){10}}
\put(25,1){\circle{2}}
\put(26,1){\line(1,0){10}}
\put(37,1){\circle{2}}
\put(38,1){\line(1,0){10}}
\put(49,1){\circle{2}}
\put(1,4){\makebox[0pt]{\scriptsize $-1$}}
\put(7,3){\makebox[0pt]{\scriptsize $-\zeta $}}
\put(13,4){\makebox[0pt]{\scriptsize $\zeta $}}
\put(19,3){\makebox[0pt]{\scriptsize $-1$}}
\put(25,4){\makebox[0pt]{\scriptsize $-1$}}
\put(31,3){\makebox[0pt]{\scriptsize $\zeta $}}
\put(37,4){\makebox[0pt]{\scriptsize $-\zeta $}}
\put(43,3){\makebox[0pt]{\scriptsize $\zeta $}}
\put(49,4){\makebox[0pt]{\scriptsize $q_{55}$}}
\end{picture}
\end{gather*}
$\roots (\chi ;\Ndbasis _4+\Ndbasis _5,\Ndbasis _3+\Ndbasis _4,
\Ndbasis _2+\Ndbasis _3,\Ndbasis _1+\Ndbasis _2)$ has
the generalized Dynkin diagram
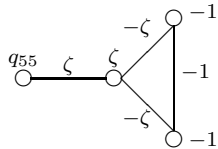
\begin{align*}
\Drightofway{}{$q_{55}$}{$\zeta $}%
{$\zeta $}{$-\zeta $}{$-\zeta $}{$-1$}{$-1$}{$-1$}
\end{align*}
which is again a contradiction to Theorem~\ref{t-classrank4}.
 
\textit{Step~2.5. Prolongations of diagram 7 in row~22 of Table~\ref{a-r4}.}
In a prolongation to the left by a simple chain of length 1
apply $s_{\Ndbasis _2,E}$, and
in a prolongation to the right by a simple chain
apply $s_{\Ndbasis _1,E}$. Then one obtains a diagram which was already
considered in Step~2.4.

\textit{Step~3.} Assume that
$(q_{11}+1)(q_{11}q_{12}q_{21}-1)=0$,
$(q_{dd}+1)(q_{dd}q_{d,d-1}q_{d-1,d}-1)=0$,
and there exists exactly
one integer $i$ with $1<i<d$ such that
$q_{ii}^2q_{i,i-1}q_{i-1,i}q_{i,i+1}q_{i+1,i}\not=1$.
Again it is sufficient to consider prolongations of those labeled
path graphs by simple chains which appear in Table~\ref{a-r4}
and satisfy equations
$(q_{11}+1)(q_{11}q_{12}q_{21}-1)=0$ and
$(q_{44}+1)(q_{44}q_{34}q_{43}-1)=0$.

\textit{Step~3.1. Prolongations of the diagram in row~3
of Table~\ref{a-r4}.}
To the left: this prolongation appears in row~7 or row~9 of Table~\ref{a-rb4}.
To the right: one gets the diagram
\begin{gather*}
%
%
\rule[-4\unitlength]{0pt}{5\unitlength}
\begin{picture}(50,4)(0,3)
\put(1,1){\circle{2}}
\put(2,1){\line(1,0){10}}
\put(13,1){\circle{2}}
\put(14,1){\line(1,0){10}}
\put(25,1){\circle{2}}
\put(26,1){\line(1,0){10}}
\put(37,1){\circle{2}}
\put(38,1){\line(1,0){10}}
\put(49,1){\circle{2}}
\put(1,4){\makebox[0pt]{\scriptsize $q$}}
\put(7,3){\makebox[0pt]{\scriptsize $q^{-1}$}}
\put(13,4){\makebox[0pt]{\scriptsize $q$}}
\put(19,3){\makebox[0pt]{\scriptsize $q^{-1}$}}
\put(25,4){\makebox[0pt]{\scriptsize $q$}}
\put(31,3){\makebox[0pt]{\scriptsize $q^{-2}$}}
\put(37,4){\makebox[0pt]{\scriptsize $q^2$}}
\put(43,3){\makebox[0pt]{\scriptsize $q^{-2}$}}
\put(49,4){\makebox[0pt]{\scriptsize $q_{55}$}}
\end{picture}
\end{gather*}
where $q_{55}\in \{q^2,-1\}$ and $q^2\not=1$.
If $q_{55}=q^2$ then $(\chi ,E)$ is of infinite Cartan type.
If $q_{55}=-1$ then by Theorem~\ref{t-classrank4}
$\roots (\chi ;\Ndbasis _2,\Ndbasis _3,\Ndbasis _4,\Ndbasis _5)$
is finite if and only if $q^2=-1$ in which case again
equation $q_{55}=q^2$ holds. 

\textit{Step~3.2. Prolongations of the diagram in row~4 of Table~\ref{a-r4}.}
To the left: one obtains the diagram
\begin{gather*}
%
%
\rule[-4\unitlength]{0pt}{5\unitlength}
\begin{picture}(50,4)(0,3)
\put(1,1){\circle{2}}
\put(2,1){\line(1,0){10}}
\put(13,1){\circle{2}}
\put(14,1){\line(1,0){10}}
\put(25,1){\circle{2}}
\put(26,1){\line(1,0){10}}
\put(37,1){\circle{2}}
\put(38,1){\line(1,0){10}}
\put(49,1){\circle{2}}
\put(1,4){\makebox[0pt]{\scriptsize $q_{11}$}}
\put(7,3){\makebox[0pt]{\scriptsize $q^{-2}$}}
\put(13,4){\makebox[0pt]{\scriptsize $q^2$}}
\put(19,3){\makebox[0pt]{\scriptsize $q^{-2}$}}
\put(25,4){\makebox[0pt]{\scriptsize $q^2$}}
\put(31,3){\makebox[0pt]{\scriptsize $q^{-2}$}}
\put(37,4){\makebox[0pt]{\scriptsize $q$}}
\put(43,3){\makebox[0pt]{\scriptsize $q^{-1}$}}
\put(49,4){\makebox[0pt]{\scriptsize $q$}}
\end{picture}
\end{gather*}
where $q_{11}\in \{q^2,-1\}$ and $q^2\not=1$.
In this case $\roots (\chi ;\Ndbasis _1+\Ndbasis _2,\Ndbasis _2+\Ndbasis _3,
\Ndbasis _3+\Ndbasis _4,\Ndbasis _4+\Ndbasis _5)$
has generalized Dynkin diagram
\begin{gather*}
%
%
\Dchainfour{}{$q_{11}$}{$q^{-2}$}{$q$}{$q^{-1}$}{$q$}{$q^{-2}$}{$q^2$}
\end{gather*}
which contradicts Theorem~\ref{t-classrank4}.

To the right: if in the diagram
\begin{gather}\label{eq-4-1-r}
%
%
\rule[-4\unitlength]{0pt}{5\unitlength}
\begin{picture}(50,4)(0,3)
\put(1,1){\circle{2}}
\put(2,1){\line(1,0){10}}
\put(13,1){\circle{2}}
\put(14,1){\line(1,0){10}}
\put(25,1){\circle{2}}
\put(26,1){\line(1,0){10}}
\put(37,1){\circle{2}}
\put(38,1){\line(1,0){10}}
\put(49,1){\circle{2}}
\put(1,4){\makebox[0pt]{\scriptsize $q^2$}}
\put(7,3){\makebox[0pt]{\scriptsize $q^{-2}$}}
\put(13,4){\makebox[0pt]{\scriptsize $q^2$}}
\put(19,3){\makebox[0pt]{\scriptsize $q^{-2}$}}
\put(25,4){\makebox[0pt]{\scriptsize $q$}}
\put(31,3){\makebox[0pt]{\scriptsize $q^{-1}$}}
\put(37,4){\makebox[0pt]{\scriptsize $q$}}
\put(43,3){\makebox[0pt]{\scriptsize $q^{-1}$}}
\put(49,4){\makebox[0pt]{\scriptsize $q_{55}$}}
\end{picture}
\end{gather}
where $(q_{55}+1)(q_{55}-q)=0$, equation $q_{55}=q$ holds,
then $(\roots, \chi ,E)$ is of infinite Cartan type.
On the other hand, if $q_{55}=-1$ then Lemma~\ref{l-5chainmid}
implies that $q\in R_3$. In this case $\gDd _{\chi ,E}$ is the last graph
in row~13 of Table~\ref{a-rb4}.

By the first part of Step~3.2 there exists no prolongation 
of (\ref{eq-4-1-r}) to the left. A prolongation of length 1 to the right
has to be of the form
\begin{gather*}
%
%
\rule[-4\unitlength]{0pt}{5\unitlength}
\begin{picture}(62,4)(0,3)
\put(1,1){\circle{2}}
\put(2,1){\line(1,0){10}}
\put(13,1){\circle{2}}
\put(14,1){\line(1,0){10}}
\put(25,1){\circle{2}}
\put(26,1){\line(1,0){10}}
\put(37,1){\circle{2}}
\put(38,1){\line(1,0){10}}
\put(49,1){\circle{2}}
\put(50,1){\line(1,0){10}}
\put(61,1){\circle{2}}
\put(1,4){\makebox[0pt]{\scriptsize $q^{-1}$}}
\put(7,3){\makebox[0pt]{\scriptsize $q$}}
\put(13,4){\makebox[0pt]{\scriptsize $q^{-1}$}}
\put(19,3){\makebox[0pt]{\scriptsize $q$}}
\put(25,4){\makebox[0pt]{\scriptsize $q$}}
\put(31,3){\makebox[0pt]{\scriptsize $q^{-1}$}}
\put(37,4){\makebox[0pt]{\scriptsize $q$}}
\put(43,3){\makebox[0pt]{\scriptsize $q^{-1}$}}
\put(49,4){\makebox[0pt]{\scriptsize $-1$}}
\put(55,3){\makebox[0pt]{\scriptsize $q$}}
\put(61,4){\makebox[0pt]{\scriptsize $q^{-1}$}}
\end{picture}
\quad ,
\end{gather*}
where $q\in R_3$,
because of the finiteness of $\roots (\chi ;\Ndbasis _1,\Ndbasis _2,
\Ndbasis _3,\Ndbasis _4,\Ndbasis _5+\Ndbasis _6)$.
In this case however one obtains a contradiction to the finiteness of
$\roots (\chi ;\Ndbasis _6,\Ndbasis _5,
\Ndbasis _3+\Ndbasis _4,\Ndbasis _2+\Ndbasis _3,\Ndbasis _1+\Ndbasis _2)$
with generalized Dynkin diagram
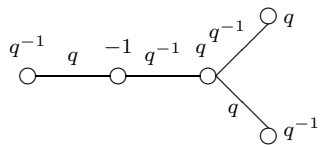
\begin{align*}
\rule[-10\unitlength]{0pt}{12\unitlength}
\begin{picture}(36,14)(0,9)
\put(1,13){\circle{2}}
\put(2,13){\line(1,0){10}}
\put(13,13){\circle{2}}
\put(14,13){\line(1,0){10}}
\put(25,13){\circle{2}}
\put(26,13){\line(1,1){7}}
\put(26,13){\line(1,-1){7}}
\put(33,21){\circle{2}}
\put(33,5){\circle{2}}
\put(1,16){\makebox[0pt]{\scriptsize $q^{-1}$}}
\put(7,15){\makebox[0pt]{\scriptsize $q$}}
\put(13,16){\makebox[0pt]{\scriptsize $-1$}}
\put(19,15){\makebox[0pt]{\scriptsize $q^{-1}$}}
\put(24,16){\makebox[0pt]{\scriptsize $q$}}
\put(30,18){\makebox[0pt][r]{\scriptsize $q^{-1}$}}
\put(29,8){\makebox[0pt][r]{\scriptsize $q$}}
\put(35,20){\makebox[0pt][l]{\scriptsize $q$}}
\put(35,5){\makebox[0pt][l]{\scriptsize $q^{-1}$}}
\end{picture}
\end{align*}
and to Theorem~\ref{t-classrank>4} with $d=5$.

\textit{Step~3.3. Prolongations of the labeled path graphs in row~8
of Table~\ref{a-r4}.}
Note that the first and second diagrams in row~8 are Weyl equivalent
via $s_{\Ndbasis _1,E}$, and the second and third are Weyl equivalent via
$s_{\Ndbasis _2,E}$. Therefore any prolongation of one of these generalized
Dynkin diagrams is Weyl equivalent to a prolongation of another one.
Thus it is sufficient to consider the first diagram in row~8.

Prolongations to the left:
these appear in rows~9 and 10 of Table~\ref{a-rb4}.

Prolongations to the right: one obtains the diagram
\begin{gather}\label{eq-8-1-r}
%
%
\rule[-4\unitlength]{0pt}{5\unitlength}
\begin{picture}(50,4)(0,3)
\put(1,1){\circle{2}}
\put(2,1){\line(1,0){10}}
\put(13,1){\circle{2}}
\put(14,1){\line(1,0){10}}
\put(25,1){\circle{2}}
\put(26,1){\line(1,0){10}}
\put(37,1){\circle{2}}
\put(38,1){\line(1,0){10}}
\put(49,1){\circle{2}}
\put(1,4){\makebox[0pt]{\scriptsize $-1$}}
\put(7,3){\makebox[0pt]{\scriptsize $q^{-1}$}}
\put(13,4){\makebox[0pt]{\scriptsize $q$}}
\put(19,3){\makebox[0pt]{\scriptsize $q^{-1}$}}
\put(25,4){\makebox[0pt]{\scriptsize $q$}}
\put(31,3){\makebox[0pt]{\scriptsize $q^{-2}$}}
\put(37,4){\makebox[0pt]{\scriptsize $q^2$}}
\put(43,3){\makebox[0pt]{\scriptsize $q^{-2}$}}
\put(49,4){\makebox[0pt]{\scriptsize $q_{55}$}}
\end{picture}
\end{gather}
where $(q_{55}+1)(q_{55}-q^2)=0$, and $q^2\not=1$.
Then Lemma~\ref{l-5chainmid}
implies that $q\in R_3$ and $q_{55}=q^{-1}$.
Thus diagram (\ref{eq-8-1-r}) and its
prolongations were already considered in Step~3.2 below graph (\ref{eq-4-1-r}).

\textit{Step~3.4. Prolongations of the labeled path graphs
in row~9 of Table~\ref{a-r4}.}
As in Step~3 Weyl equivalence allows to reduce to the
consideration of the second and fifth diagrams. Suppose that
$q^2,q^3\not=1$.

Prolongations to the left: assume that $(q_{11}+1)(q_{11}-q^2)=0$.
By Weyl equivalence one gets
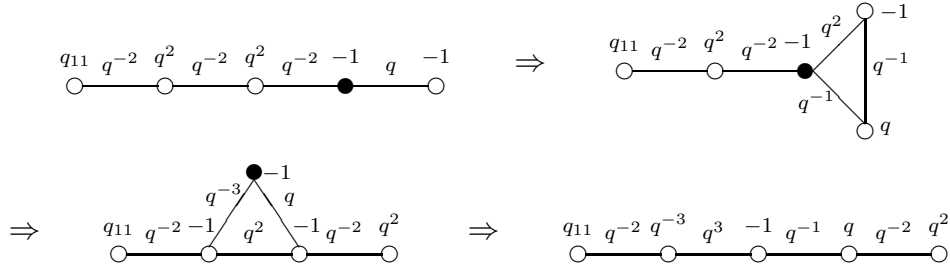
\begin{gather*}\label{eq-9-2-l}
%
%
\rule[-4\unitlength]{0pt}{5\unitlength}
\begin{picture}(50,4)(0,3)
\put(1,1){\circle{2}}
\put(2,1){\line(1,0){10}}
\put(13,1){\circle{2}}
\put(14,1){\line(1,0){10}}
\put(25,1){\circle{2}}
\put(26,1){\line(1,0){10}}
\put(37,1){\circle*{2}}
\put(38,1){\line(1,0){10}}
\put(49,1){\circle{2}}
\put(1,4){\makebox[0pt]{\scriptsize $q_{11}$}}
\put(7,3){\makebox[0pt]{\scriptsize $q^{-2}$}}
\put(13,4){\makebox[0pt]{\scriptsize $q^2$}}
\put(19,3){\makebox[0pt]{\scriptsize $q^{-2}$}}
\put(25,4){\makebox[0pt]{\scriptsize $q^2$}}
\put(31,3){\makebox[0pt]{\scriptsize $q^{-2}$}}
\put(37,4){\makebox[0pt]{\scriptsize $-1$}}
\put(43,3){\makebox[0pt]{\scriptsize $q$}}
\put(49,4){\makebox[0pt]{\scriptsize $-1$}}
\end{picture}
\qquad \Rightarrow \qquad
\rule[-10\unitlength]{0pt}{12\unitlength}
\begin{picture}(36,10)(0,9)
\put(1,9){\circle{2}}
\put(2,9){\line(1,0){10}}
\put(13,9){\circle{2}}
\put(14,9){\line(1,0){10}}
\put(25,9){\circle*{2}}
\put(26,9){\line(1,1){7}}
\put(26,9){\line(1,-1){7}}
\put(33,17){\circle{2}}
\put(33,1){\circle{2}}
\put(33,2){\line(0,1){14}}
\put(1,12){\makebox[0pt]{\scriptsize $q_{11}$}}
\put(7,11){\makebox[0pt]{\scriptsize $q^{-2}$}}
\put(13,12){\makebox[0pt]{\scriptsize $q^2$}}
\put(19,11){\makebox[0pt]{\scriptsize $q^{-2}$}}
\put(24,12){\makebox[0pt]{\scriptsize $-1$}}
\put(30,14){\makebox[0pt][r]{\scriptsize $q^2$}}
\put(29,4){\makebox[0pt][r]{\scriptsize $q^{-1}$}}
\put(35,16){\makebox[0pt][l]{\scriptsize $-1$}}
\put(34,9){\makebox[0pt][l]{\scriptsize $q^{-1}$}}
\put(35,1){\makebox[0pt][l]{\scriptsize $q$}}
\end{picture}
\\
\Rightarrow \qquad
\rule[-4\unitlength]{0pt}{5\unitlength}
\begin{picture}(38,11)(0,3)
\put(1,1){\circle{2}}
\put(2,1){\line(1,0){10}}
\put(13,1){\circle{2}}
\put(13,2){\line(2,3){6}}
\put(14,1){\line(1,0){10}}
\put(25,1){\circle{2}}
\put(25,2){\line(-2,3){6}}
\put(19,12){\circle*{2}}
\put(26,1){\line(1,0){10}}
\put(37,1){\circle{2}}
\put(1,4){\makebox[0pt]{\scriptsize $q_{11}$}}
\put(7,3){\makebox[0pt]{\scriptsize $q^{-2}$}}
\put(12,4){\makebox[0pt]{\scriptsize $-1$}}
\put(19,3){\makebox[0pt]{\scriptsize $q^2$}}
\put(26,4){\makebox[0pt]{\scriptsize $-1$}}
\put(31,3){\makebox[0pt]{\scriptsize $q^{-2}$}}
\put(37,4){\makebox[0pt]{\scriptsize $q^2$}}
\put(15,8){\makebox[0pt]{\scriptsize $q^{-3}$}}
\put(23,8){\makebox[0pt][l]{\scriptsize $q$}}
\put(22,11){\makebox[0pt]{\scriptsize $-1$}}
\end{picture}
\qquad \Rightarrow \qquad
%
%
\rule[-4\unitlength]{0pt}{5\unitlength}
\begin{picture}(50,4)(0,3)
\put(1,1){\circle{2}}
\put(2,1){\line(1,0){10}}
\put(13,1){\circle{2}}
\put(14,1){\line(1,0){10}}
\put(25,1){\circle{2}}
\put(26,1){\line(1,0){10}}
\put(37,1){\circle{2}}
\put(38,1){\line(1,0){10}}
\put(49,1){\circle{2}}
\put(1,4){\makebox[0pt]{\scriptsize $q_{11}$}}
\put(7,3){\makebox[0pt]{\scriptsize $q^{-2}$}}
\put(13,4){\makebox[0pt]{\scriptsize $q^{-3}$}}
\put(19,3){\makebox[0pt]{\scriptsize $q^3$}}
\put(25,4){\makebox[0pt]{\scriptsize $-1$}}
\put(31,3){\makebox[0pt]{\scriptsize $q^{-1}$}}
\put(37,4){\makebox[0pt]{\scriptsize $q$}}
\put(43,3){\makebox[0pt]{\scriptsize $q^{-2}$}}
\put(49,4){\makebox[0pt]{\scriptsize $q^2$}}
\end{picture}
\end{gather*}
By Lemma~\ref{l-5chain} this implies that $q\in R_5$. This diagram and its
prolongations to the left were already considered in Step~2.1.

Similarly, for a prolongation of the fifth diagram to the left
one obtains
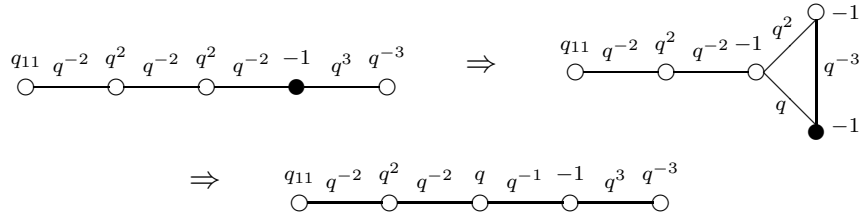
\begin{gather*}\label{eq-9-5-l}
%
%
\rule[-4\unitlength]{0pt}{5\unitlength}
\begin{picture}(50,4)(0,3)
\put(1,1){\circle{2}}
\put(2,1){\line(1,0){10}}
\put(13,1){\circle{2}}
\put(14,1){\line(1,0){10}}
\put(25,1){\circle{2}}
\put(26,1){\line(1,0){10}}
\put(37,1){\circle*{2}}
\put(38,1){\line(1,0){10}}
\put(49,1){\circle{2}}
\put(1,4){\makebox[0pt]{\scriptsize $q_{11}$}}
\put(7,3){\makebox[0pt]{\scriptsize $q^{-2}$}}
\put(13,4){\makebox[0pt]{\scriptsize $q^2$}}
\put(19,3){\makebox[0pt]{\scriptsize $q^{-2}$}}
\put(25,4){\makebox[0pt]{\scriptsize $q^2$}}
\put(31,3){\makebox[0pt]{\scriptsize $q^{-2}$}}
\put(37,4){\makebox[0pt]{\scriptsize $-1$}}
\put(43,3){\makebox[0pt]{\scriptsize $q^3$}}
\put(49,4){\makebox[0pt]{\scriptsize $q^{-3}$}}
\end{picture}
\qquad \Rightarrow \qquad
\rule[-10\unitlength]{0pt}{12\unitlength}
\begin{picture}(36,10)(0,9)
\put(1,9){\circle{2}}
\put(2,9){\line(1,0){10}}
\put(13,9){\circle{2}}
\put(14,9){\line(1,0){10}}
\put(25,9){\circle{2}}
\put(26,9){\line(1,1){7}}
\put(26,9){\line(1,-1){7}}
\put(33,17){\circle{2}}
\put(33,1){\circle*{2}}
\put(33,2){\line(0,1){14}}
\put(1,12){\makebox[0pt]{\scriptsize $q_{11}$}}
\put(7,11){\makebox[0pt]{\scriptsize $q^{-2}$}}
\put(13,12){\makebox[0pt]{\scriptsize $q^2$}}
\put(19,11){\makebox[0pt]{\scriptsize $q^{-2}$}}
\put(24,12){\makebox[0pt]{\scriptsize $-1$}}
\put(30,14){\makebox[0pt][r]{\scriptsize $q^2$}}
\put(29,4){\makebox[0pt][r]{\scriptsize $q$}}
\put(35,16){\makebox[0pt][l]{\scriptsize $-1$}}
\put(34,9){\makebox[0pt][l]{\scriptsize $q^{-3}$}}
\put(35,1){\makebox[0pt][l]{\scriptsize $-1$}}
\end{picture}
\\
\Rightarrow \qquad
%
%
\rule[-4\unitlength]{0pt}{5\unitlength}
\begin{picture}(50,4)(0,3)
\put(1,1){\circle{2}}
\put(2,1){\line(1,0){10}}
\put(13,1){\circle{2}}
\put(14,1){\line(1,0){10}}
\put(25,1){\circle{2}}
\put(26,1){\line(1,0){10}}
\put(37,1){\circle{2}}
\put(38,1){\line(1,0){10}}
\put(49,1){\circle{2}}
\put(1,4){\makebox[0pt]{\scriptsize $q_{11}$}}
\put(7,3){\makebox[0pt]{\scriptsize $q^{-2}$}}
\put(13,4){\makebox[0pt]{\scriptsize $q^2$}}
\put(19,3){\makebox[0pt]{\scriptsize $q^{-2}$}}
\put(25,4){\makebox[0pt]{\scriptsize $q$}}
\put(31,3){\makebox[0pt]{\scriptsize $q^{-1}$}}
\put(37,4){\makebox[0pt]{\scriptsize $-1$}}
\put(43,3){\makebox[0pt]{\scriptsize $q^3$}}
\put(49,4){\makebox[0pt]{\scriptsize $q^{-3}$}}
\end{picture}
\end{gather*}
In Step~2.1 it was shown that the latter diagram does not correspond to
an arithmetic root system.

Prolongations to the right: if $\gDd _{\chi ,E}$ is of the form
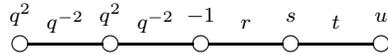
\begin{gather*}\label{eq-9-2-r}
%
%
\rule[-4\unitlength]{0pt}{5\unitlength}
\begin{picture}(50,4)(0,3)
\put(1,1){\circle{2}}
\put(2,1){\line(1,0){10}}
\put(13,1){\circle{2}}
\put(14,1){\line(1,0){10}}
\put(25,1){\circle{2}}
\put(26,1){\line(1,0){10}}
\put(37,1){\circle{2}}
\put(38,1){\line(1,0){10}}
\put(49,1){\circle{2}}
\put(1,4){\makebox[0pt]{\scriptsize $q^2$}}
\put(7,3){\makebox[0pt]{\scriptsize $q^{-2}$}}
\put(13,4){\makebox[0pt]{\scriptsize $q^2$}}
\put(19,3){\makebox[0pt]{\scriptsize $q^{-2}$}}
\put(25,4){\makebox[0pt]{\scriptsize $-1$}}
\put(31,3){\makebox[0pt]{\scriptsize $r$}}
\put(37,4){\makebox[0pt]{\scriptsize $s$}}
\put(43,3){\makebox[0pt]{\scriptsize $t$}}
\put(49,4){\makebox[0pt]{\scriptsize $u$}}
\end{picture}
\end{gather*}
where $q^2\not=1$, $q^3\not=1$, and one of the systems of equations
\begin{align*}
 r-q=s+1=t-q^{-1}=(u+1)(u-q)=0,\\
 r-q^3=s-q^{-3}=t-q^3=(u+1)(uq^3-1)=0
\end{align*}
holds,
then Lemma~\ref{l-5chainmid}
implies that $\{r,r^{-1}\}=
\{q^2,q^{-2}\}=\{q^{-2}r,q^2r^{-1}\}$. However if $r=q$ then
$r\notin \{q^2,q^{-2}\}$ and if $r=q^3$ then $q^{-2}r\notin \{q^2,q^{-2}\}$.
Therefore there exist no prolongations to the right of the
diagrams 2 and 5 in row~9 of Table~\ref{a-r4}.

\textit{Step~3.5. Prolongations of the labeled path graphs
in row~12 of Table~\ref{a-r4}.}
Again Weyl equivalence allows to reduce to the consideration
of the first diagram. The prolongations to the left appear
in rows~9 and 10 of Table~\ref{a-rb4}.
A prolongation of length 1 to the right takes the form
\begin{gather}\label{eq-12-1-r}
%
%
\rule[-4\unitlength]{0pt}{5\unitlength}
\begin{picture}(50,4)(0,3)
\put(1,1){\circle{2}}
\put(2,1){\line(1,0){10}}
\put(13,1){\circle{2}}
\put(14,1){\line(1,0){10}}
\put(25,1){\circle{2}}
\put(26,1){\line(1,0){10}}
\put(37,1){\circle{2}}
\put(38,1){\line(1,0){10}}
\put(49,1){\circle{2}}
\put(1,4){\makebox[0pt]{\scriptsize $q^{-1}$}}
\put(7,3){\makebox[0pt]{\scriptsize $q$}}
\put(13,4){\makebox[0pt]{\scriptsize $-1$}}
\put(19,3){\makebox[0pt]{\scriptsize $q^{-1}$}}
\put(25,4){\makebox[0pt]{\scriptsize $q$}}
\put(31,3){\makebox[0pt]{\scriptsize $q^{-2}$}}
\put(37,4){\makebox[0pt]{\scriptsize $q^2$}}
\put(43,3){\makebox[0pt]{\scriptsize $q^{-2}$}}
\put(49,4){\makebox[0pt]{\scriptsize $q_{55}$}}
\end{picture}
\end{gather}
where $q^2\not=1$ and $(q_{55}+1)(q_{55}-q^2)=0$. Thus Lemma~\ref{l-5chainmid}
gives that $q\in R_3$ and $q_{55}=q^2$.
This diagram coincides with the last one in row~12 of Table~\ref{a-rb4}.

A prolongation of length one to the right of the diagram (\ref{eq-12-1-r}),
where $q\in R_3$ and $q_{55}=q^2$, has to be of the form
\begin{gather*}
%
%
\rule[-4\unitlength]{0pt}{5\unitlength}
\begin{picture}(62,4)(0,3)
\put(1,1){\circle{2}}
\put(2,1){\line(1,0){10}}
\put(13,1){\circle{2}}
\put(14,1){\line(1,0){10}}
\put(25,1){\circle{2}}
\put(26,1){\line(1,0){10}}
\put(37,1){\circle{2}}
\put(38,1){\line(1,0){10}}
\put(49,1){\circle{2}}
\put(50,1){\line(1,0){10}}
\put(61,1){\circle{2}}
\put(1,4){\makebox[0pt]{\scriptsize $q^{-1}$}}
\put(7,3){\makebox[0pt]{\scriptsize $q$}}
\put(13,4){\makebox[0pt]{\scriptsize $-1$}}
\put(19,3){\makebox[0pt]{\scriptsize $q^{-1}$}}
\put(25,4){\makebox[0pt]{\scriptsize $q$}}
\put(31,3){\makebox[0pt]{\scriptsize $q$}}
\put(37,4){\makebox[0pt]{\scriptsize $q^{-1}$}}
\put(43,3){\makebox[0pt]{\scriptsize $q$}}
\put(49,4){\makebox[0pt]{\scriptsize $q^{-1}$}}
\put(55,3){\makebox[0pt]{\scriptsize $q$}}
\put(61,4){\makebox[0pt]{\scriptsize $q_{66}$}}
\end{picture}
\end{gather*}
where $q_{66}=q^{-1}$, since $\roots (\chi ;\Ndbasis _1,\Ndbasis _2,
\Ndbasis _3,\Ndbasis _4+\Ndbasis _5,\Ndbasis _6)$ is finite.
In this case $\roots (\chi ;\Ndbasis _5+\Ndbasis _6,\Ndbasis _4+\Ndbasis _5,
\Ndbasis _3+\Ndbasis _4,\Ndbasis _2+\Ndbasis _3,\Ndbasis _1)$
has the generalized Dynkin diagram
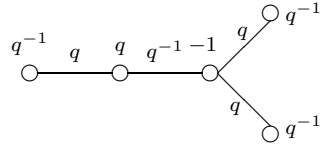
\begin{gather*}
\rule[-10\unitlength]{0pt}{12\unitlength}
\begin{picture}(36,14)(0,9)
\put(1,13){\circle{2}}
\put(2,13){\line(1,0){10}}
\put(13,13){\circle{2}}
\put(14,13){\line(1,0){10}}
\put(25,13){\circle{2}}
\put(26,13){\line(1,1){7}}
\put(26,13){\line(1,-1){7}}
\put(33,21){\circle{2}}
\put(33,5){\circle{2}}
\put(1,16){\makebox[0pt]{\scriptsize $q^{-1}$}}
\put(7,15){\makebox[0pt]{\scriptsize $q$}}
\put(13,16){\makebox[0pt]{\scriptsize $q$}}
\put(19,15){\makebox[0pt]{\scriptsize $q^{-1}$}}
\put(24,16){\makebox[0pt]{\scriptsize $-1$}}
\put(30,18){\makebox[0pt][r]{\scriptsize $q$}}
\put(29,8){\makebox[0pt][r]{\scriptsize $q$}}
\put(35,20){\makebox[0pt][l]{\scriptsize $q^{-1}$}}
\put(35,5){\makebox[0pt][l]{\scriptsize $q^{-1}$}}
\end{picture}
\end{gather*}
which is a contradiction to Theorem~\ref{t-classrank>4} for $d=5$.

It remains to consider prolongations of (\ref{eq-12-1-r}) to the left,
where $q\in R_3$ and $q_{55}=q^{-1}$. Assume that $\gDd _{\chi ,E}$
is of the form
\begin{gather}\label{eq-12-1-rl}
%
%
\rule[-4\unitlength]{0pt}{5\unitlength}
\begin{picture}(62,4)(0,3)
\put(1,1){\circle{2}}
\put(2,1){\line(1,0){10}}
\put(13,1){\circle{2}}
\put(14,1){\line(1,0){10}}
\put(25,1){\circle{2}}
\put(26,1){\line(1,0){10}}
\put(37,1){\circle{2}}
\put(38,1){\line(1,0){10}}
\put(49,1){\circle{2}}
\put(50,1){\line(1,0){10}}
\put(61,1){\circle{2}}
\put(1,4){\makebox[0pt]{\scriptsize $q_{11}$}}
\put(7,3){\makebox[0pt]{\scriptsize $q$}}
\put(13,4){\makebox[0pt]{\scriptsize $q^{-1}$}}
\put(19,3){\makebox[0pt]{\scriptsize $q$}}
\put(25,4){\makebox[0pt]{\scriptsize $-1$}}
\put(31,3){\makebox[0pt]{\scriptsize $q^{-1}$}}
\put(37,4){\makebox[0pt]{\scriptsize $q$}}
\put(43,3){\makebox[0pt]{\scriptsize $q$}}
\put(49,4){\makebox[0pt]{\scriptsize $q^{-1}$}}
\put(55,3){\makebox[0pt]{\scriptsize $q$}}
\put(61,4){\makebox[0pt]{\scriptsize $q^{-1}$}}
\end{picture}
\end{gather}
where $q\in R_3$ and $(q_{11}+1)(q_{11}q-1)=0$.
Then $\roots (\chi ;\Ndbasis _1,\Ndbasis _2,
\Ndbasis _3+\Ndbasis _4,\Ndbasis _4+\Ndbasis _5,
\Ndbasis _5+\Ndbasis _6)$ has generalized Dynkin diagram
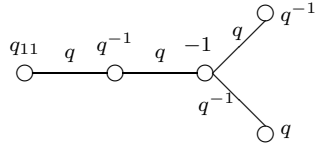
\begin{gather*}
\rule[-10\unitlength]{0pt}{12\unitlength}
\begin{picture}(36,14)(0,9)
\put(1,13){\circle{2}}
\put(2,13){\line(1,0){10}}
\put(13,13){\circle{2}}
\put(14,13){\line(1,0){10}}
\put(25,13){\circle{2}}
\put(26,13){\line(1,1){7}}
\put(26,13){\line(1,-1){7}}
\put(33,21){\circle{2}}
\put(33,5){\circle{2}}
\put(1,16){\makebox[0pt]{\scriptsize $q_{11}$}}
\put(7,15){\makebox[0pt]{\scriptsize $q$}}
\put(13,16){\makebox[0pt]{\scriptsize $q^{-1}$}}
\put(19,15){\makebox[0pt]{\scriptsize $q$}}
\put(24,16){\makebox[0pt]{\scriptsize $-1$}}
\put(30,18){\makebox[0pt][r]{\scriptsize $q$}}
\put(29,8){\makebox[0pt][r]{\scriptsize $q^{-1}$}}
\put(35,20){\makebox[0pt][l]{\scriptsize $q^{-1}$}}
\put(35,5){\makebox[0pt][l]{\scriptsize $q$}}
\end{picture}
\end{gather*}
and hence Theorem~\ref{t-classrank>4} for $d=5$ implies that $q_{11}=q^{-1}$.
Thus (\ref{eq-12-1-rl}) is the fourth last diagram in row~18 of
Table~\ref{a-rb4}.
For a prolongation of length one to the left of (\ref{eq-12-1-rl})
one again has the only possibility 
\begin{align*}
\rule[-4\unitlength]{0pt}{5\unitlength}
\begin{picture}(62,4)(0,3)
\put(1,1){\circle{2}}
\put(2,1){\line(1,0){8}}
\put(11,1){\circle{2}}
\put(12,1){\line(1,0){8}}
\put(21,1){\circle{2}}
\put(22,1){\line(1,0){8}}
\put(31,1){\circle{2}}
\put(32,1){\line(1,0){8}}
\put(41,1){\circle{2}}
\put(42,1){\line(1,0){8}}
\put(51,1){\circle{2}}
\put(52,1){\line(1,0){8}}
\put(61,1){\circle{2}}
\put(1,4){\makebox[0pt]{\scriptsize $q^{-1}$}}
\put(6,3){\makebox[0pt]{\scriptsize $q$}}
\put(11,4){\makebox[0pt]{\scriptsize $q^{-1}$}}
\put(16,3){\makebox[0pt]{\scriptsize $q$}}
\put(21,4){\makebox[0pt]{\scriptsize $q^{-1}$}}
\put(26,3){\makebox[0pt]{\scriptsize $q$}}
\put(31,4){\makebox[0pt]{\scriptsize $-1$}}
\put(36,3){\makebox[0pt]{\scriptsize $q^{-1}$}}
\put(41,4){\makebox[0pt]{\scriptsize $q$}}
\put(46,3){\makebox[0pt]{\scriptsize $q$}}
\put(51,4){\makebox[0pt]{\scriptsize $q^{-1}$}}
\put(56,3){\makebox[0pt]{\scriptsize $q$}}
\put(61,4){\makebox[0pt]{\scriptsize $q^{-1}$}}
\end{picture}
\end{align*}
since $\roots (\chi ;\Ndbasis _1,\Ndbasis _2+\Ndbasis _3,
\Ndbasis _4,\Ndbasis _5,\Ndbasis _6,\Ndbasis _7)$
has to be finite.
However in this case $\roots (\chi ;
\Ndbasis _1+\Ndbasis _2,\Ndbasis _2+\Ndbasis _3,\Ndbasis _3+\Ndbasis _4,
\Ndbasis _4+\Ndbasis _5,\Ndbasis _5+\Ndbasis _6,\Ndbasis _7)$
has to be finite but it has generalized Dynkin diagram
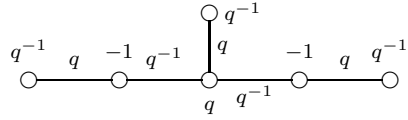
\begin{align*}
 &
\rule[-8\unitlength]{0pt}{9\unitlength}
\begin{picture}(41,8)(0,7)
\put(1,4){\circle{2}}
\put(2,4){\line(1,0){10}}
\put(13,4){\circle{2}}
\put(14,4){\line(1,0){10}}
\put(25,4){\circle{2}}
\put(26,4){\line(1,0){10}}
\put(37,4){\circle{2}}
\put(38,4){\line(1,0){10}}
\put(49,4){\circle{2}}
\put(25,5){\line(0,1){7}}
\put(25,13){\circle{2}}
\put(1,7){\makebox[0pt]{\scriptsize $q^{-1}$}}
\put(7,6){\makebox[0pt]{\scriptsize $q$}}
\put(13,7){\makebox[0pt]{\scriptsize $-1$}}
\put(19,6){\makebox[0pt]{\scriptsize $q^{-1}$}}
\put(25,0){\makebox[0pt]{\scriptsize $q$}}
\put(31,1){\makebox[0pt]{\scriptsize $q^{-1}$}}
\put(37,7){\makebox[0pt]{\scriptsize $-1$}}
\put(43,6){\makebox[0pt]{\scriptsize $q$}}
\put(49,7){\makebox[0pt]{\scriptsize $q^{-1}$}}
\put(26,8){\makebox[0pt][l]{\scriptsize $q$}}
\put(27,12){\makebox[0pt][l]{\scriptsize $q^{-1}$}}
\end{picture}
\end{align*}
which is a contradiction to Theorem~\ref{t-classrank>4} with $d=6$.

\textit{Step~3.6. Prolongations of the first diagram in row~13
of Table~\ref{a-r4}.}
The prolongations to the left appear in rows~9 and 10
of Table~\ref{a-rb4}. A prolongation of length 1 to the right takes the form
\begin{gather}\label{eq-13-1-r}
%
%
\rule[-4\unitlength]{0pt}{5\unitlength}
\begin{picture}(50,4)(0,3)
\put(1,1){\circle{2}}
\put(2,1){\line(1,0){10}}
\put(13,1){\circle{2}}
\put(14,1){\line(1,0){10}}
\put(25,1){\circle{2}}
\put(26,1){\line(1,0){10}}
\put(37,1){\circle{2}}
\put(38,1){\line(1,0){10}}
\put(49,1){\circle{2}}
\put(1,4){\makebox[0pt]{\scriptsize $q$}}
\put(7,3){\makebox[0pt]{\scriptsize $q^{-1}$}}
\put(13,4){\makebox[0pt]{\scriptsize $q$}}
\put(19,3){\makebox[0pt]{\scriptsize $q^{-1}$}}
\put(25,4){\makebox[0pt]{\scriptsize $-1$}}
\put(31,3){\makebox[0pt]{\scriptsize $q^2$}}
\put(37,4){\makebox[0pt]{\scriptsize $q^{-2}$}}
\put(43,3){\makebox[0pt]{\scriptsize $q^2$}}
\put(49,4){\makebox[0pt]{\scriptsize $q_{55}$}}
\end{picture}
\end{gather}
where $q^2\not=1$ and $(q_{55}+1)(q_{55}q^2-1)=0$. Thus Lemma~\ref{l-5chainmid}
gives that $q\in R_3$. If $q_{55}=-1$ then Weyl equivalence gives
\begin{gather}\label{eq-13-1-r-minus1}
%
%
\rule[-4\unitlength]{0pt}{5\unitlength}
\begin{picture}(50,4)(0,3)
\put(1,1){\circle{2}}
\put(2,1){\line(1,0){10}}
\put(13,1){\circle{2}}
\put(14,1){\line(1,0){10}}
\put(25,1){\circle{2}}
\put(26,1){\line(1,0){10}}
\put(37,1){\circle{2}}
\put(38,1){\line(1,0){10}}
\put(49,1){\circle*{2}}
\put(1,4){\makebox[0pt]{\scriptsize $q$}}
\put(7,3){\makebox[0pt]{\scriptsize $q^{-1}$}}
\put(13,4){\makebox[0pt]{\scriptsize $q$}}
\put(19,3){\makebox[0pt]{\scriptsize $q^{-1}$}}
\put(25,4){\makebox[0pt]{\scriptsize $-1$}}
\put(31,3){\makebox[0pt]{\scriptsize $q^{-1}$}}
\put(37,4){\makebox[0pt]{\scriptsize $q$}}
\put(43,3){\makebox[0pt]{\scriptsize $q^{-1}$}}
\put(49,4){\makebox[0pt]{\scriptsize $-1$}}
\end{picture}
\qquad \Rightarrow \qquad
%
%
\rule[-4\unitlength]{0pt}{5\unitlength}
\begin{picture}(50,4)(0,3)
\put(1,1){\circle{2}}
\put(2,1){\line(1,0){10}}
\put(13,1){\circle{2}}
\put(14,1){\line(1,0){10}}
\put(25,1){\circle{2}}
\put(26,1){\line(1,0){10}}
\put(37,1){\circle*{2}}
\put(38,1){\line(1,0){10}}
\put(49,1){\circle{2}}
\put(1,4){\makebox[0pt]{\scriptsize $q$}}
\put(7,3){\makebox[0pt]{\scriptsize $q^{-1}$}}
\put(13,4){\makebox[0pt]{\scriptsize $q$}}
\put(19,3){\makebox[0pt]{\scriptsize $q^{-1}$}}
\put(25,4){\makebox[0pt]{\scriptsize $-1$}}
\put(31,3){\makebox[0pt]{\scriptsize $q^{-1}$}}
\put(37,4){\makebox[0pt]{\scriptsize $-1$}}
\put(43,3){\makebox[0pt]{\scriptsize $q$}}
\put(49,4){\makebox[0pt]{\scriptsize $-1$}}
\end{picture}\\
\label{eq-13-1-r-minus1W}
\Rightarrow \qquad
%
%
\rule[-4\unitlength]{0pt}{5\unitlength}
\begin{picture}(50,4)(0,3)
\put(1,1){\circle{2}}
\put(2,1){\line(1,0){10}}
\put(13,1){\circle{2}}
\put(14,1){\line(1,0){10}}
\put(25,1){\circle{2}}
\put(26,1){\line(1,0){10}}
\put(37,1){\circle{2}}
\put(38,1){\line(1,0){10}}
\put(49,1){\circle{2}}
\put(1,4){\makebox[0pt]{\scriptsize $q$}}
\put(7,3){\makebox[0pt]{\scriptsize $q^{-1}$}}
\put(13,4){\makebox[0pt]{\scriptsize $q$}}
\put(19,3){\makebox[0pt]{\scriptsize $q^{-1}$}}
\put(25,4){\makebox[0pt]{\scriptsize $q^{-1}$}}
\put(31,3){\makebox[0pt]{\scriptsize $q$}}
\put(37,4){\makebox[0pt]{\scriptsize $-1$}}
\put(43,3){\makebox[0pt]{\scriptsize $q^{-1}$}}
\put(49,4){\makebox[0pt]{\scriptsize $q$}}
\end{picture}
\end{gather}
This diagram and all of its prolongations were already considered in
Step~3.5 and hence for $q_{55}=-1$
(\ref{eq-13-1-r}) and all of its prolongations are
Weyl equivalent to generalized Dynkin diagrams appearing in Step~3.5.

If $q_{55}=q\in R_3$ then (\ref{eq-13-1-r}) is the first diagram in row~11
of Table~\ref{a-rb4}. Because of its symmetry it is sufficient to consider ist
prolongations to the left and to both directions, respectively.
So assume next that $\gDd _{\chi ,E}$ takes the form
\begin{gather}\label{eq-13-1-rl}
%
%
\rule[-4\unitlength]{0pt}{5\unitlength}
\begin{picture}(62,4)(0,3)
\put(1,1){\circle{2}}
\put(2,1){\line(1,0){10}}
\put(13,1){\circle{2}}
\put(14,1){\line(1,0){10}}
\put(25,1){\circle{2}}
\put(26,1){\line(1,0){10}}
\put(37,1){\circle{2}}
\put(38,1){\line(1,0){10}}
\put(49,1){\circle{2}}
\put(50,1){\line(1,0){10}}
\put(61,1){\circle{2}}
\put(1,4){\makebox[0pt]{\scriptsize $q_{11}$}}
\put(7,3){\makebox[0pt]{\scriptsize $q^{-1}$}}
\put(13,4){\makebox[0pt]{\scriptsize $q$}}
\put(19,3){\makebox[0pt]{\scriptsize $q^{-1}$}}
\put(25,4){\makebox[0pt]{\scriptsize $q$}}
\put(31,3){\makebox[0pt]{\scriptsize $q^{-1}$}}
\put(37,4){\makebox[0pt]{\scriptsize $-1$}}
\put(43,3){\makebox[0pt]{\scriptsize $q^{-1}$}}
\put(49,4){\makebox[0pt]{\scriptsize $q$}}
\put(55,3){\makebox[0pt]{\scriptsize $q^{-1}$}}
\put(61,4){\makebox[0pt]{\scriptsize $q$}}
\end{picture}
\end{gather}
where $q\in R_3$ and $(q_{11}+1)(q_{11}-q)=0$. If $q_{11}=-1$ then
similarly to the transformations in (\ref{eq-13-1-r-minus1})
and (\ref{eq-13-1-r-minus1W}) one
obtains a prolongation
of (\ref{eq-13-1-r-minus1W}) which was already considered
in Step~3.5. On the other hand, if $q_{11}=q$ then $\gDd _{\chi ,E}$ appears
in row~17 of Table~\ref{a-rb4}.

The diagram (\ref{eq-13-1-rl}) has no prolongation to the right. Indeed,
if $\gDd _{\chi ,E}$ is of the form
\begin{align*}
\rule[-4\unitlength]{0pt}{5\unitlength}
\begin{picture}(62,4)(0,3)
\put(1,1){\circle{2}}
\put(2,1){\line(1,0){8}}
\put(11,1){\circle{2}}
\put(12,1){\line(1,0){8}}
\put(21,1){\circle{2}}
\put(22,1){\line(1,0){8}}
\put(31,1){\circle{2}}
\put(32,1){\line(1,0){8}}
\put(41,1){\circle{2}}
\put(42,1){\line(1,0){8}}
\put(51,1){\circle{2}}
\put(52,1){\line(1,0){8}}
\put(61,1){\circle{2}}
\put(1,4){\makebox[0pt]{\scriptsize $q$}}
\put(6,3){\makebox[0pt]{\scriptsize $q^{-1}$}}
\put(11,4){\makebox[0pt]{\scriptsize $q$}}
\put(16,3){\makebox[0pt]{\scriptsize $q^{-1}$}}
\put(21,4){\makebox[0pt]{\scriptsize $q$}}
\put(26,3){\makebox[0pt]{\scriptsize $q^{-1}$}}
\put(31,4){\makebox[0pt]{\scriptsize $-1$}}
\put(36,3){\makebox[0pt]{\scriptsize $q^{-1}$}}
\put(41,4){\makebox[0pt]{\scriptsize $q$}}
\put(46,3){\makebox[0pt]{\scriptsize $q^{-1}$}}
\put(51,4){\makebox[0pt]{\scriptsize $q$}}
\put(56,3){\makebox[0pt]{\scriptsize $q^{-1}$}}
\put(61,4){\makebox[0pt]{\scriptsize $q_{77}$}}
\end{picture}
\end{align*}
with $q\in R_3$ and $(q_{77}+1)(q_{77}-q)=0$ then the finiteness of
$\roots (\chi ;\Ndbasis _1+\Ndbasis _2,\Ndbasis _2+\Ndbasis _3,
\Ndbasis _3+\Ndbasis _4,\Ndbasis _4+\Ndbasis _5,\Ndbasis _5+\Ndbasis _6,
\Ndbasis _6+\Ndbasis _7)$
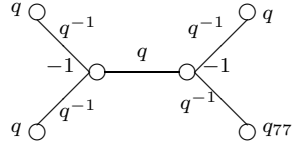
\begin{gather*}
\rule[-10\unitlength]{0pt}{12\unitlength}
\begin{picture}(36,14)(0,9)
\put(5,21){\circle{2}}
\put(5,5){\circle{2}}
\put(5,6){\line(1,1){7}}
\put(5,20){\line(1,-1){7}}
\put(13,13){\circle{2}}
\put(14,13){\line(1,0){10}}
\put(25,13){\circle{2}}
\put(26,13){\line(1,1){7}}
\put(26,13){\line(1,-1){7}}
\put(33,21){\circle{2}}
\put(33,5){\circle{2}}
\put(3,21){\makebox[0pt][r]{\scriptsize $q$}}
\put(3,5){\makebox[0pt][r]{\scriptsize $q$}}
\put(8,7){\makebox[0pt][l]{\scriptsize $q^{-1}$}}
\put(8,18){\makebox[0pt][l]{\scriptsize $q^{-1}$}}
\put(10,13){\makebox[0pt][r]{\scriptsize $-1$}}
\put(19,15){\makebox[0pt]{\scriptsize $q$}}
\put(27,13){\makebox[0pt][l]{\scriptsize $-1$}}
\put(30,18){\makebox[0pt][r]{\scriptsize $q^{-1}$}}
\put(29,8){\makebox[0pt][r]{\scriptsize $q^{-1}$}}
\put(35,20){\makebox[0pt][l]{\scriptsize $q$}}
\put(35,5){\makebox[0pt][l]{\scriptsize $q_{77}$}}
\end{picture}
\end{gather*}
is a contradiction to Lemma~\ref{l-nooctopus}. Thus it remains to
consider prolongations of (\ref{eq-13-1-rl}) to the left.

Suppose that $q\in R_3$, $(q_{11}+1)(q_{11}-q)=0$ and $\gDd _{\chi ,E}$
is of the form
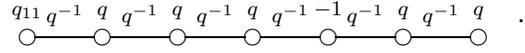
\begin{align*}
\rule[-4\unitlength]{0pt}{5\unitlength}
\begin{picture}(62,4)(0,3)
\put(1,1){\circle{2}}
\put(2,1){\line(1,0){8}}
\put(11,1){\circle{2}}
\put(12,1){\line(1,0){8}}
\put(21,1){\circle{2}}
\put(22,1){\line(1,0){8}}
\put(31,1){\circle{2}}
\put(32,1){\line(1,0){8}}
\put(41,1){\circle{2}}
\put(42,1){\line(1,0){8}}
\put(51,1){\circle{2}}
\put(52,1){\line(1,0){8}}
\put(61,1){\circle{2}}
\put(1,4){\makebox[0pt]{\scriptsize $q_{11}$}}
\put(6,3){\makebox[0pt]{\scriptsize $q^{-1}$}}
\put(11,4){\makebox[0pt]{\scriptsize $q$}}
\put(16,3){\makebox[0pt]{\scriptsize $q^{-1}$}}
\put(21,4){\makebox[0pt]{\scriptsize $q$}}
\put(26,3){\makebox[0pt]{\scriptsize $q^{-1}$}}
\put(31,4){\makebox[0pt]{\scriptsize $q$}}
\put(36,3){\makebox[0pt]{\scriptsize $q^{-1}$}}
\put(41,4){\makebox[0pt]{\scriptsize $-1$}}
\put(46,3){\makebox[0pt]{\scriptsize $q^{-1}$}}
\put(51,4){\makebox[0pt]{\scriptsize $q$}}
\put(56,3){\makebox[0pt]{\scriptsize $q^{-1}$}}
\put(61,4){\makebox[0pt]{\scriptsize $q$}}
\end{picture}\quad .
\end{align*}
If $q_{11}=-1$ then
similarly to the transformations in (\ref{eq-13-1-r-minus1})
and (\ref{eq-13-1-r-minus1W}) one obtains
a prolongation of (\ref{eq-13-1-r-minus1W})
which was already considered
in Step~3.5. On the other hand, if $q_{11}=q$ then $\gDd _{\chi ,E}$ appears
in row~21 of Table~\ref{a-rb4}.
For an additional prolongation to the left one can again assume
that $\gDd _{\chi ,E}$ is of the form
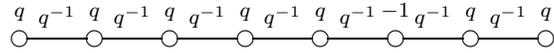
\begin{align*}
\rule[-4\unitlength]{0pt}{5\unitlength}
\begin{picture}(72,4)(0,3)
\put(1,1){\circle{2}}
\put(2,1){\line(1,0){8}}
\put(11,1){\circle{2}}
\put(12,1){\line(1,0){8}}
\put(21,1){\circle{2}}
\put(22,1){\line(1,0){8}}
\put(31,1){\circle{2}}
\put(32,1){\line(1,0){8}}
\put(41,1){\circle{2}}
\put(42,1){\line(1,0){8}}
\put(51,1){\circle{2}}
\put(52,1){\line(1,0){8}}
\put(61,1){\circle{2}}
\put(62,1){\line(1,0){8}}
\put(71,1){\circle{2}}
\put(1,4){\makebox[0pt]{\scriptsize $q$}}
\put(6,3){\makebox[0pt]{\scriptsize $q^{-1}$}}
\put(11,4){\makebox[0pt]{\scriptsize $q$}}
\put(16,3){\makebox[0pt]{\scriptsize $q^{-1}$}}
\put(21,4){\makebox[0pt]{\scriptsize $q$}}
\put(26,3){\makebox[0pt]{\scriptsize $q^{-1}$}}
\put(31,4){\makebox[0pt]{\scriptsize $q$}}
\put(36,3){\makebox[0pt]{\scriptsize $q^{-1}$}}
\put(41,4){\makebox[0pt]{\scriptsize $q$}}
\put(46,3){\makebox[0pt]{\scriptsize $q^{-1}$}}
\put(51,4){\makebox[0pt]{\scriptsize $-1$}}
\put(56,3){\makebox[0pt]{\scriptsize $q^{-1}$}}
\put(61,4){\makebox[0pt]{\scriptsize $q$}}
\put(66,3){\makebox[0pt]{\scriptsize $q^{-1}$}}
\put(71,4){\makebox[0pt]{\scriptsize $q$}}
\end{picture}
\end{align*}
where $q\in R_3$. In this case $\roots (\chi ;\Ndbasis _1,
\Ndbasis _2+\Ndbasis _3,\Ndbasis _3+\Ndbasis _4,\Ndbasis _4+\Ndbasis _5,
\Ndbasis _5+\Ndbasis _6,\Ndbasis _6+\Ndbasis _7,\Ndbasis _8)$
has generalized Dynkin diagram
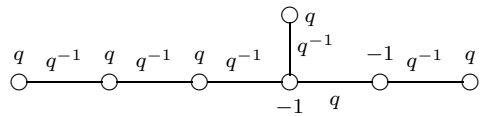
\begin{gather*}
\rule[-8\unitlength]{0pt}{9\unitlength}
\begin{picture}(62,8)(0,7)
\put(1,4){\circle{2}}
\put(2,4){\line(1,0){10}}
\put(13,4){\circle{2}}
\put(14,4){\line(1,0){10}}
\put(25,4){\circle{2}}
\put(26,4){\line(1,0){10}}
\put(37,4){\circle{2}}
\put(38,4){\line(1,0){10}}
\put(49,4){\circle{2}}
\put(50,4){\line(1,0){10}}
\put(61,4){\circle{2}}
\put(37,5){\line(0,1){7}}
\put(37,13){\circle{2}}
\put(1,7){\makebox[0pt]{\scriptsize $q$}}
\put(7,6){\makebox[0pt]{\scriptsize $q^{-1}$}}
\put(13,7){\makebox[0pt]{\scriptsize $q$}}
\put(19,6){\makebox[0pt]{\scriptsize $q^{-1}$}}
\put(25,7){\makebox[0pt]{\scriptsize $q$}}
\put(31,6){\makebox[0pt]{\scriptsize $q^{-1}$}}
\put(37,0){\makebox[0pt]{\scriptsize $-1$}}
\put(43,1){\makebox[0pt]{\scriptsize $q$}}
\put(49,7){\makebox[0pt]{\scriptsize $-1$}}
\put(55,6){\makebox[0pt]{\scriptsize $q^{-1}$}}
\put(61,7){\makebox[0pt]{\scriptsize $q$}}
\put(38,8){\makebox[0pt][l]{\scriptsize $q^{-1}$}}
\put(39,12){\makebox[0pt][l]{\scriptsize $q$}}
\end{picture}
\end{gather*}
which is a contradiction to Theorem~\ref{t-classrank>4} for $d=7$.

\textit{Step~3.7. Prolongations of the first and last diagrams in row~14
of Table~\ref{a-r4}.}
The last diagram in row~14 of Table~\ref{a-r4}
can be obtained from the first one
by replacing $q$ by $-q^{-1}$, and hence it is sufficient to consider
prolongations of the first graph.

By Lemma~\ref{l-5chainmid} for any prolongation to the right one has
to have the relation $\{-q,-q^{-1}\}=\{q,q^{-1}\}=\{-1\}$ which is
a contradiction to $q^2\not=1$. Moreover, since the first diagram
of row~14 is Weyl equivalent to the fourth one via
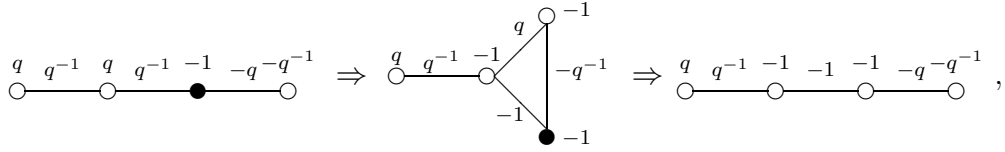
\begin{align*}
\Dchainfour{3}{$q$}{$q^{-1}$}{$q$}{$q^{-1}$}{$-1$}{$-q$}{$-q^{-1}$}
\quad \Rightarrow \,
\Drightofway{b}{$q$}{$q^{-1}$}{$-1$}{$q$}{$-1$}{$-1$}{$-q^{-1}$}{$-1$}
\quad \Rightarrow \,
\Dchainfour{}{$q$}{$q^{-1}$}{$-1$}{$-1$}{$-1$}{$-q$}{$-q^{-1}$}\quad ,
\end{align*}
any prolongation of the first diagram to the left is Weyl equivalent
to a diagram which was already considered in Step~2.2.

\textit{Step~3.8. Prolongations of the last diagram in row~17
of Table~\ref{a-r4}.}
By Lemma~\ref{l-5chainmid} for any prolongation to the right one has
to have the relation $\{-\zeta ,-\zeta ^{-1}\}=\{\zeta ,\zeta ^{-1}\}$
which is a contradiction to $\zeta \in R_3$. Moreover, since the last diagram
of row~17 is Weyl equivalent to the fifth one via
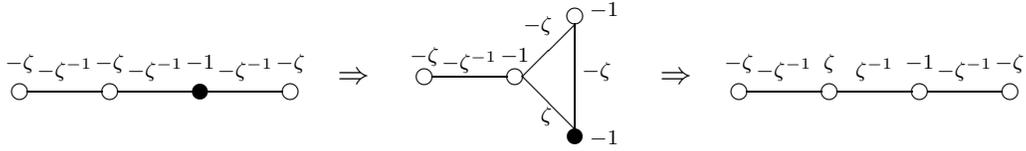
\begin{align*}
\Dchainfour{3}{$-\zeta $}{$-\zeta ^{-1}$}{$-\zeta $}{$-\zeta ^{-1}$}%
{$-1$}{$-\zeta ^{-1}$}{$-\zeta $}
\quad \Rightarrow \quad
\Drightofway{b}{$-\zeta $}{$-\zeta ^{-1}$}{$-1$}{$-\zeta $}{$\zeta $}%
{$-1$}{$-\zeta $}{$-1$}
\quad \Rightarrow \quad
\Dchainfour{}{$-\zeta $}{$-\zeta ^{-1}$}{$\zeta $}{$\zeta ^{-1}$}%
{$-1$}{$-\zeta ^{-1}$}{$-\zeta $}
\end{align*}
any prolongation of the last diagram of row~17 to the left is Weyl equivalent
to a diagram which was already considered in Step~2.3.

\textit{Step~3.9. Prolongations of the first two diagrams in row~18
of Table~\ref{a-r4}.}
By Weyl equivalence it is sufficient to consider one of these
diagrams. Prolongations to the right of the first diagram:
assume first that $\gDd _{\chi ,E}$
takes the form
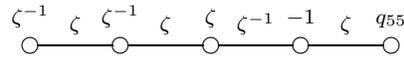
\begin{gather*}\label{eq-18-1-r}
%
%
\rule[-4\unitlength]{0pt}{5\unitlength}
\begin{picture}(50,4)(0,3)
\put(1,1){\circle{2}}
\put(2,1){\line(1,0){10}}
\put(13,1){\circle{2}}
\put(14,1){\line(1,0){10}}
\put(25,1){\circle{2}}
\put(26,1){\line(1,0){10}}
\put(37,1){\circle{2}}
\put(38,1){\line(1,0){10}}
\put(49,1){\circle{2}}
\put(1,4){\makebox[0pt]{\scriptsize $\zeta ^{-1}$}}
\put(7,3){\makebox[0pt]{\scriptsize $\zeta $}}
\put(13,4){\makebox[0pt]{\scriptsize $\zeta ^{-1}$}}
\put(19,3){\makebox[0pt]{\scriptsize $\zeta $}}
\put(25,4){\makebox[0pt]{\scriptsize $\zeta $}}
\put(31,3){\makebox[0pt]{\scriptsize $\zeta ^{-1}$}}
\put(37,4){\makebox[0pt]{\scriptsize $-1$}}
\put(43,3){\makebox[0pt]{\scriptsize $\zeta $}}
\put(49,4){\makebox[0pt]{\scriptsize $q_{55}$}}
\end{picture}
\end{gather*}
where $\zeta \in R_3$ and $(q_{55}+1)(q_{55}\zeta -1)=0$.
If $q_{55}=\zeta ^{-1}$ then this diagram was already considered in
Step~3.5. If $q_{55}=-1$ then up to Weyl equivalence this diagram and its
prolongations were already considered in Step~3.2, see (\ref{eq-4-1-r}) and
below.

Assume now that $\zeta \in R_3$ and $\gDd _{\chi ,E}$ is a prolongation of
length 1 of the second diagram in row~18 of Table~\ref{a-r4}.
Then Weyl equivalence gives
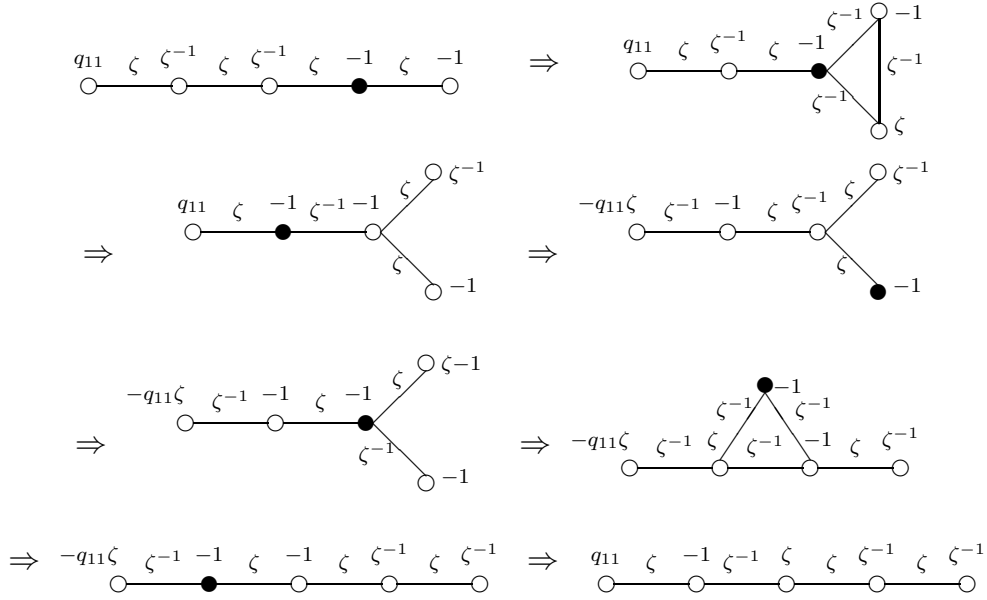
\begin{gather*}\label{eq-18-1-l}
%
%
\rule[-4\unitlength]{0pt}{5\unitlength}
\begin{picture}(50,4)(0,3)
\put(1,1){\circle{2}}
\put(2,1){\line(1,0){10}}
\put(13,1){\circle{2}}
\put(14,1){\line(1,0){10}}
\put(25,1){\circle{2}}
\put(26,1){\line(1,0){10}}
\put(37,1){\circle*{2}}
\put(38,1){\line(1,0){10}}
\put(49,1){\circle{2}}
\put(1,4){\makebox[0pt]{\scriptsize $q_{11}$}}
\put(7,3){\makebox[0pt]{\scriptsize $\zeta $}}
\put(13,4){\makebox[0pt]{\scriptsize $\zeta ^{-1}$}}
\put(19,3){\makebox[0pt]{\scriptsize $\zeta $}}
\put(25,4){\makebox[0pt]{\scriptsize $\zeta ^{-1}$}}
\put(31,3){\makebox[0pt]{\scriptsize $\zeta $}}
\put(37,4){\makebox[0pt]{\scriptsize $-1$}}
\put(43,3){\makebox[0pt]{\scriptsize $\zeta $}}
\put(49,4){\makebox[0pt]{\scriptsize $-1$}}
\end{picture}
\qquad \Rightarrow \qquad
\rule[-10\unitlength]{0pt}{12\unitlength}
\begin{picture}(36,10)(0,9)
\put(1,9){\circle{2}}
\put(2,9){\line(1,0){10}}
\put(13,9){\circle{2}}
\put(14,9){\line(1,0){10}}
\put(25,9){\circle*{2}}
\put(26,9){\line(1,1){7}}
\put(26,9){\line(1,-1){7}}
\put(33,17){\circle{2}}
\put(33,1){\circle{2}}
\put(33,2){\line(0,1){14}}
\put(1,12){\makebox[0pt]{\scriptsize $q_{11}$}}
\put(7,11){\makebox[0pt]{\scriptsize $\zeta $}}
\put(13,12){\makebox[0pt]{\scriptsize $\zeta ^{-1}$}}
\put(19,11){\makebox[0pt]{\scriptsize $\zeta $}}
\put(24,12){\makebox[0pt]{\scriptsize $-1$}}
\put(31,15){\makebox[0pt][r]{\scriptsize $\zeta ^{-1}$}}
\put(29,4){\makebox[0pt][r]{\scriptsize $\zeta ^{-1}$}}
\put(35,16){\makebox[0pt][l]{\scriptsize $-1$}}
\put(34,9){\makebox[0pt][l]{\scriptsize $\zeta ^{-1}$}}
\put(35,1){\makebox[0pt][l]{\scriptsize $\zeta $}}
\end{picture}
\\
\Rightarrow \qquad
\rule[-10\unitlength]{0pt}{12\unitlength}
\begin{picture}(36,14)(0,9)
\put(1,13){\circle{2}}
\put(2,13){\line(1,0){10}}
\put(13,13){\circle*{2}}
\put(14,13){\line(1,0){10}}
\put(25,13){\circle{2}}
\put(26,13){\line(1,1){7}}
\put(26,13){\line(1,-1){7}}
\put(33,21){\circle{2}}
\put(33,5){\circle{2}}
\put(1,16){\makebox[0pt]{\scriptsize $q_{11}$}}
\put(7,15){\makebox[0pt]{\scriptsize $\zeta $}}
\put(13,16){\makebox[0pt]{\scriptsize $-1$}}
\put(19,15){\makebox[0pt]{\scriptsize $\zeta ^{-1}$}}
\put(24,16){\makebox[0pt]{\scriptsize $-1$}}
\put(30,18){\makebox[0pt][r]{\scriptsize $\zeta $}}
\put(29,8){\makebox[0pt][r]{\scriptsize $\zeta $}}
\put(35,20){\makebox[0pt][l]{\scriptsize $\zeta ^{-1}$}}
\put(35,5){\makebox[0pt][l]{\scriptsize $-1$}}
\end{picture}
\qquad \Rightarrow \qquad
\rule[-10\unitlength]{0pt}{12\unitlength}
\begin{picture}(36,14)(0,9)
\put(1,13){\circle{2}}
\put(2,13){\line(1,0){10}}
\put(13,13){\circle{2}}
\put(14,13){\line(1,0){10}}
\put(25,13){\circle{2}}
\put(26,13){\line(1,1){7}}
\put(26,13){\line(1,-1){7}}
\put(33,21){\circle{2}}
\put(33,5){\circle*{2}}
\put(1,16){\makebox[0pt][r]{\scriptsize $-q_{11}\zeta $}}
\put(7,15){\makebox[0pt]{\scriptsize $\zeta ^{-1}$}}
\put(13,16){\makebox[0pt]{\scriptsize $-1$}}
\put(19,15){\makebox[0pt]{\scriptsize $\zeta $}}
\put(24,16){\makebox[0pt]{\scriptsize $\zeta ^{-1}$}}
\put(30,18){\makebox[0pt][r]{\scriptsize $\zeta $}}
\put(29,8){\makebox[0pt][r]{\scriptsize $\zeta $}}
\put(35,20){\makebox[0pt][l]{\scriptsize $\zeta ^{-1}$}}
\put(35,5){\makebox[0pt][l]{\scriptsize $-1$}}
\end{picture}
\\
\Rightarrow \qquad
\rule[-10\unitlength]{0pt}{12\unitlength}
\begin{picture}(36,14)(0,9)
\put(1,13){\circle{2}}
\put(2,13){\line(1,0){10}}
\put(13,13){\circle{2}}
\put(14,13){\line(1,0){10}}
\put(25,13){\circle*{2}}
\put(26,13){\line(1,1){7}}
\put(26,13){\line(1,-1){7}}
\put(33,21){\circle{2}}
\put(33,5){\circle{2}}
\put(1,16){\makebox[0pt][r]{\scriptsize $-q_{11}\zeta $}}
\put(7,15){\makebox[0pt]{\scriptsize $\zeta ^{-1}$}}
\put(13,16){\makebox[0pt]{\scriptsize $-1$}}
\put(19,15){\makebox[0pt]{\scriptsize $\zeta $}}
\put(24,16){\makebox[0pt]{\scriptsize $-1$}}
\put(30,18){\makebox[0pt][r]{\scriptsize $\zeta $}}
\put(29,8){\makebox[0pt][r]{\scriptsize $\zeta ^{-1}$}}
\put(35,20){\makebox[0pt][l]{\scriptsize $\zeta {-1}$}}
\put(35,5){\makebox[0pt][l]{\scriptsize $-1$}}
\end{picture}
\qquad \Rightarrow \qquad
\rule[-4\unitlength]{0pt}{5\unitlength}
\begin{picture}(38,11)(0,3)
\put(1,1){\circle{2}}
\put(2,1){\line(1,0){10}}
\put(13,1){\circle{2}}
\put(13,2){\line(2,3){6}}
\put(14,1){\line(1,0){10}}
\put(25,1){\circle{2}}
\put(25,2){\line(-2,3){6}}
\put(19,12){\circle*{2}}
\put(26,1){\line(1,0){10}}
\put(37,1){\circle{2}}
\put(1,4){\makebox[0pt][r]{\scriptsize $-q_{11}\zeta $}}
\put(7,3){\makebox[0pt]{\scriptsize $\zeta ^{-1}$}}
\put(12,4){\makebox[0pt]{\scriptsize $\zeta $}}
\put(19,3){\makebox[0pt]{\scriptsize $\zeta ^{-1}$}}
\put(26,4){\makebox[0pt]{\scriptsize $-1$}}
\put(31,3){\makebox[0pt]{\scriptsize $\zeta $}}
\put(37,4){\makebox[0pt]{\scriptsize $\zeta ^{-1}$}}
\put(15,8){\makebox[0pt]{\scriptsize $\zeta ^{-1}$}}
\put(23,8){\makebox[0pt][l]{\scriptsize $\zeta ^{-1}$}}
\put(22,11){\makebox[0pt]{\scriptsize $-1$}}
\end{picture}
\\
\Rightarrow \qquad
%
%
\rule[-4\unitlength]{0pt}{5\unitlength}
\begin{picture}(50,4)(0,3)
\put(1,1){\circle{2}}
\put(2,1){\line(1,0){10}}
\put(13,1){\circle*{2}}
\put(14,1){\line(1,0){10}}
\put(25,1){\circle{2}}
\put(26,1){\line(1,0){10}}
\put(37,1){\circle{2}}
\put(38,1){\line(1,0){10}}
\put(49,1){\circle{2}}
\put(1,4){\makebox[0pt][r]{\scriptsize $-q_{11}\zeta $}}
\put(7,3){\makebox[0pt]{\scriptsize $\zeta ^{-1}$}}
\put(13,4){\makebox[0pt]{\scriptsize $-1$}}
\put(19,3){\makebox[0pt]{\scriptsize $\zeta $}}
\put(25,4){\makebox[0pt]{\scriptsize $-1$}}
\put(31,3){\makebox[0pt]{\scriptsize $\zeta $}}
\put(37,4){\makebox[0pt]{\scriptsize $\zeta ^{-1}$}}
\put(43,3){\makebox[0pt]{\scriptsize $\zeta $}}
\put(49,4){\makebox[0pt]{\scriptsize $\zeta ^{-1}$}}
\end{picture}
\quad \Rightarrow \quad
%
%
\rule[-4\unitlength]{0pt}{5\unitlength}
\begin{picture}(50,4)(0,3)
\put(1,1){\circle{2}}
\put(2,1){\line(1,0){10}}
\put(13,1){\circle{2}}
\put(14,1){\line(1,0){10}}
\put(25,1){\circle{2}}
\put(26,1){\line(1,0){10}}
\put(37,1){\circle{2}}
\put(38,1){\line(1,0){10}}
\put(49,1){\circle{2}}
\put(1,4){\makebox[0pt]{\scriptsize $q_{11}$}}
\put(7,3){\makebox[0pt]{\scriptsize $\zeta $}}
\put(13,4){\makebox[0pt]{\scriptsize $-1$}}
\put(19,3){\makebox[0pt]{\scriptsize $\zeta ^{-1}$}}
\put(25,4){\makebox[0pt]{\scriptsize $\zeta $}}
\put(31,3){\makebox[0pt]{\scriptsize $\zeta $}}
\put(37,4){\makebox[0pt]{\scriptsize $\zeta ^{-1}$}}
\put(43,3){\makebox[0pt]{\scriptsize $\zeta $}}
\put(49,4){\makebox[0pt]{\scriptsize $\zeta ^{-1}$}}
\end{picture}
\end{gather*}
The latter diagram was already considered in the first part of this step.
In the same way one can see that any prolongation to the left of the
second diagram in row~18 of Table~\ref{a-r4} is Weyl equivalent to a
prolongation of the first diagram to the right and hence it was
already considered in the first part of this step.

\textit{Step~3.10. Prolongations of the first six diagrams in row~20
of Table~\ref{a-r4}.} One can see that such prolongations are Weyl equivalent
to prolongations of the first diagram.

Assume first that $\gDd _{\chi ,E}$ is of the form
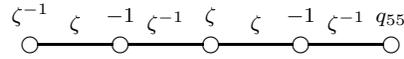
\begin{align*}
\Dchainfive{$\zeta ^{-1}$}{$\zeta $}{$-1$}{$\zeta ^{-1}$}{$\zeta $}%
{$\zeta $}{$-1$}{$\zeta ^{-1}$}{$q_{55}$}
\end{align*}
where $\zeta \in R_3$. Then one obtains a contradiction to
(\ref{eq-5chainmidvrel}).

Now suppose that $\gDd _{\chi ,E}$ is a prolongation of length one
to the left, that is
\begin{gather}\label{eq-20-1-l}
%
%
\rule[-4\unitlength]{0pt}{5\unitlength}
\begin{picture}(50,4)(0,3)
\put(1,1){\circle{2}}
\put(2,1){\line(1,0){10}}
\put(13,1){\circle{2}}
\put(14,1){\line(1,0){10}}
\put(25,1){\circle*{2}}
\put(26,1){\line(1,0){10}}
\put(37,1){\circle{2}}
\put(38,1){\line(1,0){10}}
\put(49,1){\circle{2}}
\put(1,4){\makebox[0pt]{\scriptsize $q_{11}$}}
\put(7,3){\makebox[0pt]{\scriptsize $\zeta $}}
\put(13,4){\makebox[0pt]{\scriptsize $\zeta ^{-1}$}}
\put(19,3){\makebox[0pt]{\scriptsize $\zeta $}}
\put(25,4){\makebox[0pt]{\scriptsize $-1$}}
\put(31,3){\makebox[0pt]{\scriptsize $\zeta ^{-1}$}}
\put(37,4){\makebox[0pt]{\scriptsize $\zeta $}}
\put(43,3){\makebox[0pt]{\scriptsize $\zeta $}}
\put(49,4){\makebox[0pt]{\scriptsize $-1$}}
\end{picture}
\qquad \Rightarrow \quad
%
%
\rule[-4\unitlength]{0pt}{5\unitlength}
\begin{picture}(50,4)(0,3)
\put(1,1){\circle{2}}
\put(2,1){\line(1,0){10}}
\put(13,1){\circle{2}}
\put(14,1){\line(1,0){10}}
\put(25,1){\circle{2}}
\put(26,1){\line(1,0){10}}
\put(37,1){\circle{2}}
\put(38,1){\line(1,0){10}}
\put(49,1){\circle{2}}
\put(1,4){\makebox[0pt]{\scriptsize $q_{11}$}}
\put(7,3){\makebox[0pt]{\scriptsize $\zeta $}}
\put(13,4){\makebox[0pt]{\scriptsize $-1$}}
\put(19,3){\makebox[0pt]{\scriptsize $\zeta ^{-1}$}}
\put(25,4){\makebox[0pt]{\scriptsize $-1$}}
\put(31,3){\makebox[0pt]{\scriptsize $\zeta $}}
\put(37,4){\makebox[0pt]{\scriptsize $-1$}}
\put(43,3){\makebox[0pt]{\scriptsize $\zeta $}}
\put(49,4){\makebox[0pt]{\scriptsize $-1$}}
\end{picture}\quad
\end{gather}
where $\zeta \in R_3$.
Then $\roots (\chi ;\Ndbasis _1,\Ndbasis _2+\Ndbasis _3,
\Ndbasis _3+\Ndbasis _4,\Ndbasis _4+\Ndbasis _5)$
with respect to the last
diagram has generalized Dynkin diagram
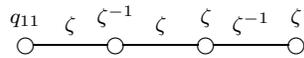
\begin{align*}
\Dchainfour{}{$q_{11}$}{$\zeta $}{$\zeta ^{-1}$}{$\zeta $}{$\zeta $}%
{$\zeta ^{-1}$}{$\zeta $}
\end{align*}
and hence Theorem~\ref{t-classrank4} implies that $q_{11}=\zeta ^{-1}$.
In this case the first diagram of (\ref{eq-20-1-l}) coincides with the
eleventh diagram in row~13 of Table~\ref{a-rb4}.

By the above argumentations a
prolongation to the left of length 2 of the first diagram in row~20
of Table~\ref{a-r4} has to take the form
\begin{gather*}\label{eq-20-1-ll}
%
%
\rule[-4\unitlength]{0pt}{5\unitlength}
\begin{picture}(62,4)(0,3)
\put(1,1){\circle{2}}
\put(2,1){\line(1,0){10}}
\put(13,1){\circle{2}}
\put(14,1){\line(1,0){10}}
\put(25,1){\circle{2}}
\put(26,1){\line(1,0){10}}
\put(37,1){\circle{2}}
\put(38,1){\line(1,0){10}}
\put(49,1){\circle{2}}
\put(50,1){\line(1,0){10}}
\put(61,1){\circle{2}}
\put(1,4){\makebox[0pt]{\scriptsize $q_{11}$}}
\put(7,3){\makebox[0pt]{\scriptsize $\zeta $}}
\put(13,4){\makebox[0pt]{\scriptsize $\zeta ^{-1}$}}
\put(19,3){\makebox[0pt]{\scriptsize $\zeta $}}
\put(25,4){\makebox[0pt]{\scriptsize $\zeta ^{-1}$}}
\put(31,3){\makebox[0pt]{\scriptsize $\zeta $}}
\put(37,4){\makebox[0pt]{\scriptsize $-1$}}
\put(43,3){\makebox[0pt]{\scriptsize $\zeta ^{-1}$}}
\put(49,4){\makebox[0pt]{\scriptsize $\zeta $}}
\put(55,3){\makebox[0pt]{\scriptsize $\zeta $}}
\put(61,4){\makebox[0pt]{\scriptsize $-1$}}
\end{picture}
\end{gather*}
where $\zeta \in R_3$ and $(q_{11}+1)(q_{11}\zeta -1)=0$.
However in this case the finiteness of
$\roots (\chi ;\Ndbasis _1+\Ndbasis _2,\Ndbasis _2+\Ndbasis _3,
\Ndbasis _3+\Ndbasis _4,\Ndbasis _4+\Ndbasis _5,\Ndbasis _5+\Ndbasis _6)$
\begin{align*}
\Dchainfive{$q_{11}$}{$\zeta $}{$-1$}{$\zeta ^{-1}$}{$-\zeta ^{-1}$}%
{$\zeta ^{-1}$}{$-1$}{$\zeta $}{$\zeta ^{-1}$}
\end{align*}
gives a contradiction to Theorem~\ref{t-classrank4}.

\textit{Step~3.11. Prolongations of the first six diagrams in row~21
of Table~\ref{a-r4}.} One can see that such prolongations are Weyl equivalent
to prolongations of the first diagram.

Assume first that $\gDd _{\chi ,E}$ is of the form
\begin{align*}
\Dchainfive{$-1$}{$\zeta ^{-1}$}{$\zeta $}{$\zeta ^{-1}$}{$\zeta $}%
{$\zeta $}{$-1$}{$\zeta ^{-1}$}{$q_{55}$}
\end{align*}
where $\zeta \in R_3$. Then one obtains a contradiction to
(\ref{eq-5chainmidvrel}).

Now suppose that $\gDd _{\chi ,E}$ is a prolongation of length one
to the left of the first diagram in row~21 of Table~\ref{a-r4}, that is
\begin{gather}\label{eq-21-1-l}
%
%
\rule[-4\unitlength]{0pt}{5\unitlength}
\begin{picture}(50,4)(0,3)
\put(1,1){\circle{2}}
\put(2,1){\line(1,0){10}}
\put(13,1){\circle{2}}
\put(14,1){\line(1,0){10}}
\put(25,1){\circle{2}}
\put(26,1){\line(1,0){10}}
\put(37,1){\circle{2}}
\put(38,1){\line(1,0){10}}
\put(49,1){\circle{2}}
\put(1,4){\makebox[0pt]{\scriptsize $q_{11}$}}
\put(7,3){\makebox[0pt]{\scriptsize $\zeta $}}
\put(13,4){\makebox[0pt]{\scriptsize $-1$}}
\put(19,3){\makebox[0pt]{\scriptsize $\zeta ^{-1}$}}
\put(25,4){\makebox[0pt]{\scriptsize $\zeta $}}
\put(31,3){\makebox[0pt]{\scriptsize $\zeta ^{-1}$}}
\put(37,4){\makebox[0pt]{\scriptsize $\zeta $}}
\put(43,3){\makebox[0pt]{\scriptsize $\zeta $}}
\put(49,4){\makebox[0pt]{\scriptsize $-1$}}
\end{picture}
\end{gather}
where $\zeta \in R_3$ and $(q_{11}+1)(q_{11}\zeta -1)=0$.
If $q_{11}=-1$ then the finiteness of
$\roots (\chi ;\Ndbasis _1+\Ndbasis _2,
\Ndbasis _3,\Ndbasis _4,\Ndbasis _5)$
gives a contradiction to Theorem~\ref{t-classrank4}.
On the other hand, if $q_{11}=\zeta ^{-1}$ then
the transformation $s_{\Ndbasis '_1,E'}s_{\Ndbasis _2,E}$,
where $E'=s_{\Ndbasis _2,E}(E)$ and
$\Ndbasis '_1=\Ndbasis _1+\Ndbasis _2
=s_{\Ndbasis _2,E}(\Ndbasis _1)$,
gives a prolongation to the left of the first diagram in row~20 of
Table~\ref{a-r4}.
Therefore (\ref{eq-21-1-l}) with $q_{11}=\zeta ^{-1}$
and all of its prolongations to the left
were already considered in Step~3.10.

\textit{Step~3.12. Prolongations of the first three diagrams in row~22
of Table~\ref{a-r4}.} Again all such prolongations are Weyl equivalent
to prolongations of the first diagram. Moreover, a prolongation of length one
to the right of the first diagram in row~22 of Table~\ref{a-r4}
would be of the form
\begin{align*}
\Dchainfive{$-\zeta $}{$\zeta $}{$-1$}{$-\zeta $}{$\zeta $}{$\zeta $}%
{$-\zeta $}{$\zeta $}{$q_{55}$}
\end{align*}
where $\zeta \in R_4$, which is a contradiction to (\ref{eq-5chainmiderel}).

Finally, a prolongation of length one
to the left of the first diagram in row~22 of Table~\ref{a-r4}
is via the transformations given below Weyl equivalent to a diagram
which was already considered in Step~2.4.
\begin{gather*}
%
%
\rule[-4\unitlength]{0pt}{5\unitlength}
\begin{picture}(50,4)(0,3)
\put(1,1){\circle{2}}
\put(2,1){\line(1,0){10}}
\put(13,1){\circle{2}}
\put(14,1){\line(1,0){10}}
\put(25,1){\circle*{2}}
\put(26,1){\line(1,0){10}}
\put(37,1){\circle{2}}
\put(38,1){\line(1,0){10}}
\put(49,1){\circle{2}}
\put(1,4){\makebox[0pt]{\scriptsize $q_{11}$}}
\put(7,3){\makebox[0pt]{\scriptsize $\zeta $}}
\put(13,4){\makebox[0pt]{\scriptsize $-\zeta $}}
\put(19,3){\makebox[0pt]{\scriptsize $\zeta $}}
\put(25,4){\makebox[0pt]{\scriptsize $-1$}}
\put(31,3){\makebox[0pt]{\scriptsize $-\zeta $}}
\put(37,4){\makebox[0pt]{\scriptsize $\zeta $}}
\put(43,3){\makebox[0pt]{\scriptsize $\zeta $}}
\put(49,4){\makebox[0pt]{\scriptsize $-\zeta $}}
\end{picture}
\qquad \Rightarrow \qquad
%
%
\rule[-4\unitlength]{0pt}{5\unitlength}
\begin{picture}(50,4)(0,3)
\put(1,1){\circle{2}}
\put(2,1){\line(1,0){10}}
\put(13,1){\circle{2}}
\put(14,1){\line(1,0){10}}
\put(25,1){\circle{2}}
\put(26,1){\line(1,0){10}}
\put(37,1){\circle*{2}}
\put(38,1){\line(1,0){10}}
\put(49,1){\circle{2}}
\put(1,4){\makebox[0pt]{\scriptsize $q_{11}$}}
\put(7,3){\makebox[0pt]{\scriptsize $\zeta $}}
\put(13,4){\makebox[0pt]{\scriptsize $-1$}}
\put(19,3){\makebox[0pt]{\scriptsize $-\zeta $}}
\put(25,4){\makebox[0pt]{\scriptsize $-1$}}
\put(31,3){\makebox[0pt]{\scriptsize $\zeta $}}
\put(37,4){\makebox[0pt]{\scriptsize $-1$}}
\put(43,3){\makebox[0pt]{\scriptsize $\zeta $}}
\put(49,4){\makebox[0pt]{\scriptsize $-\zeta $}}
\end{picture}\\
\Rightarrow \qquad
\rule[-10\unitlength]{0pt}{12\unitlength}
\begin{picture}(36,10)(0,9)
\put(1,9){\circle{2}}
\put(2,9){\line(1,0){10}}
\put(13,9){\circle{2}}
\put(14,9){\line(1,0){10}}
\put(25,9){\circle{2}}
\put(26,9){\line(1,1){7}}
\put(26,9){\line(1,-1){7}}
\put(33,17){\circle{2}}
\put(33,1){\circle*{2}}
\put(33,2){\line(0,1){14}}
\put(1,12){\makebox[0pt]{\scriptsize $q_{11}$}}
\put(7,11){\makebox[0pt]{\scriptsize $\zeta $}}
\put(13,12){\makebox[0pt]{\scriptsize $-1$}}
\put(19,11){\makebox[0pt]{\scriptsize $-\zeta $}}
\put(24,12){\makebox[0pt]{\scriptsize $\zeta $}}
\put(31,15){\makebox[0pt][r]{\scriptsize $-\zeta $}}
\put(29,4){\makebox[0pt][r]{\scriptsize $-1$}}
\put(35,16){\makebox[0pt][l]{\scriptsize $-1$}}
\put(34,9){\makebox[0pt][l]{\scriptsize $-\zeta $}}
\put(35,1){\makebox[0pt][l]{\scriptsize $-1$}}
\end{picture}
\qquad \Rightarrow \qquad
\Dchainfive{$q_{11}$}{$\zeta $}{$-1$}{$-\zeta $}{$\zeta $}{$-1$}{$-1$}%
{$\zeta $}{$-\zeta $}\quad .
\end{gather*}
Prolongations to the left of length bigger than one can be handled
with the analogous transformations.

\textit{Step~4.} Assume that equations
$(q_{11}+1)(q_{11}q_{12}q_{21}-1)=0$ and
$(q_{dd}+1)(q_{dd}q_{d,d-1}q_{d-1,d}-1)=0$ hold and one has
$q_{ii}^2q_{i,i-1}q_{i-1,i}q_{i,i+1}q_{i+1,i}=1$ for all $i$
with $1<i<d$. Then $\gDd _{\chi ,E}$ is a simple chain and it
appears in row~1 or row~2 of Table~\ref{a-rb4}.

With Steps 1--4 the analysis of generalized Dynkin diagrams which are
labeled path graphs is finished and the theorem is proven.
\end{bew}

\begin{appendix}

\section{On the finiteness of the Weyl groupoid}
\label{a-(T,E)}

According to the proofs of Theorem~\ref{t-classrank4} and
\ref{t-classrank>4} for generalized Dynkin diagrams $\gDd _{\chi ,E}$
of Tables \ref{a-r4} and \ref{a-rb4} elements $(T,E)\in \extWBG _{\chi ,E}$
are given such that $T(E)$ consists of $d-1$ elements of $E$ and one element
of $-\roots ^+_E$. It suffices to consider one single representant in
each Weyl equivalence class, and generalized Dynkin diagrams
of Cartan type can be ignored.

The starting point in each row is the first diagram which is a labeled
path graph. 
The numbering of the vertices of this diagram
is from left to right by 1 to $d$.
%

Abbreviate $m_1\Ndbasis _1+m_2\Ndbasis _2+m_3\Ndbasis _3+m_4\Ndbasis _4$ by
$1^{m_1}2^{m_2}3^{m_3}4^{m_4}$, where $i^{m_i}$ is omitted if $m_i=0$ and
$i^{m_i}$ is replaced by $i$ if $m_i=1$. For the map $s_{\aNdbasis _i,F}$,
where $\aNdbasis _i$ is the $i$th element of the basis $F$,
the symbol $s_i$ will be used. Lower indices to bases of $\ndZ ^d$
indicate the diagram corresponding to this basis.

\subsection{Computations for Table~\ref{a-r4}}

\noindent
Rows 6 and 10: $T:=s_4s_3s_2s_1$.
\begin{center}
$E=(1,2,3,4)$ $\mapsto $ $(-1,12,3,4)$ $\mapsto $ $(2,-12,123,4)$
$\mapsto $ $(2,3,-123,1234)$ $\mapsto $ $(2,3,4,-1234)$
\end{center}

\noindent
Rows 7,11,15,16 and 19: $T=s_1s_2s_3s_4s_3s_2s_1$.
\begin{center}
$(1,2,3,4)$ $\mapsto $ $(-1,12,3,4)$ $\mapsto $ $(2,-12,123,4)$
$\mapsto $ $(2,3,-123,1234)$ $\mapsto $ $(2,3,1234^2,-1234)$
$\mapsto $ $(2,123^24^2,-1234^2,4)$ $\mapsto $ $(12^23^24^2,-123^24^2,3,4)$
$\mapsto $ $(-12^23^24^2,2,3,4)$
\end{center}

\noindent
Row 8: $T=s_1 s_2 s_4 s_3 s_2 s_1$.
\begin{center}
$(1,2,3,4)_1$
$\mapsto $ $(-1,12,3,4)_2$
$\mapsto $ $(2,-12,123,4)_3$
$\mapsto $ $(2,3,-123,1234)_4$
$\mapsto $ $(2,123^24,4,-1234)_3$
$\mapsto $ $(12^23^24,-123^24,4,3)_2$
$\mapsto $ $(-12^23^24,2,4,3)_1$
\end{center}

\noindent
Row 9: $T=s_1s_4s_2s_3s_4s_1s_2s_3s_2s_1$.
\begin{center}
$(1,2,3,4)_1$
$\mapsto $ $(-1,12,3,4)_1$
$\mapsto $ $(2,-12,123,4)_1$
$\mapsto $ $(2,123^2,-123,1234)_1$
$\mapsto $ $(12^23^2,-123^2,3,1234)_1$
$\mapsto $ $(-12^23^2,2,3,1234)_1$
$\mapsto $ $(-12^23^2,2,123^24,-1234)_2$
$\mapsto $ $(-12^23^2,12^23^24,-123^24,3)_3$
$\mapsto $ $(4,-12^23^24,2,12^23^34)_4$
$\mapsto $ $(12^23^34^2,3,2,-12^23^34)_6$
$\mapsto $ $(-12^23^34^2,3,2,4)_6$
\end{center}

\noindent
Row 12: $T=s_1 s_2 s_4 s_3 s_2 s_1$.
\begin{center}
$(1,2,3,4)_1$
$\mapsto $ $(-1,12,3,4)_1$
$\mapsto $ $(2,-12,123,4)_2$
$\mapsto $ $(2,3,-123,1234)_4$
$\mapsto $ $(2,123^24,4,-1234)_2$
$\mapsto $ $(12^23^24,-123^24,4,3)_1$
$\mapsto $ $(-12^23^24,2,4,3)_1$
\end{center}

\noindent
Row 13: $T=s_1 s_2 s_4 s_3 s_2 s_1$.
\begin{center}
$(1,2,3,4)_1$
$\mapsto $ $(-1,12,3,4)_1$
$\mapsto $ $(2,-12,123,4)_1$
$\mapsto $ $(2,3,-123,1234)_2$
$\mapsto $ $(2,123^24,4,-1234)_1$
$\mapsto $ $(12^23^24,-123^24,4,3)_1$
$\mapsto $ $(-12^23^24,2,4,3)_1$
\end{center}

\noindent
Row 14: $T= s_3 s_4 s_1 s_2 s_4 s_3 s_2 s_1$.
\begin{center}
$(1,2,3,4)_1$
$\mapsto $ $(-1,12,3,4)_1$
$\mapsto $ $(2,-12,123,4)_1$
$\mapsto $ $(2,3,-123,1234)_2$
$\mapsto $ $(2,123^24,4,-1234)_4$
$\mapsto $ $(12^23^24,-123^24,4,3)_2$
$\mapsto $ $(-12^23^24,2,4,12^23^34)_5$
$\mapsto $ $(3,2,12^23^34^2,-12^23^34)_5$
$\mapsto $ $(3,2,-12^23^34^2,4)_5$
\end{center}

\noindent
Row 17: $T= s_1 s_2 s_3 s_4 s_3 s_2 s_1 s_4 s_1 s_2 s_3 s_4 s_3 s_2 s_1$.
\begin{center}
$(1,2,3,4)_1$
$\mapsto $ $(-1,12,3,4)_1$
$\mapsto $ $(2,-12,123,4)_1$
$\mapsto $ $(2,3,-123,1234)_1$
$\mapsto $ $(2,3,1234^2,-1234)_2$
$\mapsto $ $(2,123^24^2,-1234^2,4)_3$
$\mapsto $ $(12^23^24^2,-123^24^2,3,123^24^3)_4$
$\mapsto $ $(-12^23^24^2,2,3,1^22^33^44^5)_6$
$\mapsto $ $(123^24^3,2,3,-1^22^33^44^5)_6$
$\mapsto $ $(-123^24^3,12^23^24^3,3,-12^23^24^2)_4$
$\mapsto $ $(2,-12^23^24^3,12^23^34^3,4)_3$
$\mapsto $ $(2,3,-12^23^34^3,12^23^34^4)_2$
$\mapsto $ $(2,3,12^23^34^5,-12^23^34^4)_1$
$\mapsto $ $(2,12^23^44^5,-12^23^34^5,4)_1$
$\mapsto $ $(12^33^44^5,-12^23^44^5,3,4)_1$
$\mapsto $ $(-12^33^44^5,2,3,4)_1$
\end{center}

\noindent
Row 18: $T= s_4 s_2 s_3 s_4 s_1 s_2 s_3 s_2 s_1$.
\begin{center}
$(1,2,3,4)_1$
$\mapsto $ $(-1,12,3,4)_1$
$\mapsto $ $(2,-12,123,4)_1$
$\mapsto $ $(2,123^2,-123,1234)_1$
$\mapsto $ $(12^23^2,-123^2,3,1234)_1$
$\mapsto $ $(-12^23^2,2,3,1234)_1$
$\mapsto $ $(-12^23^2,2,123^24,-1234)_2$
$\mapsto $ $(-12^23^2,12^23^24,-123^24,3)_3$
$\mapsto $ $(4,-12^23^24,2,12^23^34)_4$
$\mapsto $ $(4,3,2,-12^23^34)_5$
\end{center}

\noindent
Row 20: $T= s_1 s_2 s_3 s_4 s_2 s_3 s_1 s_2 s_4 s_3 s_2 s_1$.
\begin{center}
$(1,2,3,4)_1$
$\mapsto $ $(-1,12,3,4)_1$
$\mapsto $ $(2,-12,123,4)_3$
$\mapsto $ $(2,3,-123,1234)_7$
$\mapsto $ $(2,123^24,4,-1234)_7$
$\mapsto $ $(12^23^24,-123^24,123^24^2,3)_7$
$\mapsto $ $(-12^23^24,2,123^24^2,3)_8$
$\mapsto $ $(-12^23^24,12^23^24^2,-123^24^2,123^34^2)_4$
$\mapsto $ $(4,-12^23^24^2,2,123^34^2)_4$
$\mapsto $ $(4,-12^23^24^2,12^23^34^2,-123^34^2)_6$
$\mapsto $ $(4,12^23^44^2,-12^23^34^2,2)_6$
$\mapsto $ $(12^23^44^3,-12^23^44^2,3,2)_6$
$\mapsto $ $(-12^23^44^3,4,3,2)_5$
\end{center}

\noindent
Row 21: $T= s_1 s_2 s_3 s_4 s_2 s_3 s_4 s_1 s_2 s_3 s_4 s_2 s_3 s_1
 s_2 s_4 s_3 s_2 s_1$.
\begin{center}
$(1,2,3,4)_1$
$\mapsto $ $(-1,12,3,4)_3$
$\mapsto $ $(2,-12,123,4)_6$
$\mapsto $ $(2,3,-123,1234)_7$
$\mapsto $ $(2,123^24,4,-1234)_7$
$\mapsto $ $(12^23^24,-123^24,123^24^2,3)_7$
$\mapsto $ $(-12^23^24,2,123^24^2,3)_7$
$\mapsto $ $(-12^23^24,12^23^24^2,-123^24^2,123^34^2)_6$
$\mapsto $ $(4,-12^23^24^2,2,123^34^2)_3$
$\mapsto $ $(4,-12^23^24^2,12^23^34^2,-123^34^2)_4$
$\mapsto $ $(4,12^23^44^2,-12^23^34^2,12^33^34^2)_4$
$\mapsto $ $(12^23^44^3,-12^23^44^2,3,12^33^34^2)_5$
$\mapsto $ $(-12^23^44^3,4,3,12^33^34^2)_5$
$\mapsto $ $(-12^23^44^3,4,12^33^44^2,-12^33^34^2)_6$
$\mapsto $ $(-12^23^44^3,12^33^44^3,-12^33^44^2,3)_7$
$\mapsto $ $(-12^23^44^3,12^33^44^3,-12^33^44^2,3)_7$
$\mapsto $ $(2,-12^33^44^3,4,12^33^54^3)_7$
$\mapsto $ $(2,3,12^33^54^4,-12^33^54^3)_7$
$\mapsto $ $(2,12^33^64^4,-12^33^54^4,4)_6$
$\mapsto $ $(12^43^64^4,-12^33^64^4,3,4)_3$
$\mapsto $ $(-12^43^64^4,2,3,4)_1$
\end{center}

\noindent
Row 22: $T= s_1 s_2 s_3 s_4 s_2 s_4 s_1 s_2 s_3 s_2 s_4 s_1
 s_2 s_4 s_3 s_2 s_1$.
\begin{center}
$(1,2,3,4)_1$
$\mapsto $ $(-1,12,3,4)_1$
$\mapsto $ $(2,-12,123,4)_2$
$\mapsto $ $(2,3,-123,1234)_4$
$\mapsto $ $(2,123^24,4,-1234)_6$
$\mapsto $ $(12^23^24,-123^24,4,123^34)_6$
$\mapsto $ $(-12^23^24,2,4,123^34)_7$
$\mapsto $ $(-12^23^24,12^23^34,123^34^2,-123^34)_5$
$\mapsto $ $(3,-12^23^34,1^22^33^64^3,2)_8$
$\mapsto $ $(3,123^34^2,-1^22^33^64^3,2)_8$
$\mapsto $ $(123^44^2,-123^34^2,-12^23^34,12^23^34^2)_5$
$\mapsto $ $(-123^44^2,3,-12^23^34,12^23^34^2)_4$
$\mapsto $ $(-123^44^2,12^23^44^2,4,-12^23^34^2)_6$
$\mapsto $ $(2,-12^23^44^2,4,12^23^54^2)_6$
$\mapsto $ $(2,3,12^23^54^3,-12^23^54^2)_4$
$\mapsto $ $(2,12^23^64^3,-12^23^54^3,4)_2$
$\mapsto $ $(12^33^64^3,-12^23^64^3,3,4)_1$
$\mapsto $ $(-12^33^64^3,2,3,4)_1$
\end{center}

\subsection{Computations for Table~\ref{a-rb4}}

\noindent
Row 2: $T:=s_d\cdots s_2s_1$.

\noindent
Rows 4,5 and 6: $T:=s_1s_2\cdots \cdots s_{d-1}s_ds_{d-1}\cdots s_2s_1$.

\noindent
Row 8, $i_1=1$: $T:=s_1s_2\cdots s_{d-3}s_{d-2}s_ds_{d-1}s_{d-2}\cdots s_2s_1$.

\noindent
Row 9, $(i_1,\ldots ,i_j)=(1,\ldots ,j)$:
$T:=s_1s_2\cdots s_{d-3}s_{d-2}s_ds_{d-1}s_{d-2}\cdots s_2s_1$.

\noindent
Row 10:
$T:= s_5 s_4 s_3 s_2 s_4 s_5 s_1 s_2 s_4 s_3 s_2 s_1$.

\noindent
Row 11:
$T:= s_1 s_2 s_4 s_3 s_5 s_2 s_3 s_4 s_5 s_1 s_2 s_3 s_4 s_3 s_2 s_1$.

\noindent
Row 12:
$T:= s_1 s_2 s_4 s_3 s_5 s_4 s_3 s_2 s_4 s_3 s_5 s_2 s_4 s_3 s_1 s_2 s_4
 s_2 s_3 s_5 s_1 s_2 s_3 s_4 s_5 s_1 s_2 s_3 s_4 \times $
$s_3 s_2 s_1$.

\noindent
Row 13:
$T:= s_1 s_2 s_3 s_4 s_5 s_2 s_4 s_3 s_5 s_1 s_2 s_4 s_3 s_2 s_1$.

\noindent
Row 14:
$T:= s_1 s_2 s_3 s_4 s_5 s_2 s_4 s_3 s_5 s_1 s_2 s_4 s_5 s_2 s_4 s_3 s_5
 s_1 s_2 s_4 s_3 s_2 s_1$.

\noindent
Row 16:
$T:= s_1 s_2 s_3 s_4 s_5 s_6 s_3 s_5 s_4 s_2 s_3 s_5 s_6 s_1 s_2 s_3 s_5
 s_4 s_3 s_2 s_1$.

\noindent
Row 17:
$T:= s_1 s_2 s_3 s_4 s_5 s_6 s_4 s_5 s_3 s_4 s_6 s_2 s_3 s_4 s_1 s_2 s_3
 s_5 s_3 s_4 s_6 s_2 s_3 s_4 s_5 s_6 s_1 s_2 s_3 \times $
$s_4 s_5 s_4 s_3 s_2 s_1$.

\noindent
Row 18:
$T:= s_1 s_2 s_3 s_4 s_5 s_6 s_3 s_5 s_4 s_6 s_2 s_3 s_5 s_1 s_2 s_3 s_4
 s_3 s_5 s_6 s_2 s_3 s_5 s_4 s_6 s_1 s_2 s_3 s_5 \times $
$s_4 s_3 s_2 s_1$.

\noindent
Row 20:
$T:= s_1 s_2 s_3 s_4 s_5 s_6 s_7 s_4 s_6 s_5 s_3 s_4 s_6 s_7 s_2 s_3 s_4
 s_6 s_1 s_2 s_3 s_4 s_5 s_4 s_6 s_7 s_3 s_4 s_6 \times $
$s_5 s_2 s_3 s_4 s_6 s_7 s_1 s_2 s_3 s_4 s_6 s_5 s_4 s_3 s_2 s_1$.

\section{Connected arithmetic root systems of rank four}
\label{a-r4}


\setlength{\unitlength}{1mm}
\begin{tabular}{r|l|l}
row & gener. Dynkin diagrams & fixed parameters \\
\hline \hline
 1 & \Dchainfour{}{$q$}{$q^{-1}$}{$q$}{$q^{-1}$}{$q$}{$q^{-1}$}{$q$}
& $q\in k^\ast \setminus \{1\}$ \\
\hline
 2 & \Dchainfour{}{$q^2$}{$q^{-2}$}{$q^2$}{$q^{-2}$}{$q^2$}{$q^{-2}$}{$q$}
& $q\in k^\ast \setminus \{-1,1\}$ \\
\hline
 3 & \Dchainfour{}{$q$}{$q^{-1}$}{$q$}{$q^{-1}$}{$q$}{$q^{-2}$}{$q^2$}
& $q\in k^\ast \setminus \{-1,1\}$ \\
\hline
 4 & \Dchainfour{}{$q^2$}{$q^{-2}$}{$q^2$}{$q^{-2}$}{$q$}{$q^{-1}$}{$q$}
& $q\in k^\ast \setminus \{-1,1\}$ \\
\hline
 5 & \Dthreefork{}{$q$}{$q^{-1}$}{$q$}{$q^{-1}$}{$q^{-1}$}{$q$}{$q$}
& $q\in k^\ast \setminus \{1\}$ \\
\hline
6 & \Dchainfour{}{$-1$}{$q^{-1}$}{$q$}{$q^{-1}$}{$q$}{$q^{-1}$}{$q$} &
$q\in k^\ast \setminus \{-1,1\}$ \\
  & \Dchainfour{}{$-1$}{$q$}{$-1$}{$q^{-1}$}{$q$}{$q^{-1}$}{$q$} & \\
  & \Dchainfour{}{$q$}{$q^{-1}$}{$-1$}{$q$}{$-1$}{$q^{-1}$}{$q$} & \\
\hline
7 & \Dchainfour{}{$-1$}{$q^{-2}$}{$q^2$}{$q^{-2}$}{$q^2$}{$q^{-2}$}{$q$} &
$q\in k^\ast$, $q^4\not=1$ \\
  & \Dchainfour{}{$-1$}{$q^2$}{$-1$}{$q^{-2}$}{$q^2$}{$q^{-2}$}{$q$} &\\
  & \Dchainfour{}{$q^2$}{$q^{-2}$}{$-1$}{$q^2$}{$-1$}{$q^{-2}$}{$q$} &\\
  & \Dchainfour{}{$q^2$}{$q^{-2}$}{$q^2$}{$q^{-2}$}{$-1$}{$q^2$}{$-q^{-1}$}
& \\
\hline
8 & \Dchainfour{}{$-1$}{$q^{-1}$}{$q$}{$q^{-1}$}{$q$}{$q^{-2}$}{$q^2$} &
$q\in k^\ast \setminus \{-1,1\}$ \\
  & \Dchainfour{}{$-1$}{$q$}{$-1$}{$q^{-1}$}{$q$}{$q^{-2}$}{$q^2$} & \\
  & \Dchainfour{}{$q$}{$q^{-1}$}{$-1$}{$q$}{$-1$}{$q^{-2}$}{$q^2$} & \\
  & \Drightofway{}{$q$}{$q^{-1}$}{$q$}{$q^{-1}$}{$q^{-1}$}{$-1$}{$q^2$}{$-1$}
& \\
\hline
\end{tabular}


\begin{tabular}{r|l|l}
row & gener. Dynkin diagrams & fixed param. \\
\hline \hline
9 & \Dchainfour{}{$q^2$}{$q^{-2}$}{$q^2$}{$q^{-2}$}{$q$}{$q^{-1}$}{$-1$} &
$q\in k^\ast $, $q^2,q^3\not=1$ \\
 & \Dchainfour{}{$q^2$}{$q^{-2}$}{$q^2$}{$q^{-2}$}{$-1$}{$q$}{$-1$} & \\
& \Drightofway{}{$q^2$}{$q^{-2}$}{$-1$}{$q^2$}{$q^{-1}$}{$-1$}{$q^{-1}$}{$q$}
\Drightofway{}{$q^2$}{$q^{-2}$}{$-1$}{$q^2$}{$q$}{$-1$}{$q^{-3}$}{$-1$}
& \\
& \Dchainfour{}{$q^2$}{$q^{-2}$}{$q^2$}{$q^{-2}$}{$-1$}{$q^3$}{$q^{-3}$} & \\
& \Dchainfour{}{$q^2$}{$q^{-2}$}{$q$}{$q^{-1}$}{$-1$}{$q^3$}{$q^{-3}$} & \\
\hline
10 & \Dchainfour{}{$q^{-1}$}{$q$}{$-1$}{$q^{-1}$}{$q$}{$q^{-1}$}{$q$} &
$q\in k^\ast \setminus \{-1,1\}$ \\
 & \Dchainfour{}{$-1$}{$q^{-1}$}{$-1$}{$q$}{$-1$}{$q^{-1}$}{$q$} & \\
 & \Dchainfour{}{$-1$}{$q$}{$q^{-1}$}{$q$}{$-1$}{$q^{-1}$}{$q$} & \\
 & \Dchainfour{}{$-1$}{$q^{-1}$}{$q$}{$q^{-1}$}{$-1$}{$q$}{$-1$} &\\
 & \Dchainfour{}{$-1$}{$q^{-1}$}{$q$}{$q^{-1}$}{$q$}{$q^{-1}$}{$-1$} & \\
 & \Dchainfour{}{$-1$}{$q$}{$-1$}{$q^{-1}$}{$-1$}{$q$}{$-1$} & \\
\hline
11 & \Dchainfour{}{$q^{-2}$}{$q^2$}{$-1$}{$q^{-2}$}{$q^2$}{$q^{-2}$}{$q$} &
$q\in k^\ast $, $q^4\not=1$ \\
 & \Dchainfour{}{$-1$}{$q^{-2}$}{$-1$}{$q^2$}{$-1$}{$q^{-2}$}{$q$} & \\
 & \Dchainfour{}{$-1$}{$q^2$}{$q^{-2}$}{$q^2$}{$-1$}{$q^{-2}$}{$q$} & \\
 & \Dchainfour{}{$-1$}{$q^{-2}$}{$q^2$}{$q^{-2}$}{$-1$}{$q^2$}{$-q^{-1}$} & \\
 & \Dchainfour{}{$-1$}{$q^2$}{$-1$}{$q^{-2}$}{$-1$}{$q^2$}{$-q^{-1}$} & \\
 & \Dchainfour{}{$q^2$}{$q^{-2}$}{$-1$}{$q^2$}{$q^{-2}$}{$q^2$}{$-q^{-1}$} & \\
\hline
\end{tabular}


\begin{tabular}{r|l|l}
row & gener. Dynkin diagrams & fixed param. \\
\hline \hline
12 & \Dchainfour{}{$q^{-1}$}{$q$}{$-1$}{$q^{-1}$}{$q$}{$q^{-2}$}{$q^2$} &
$q\in k^\ast \setminus \{-1,1\}$ \\
 & \Dchainfour{}{$-1$}{$q^{-1}$}{$-1$}{$q$}{$-1$}{$q^{-2}$}{$q^2$} & \\
 & \Dchainfour{}{$-1$}{$q$}{$q^{-1}$}{$q$}{$-1$}{$q^{-2}$}{$q^2$} & \\
 & \Drightofway{}{$-1$}{$q^{-1}$}{$q$}{$q^{-1}$}{$q^{-1}$}{$-1$}{$q^2$}{$-1$}\
\Drightofway{}{$-1$}{$q$}{$-1$}{$q^{-1}$}{$q^{-1}$}{$-1$}{$q^2$}{$-1$} & \\
 & \Dthreefork{}{$q$}{$q^{-1}$}{$-1$}{$q$}{$q$}{$q^{-1}$}{$q^{-1}$} & \\
\hline
13 & \Dchainfour{}{$q$}{$q^{-1}$}{$q$}{$q^{-1}$}{$-1$}{$q^2$}{$q^{-2}$} &
$q\in k^\ast \setminus \{-1,1\}$ \\
 & \Drightofway{}{$q$}{$q^{-1}$}{$-1$}{$q$}{$q$}{$-1$}{$q^{-2}$}{$-1$}\
\Dthreefork{}{$-1$}{$q$}{$-1$}{$q^{-1}$}{$q^{-1}$}{$q$}{$q$} & \\
 & \Dthreefork{}{$-1$}{$q^{-1}$}{$q$}{$q^{-1}$}{$q^{-1}$}{$q$}{$q$} & \\
\hline
14 & \Dchainfour{}{$q$}{$q^{-1}$}{$q$}{$q^{-1}$}{$-1$}{$-q$}{$-q^{-1}$} &
$q\in k^\ast \setminus \{-1,1\}$ \\
 & \Drightofway{}{$q$}{$q^{-1}$}{$-1$}{$q$}{$-1$}{$-1$}{$-q^{-1}$}{$-1$}\
\Drightofway{}{$-q^{-1}$}{$-q$}{$-1$}{$-q^{-1}$}{$-1$}{$-1$}{$q$}{$-1$} & \\
 & \Dchainfour{}{$q$}{$q^{-1}$}{$-1$}{$-1$}{$-1$}{$-q$}{$-q^{-1}$} & \\
 & \Dchainfour{}{$-q^{-1}$}{$-q$}{$-q^{-1}$}{$-q$}{$-1$}{$q^{-1}$}{$q$} & \\
\hline
\end{tabular}


\begin{tabular}{r|l|l}
row & gener. Dynkin diagrams & fixed param. \\
\hline \hline
15 & \Dchainfour{}{$-\zeta ^{-1}$}{$-\zeta $}{$-\zeta ^{-1}$}%
{$-\zeta $}{$-\zeta ^{-1}$}{$-\zeta $}{$\zeta $}
& $\zeta \in R_3$ \\
\hline
16 & \Dchainfour{}{$-1$}{$-\zeta $}{$-\zeta ^{-1}$}{$-\zeta $}%
{$-\zeta ^{-1}$}{$-\zeta $}{$\zeta $} & $\zeta \in R_3$ \\
  & \Dchainfour{}{$-1$}{$-\zeta ^{-1}$}{$-1$}{$-\zeta $}%
{$-\zeta ^{-1}$}{$-\zeta $}{$\zeta $} & \\
  & \Dchainfour{}{$-\zeta ^{-1}$}{$-\zeta $}{$-1$}{$-\zeta ^{-1}$}%
{$-1$}{$-\zeta $}{$\zeta $} & \\
  & \Dchainfour{}{$-\zeta ^{-1}$}{$-\zeta $}{$-\zeta ^{-1}$}{$-\zeta $}%
{$-1$}{$-\zeta ^{-1}$}{$\zeta ^{-1}$} & \\
\hline
17 & \Dchainfour{}{$-\zeta $}{$-\zeta ^{-1}$}{$-\zeta $}{$-\zeta ^{-1}$}%
{$-\zeta $}{$-\zeta ^{-1}$}{$\zeta $} & $\zeta \in R_3$ \\
  & \Dchainfour{}{$-\zeta $}{$-\zeta ^{-1}$}{$-\zeta $}{$-\zeta ^{-1}$}%
{$-1$}{$-1$}{$\zeta $} & \\
  & \Drightofway{}{$-\zeta $}{$-\zeta ^{-1}$}{$-1$}{$-\zeta $}%
{$\zeta ^{-1}$}{$-1$}{$-1$}{$\zeta $}
 \Drightofway{}{$-\zeta $}{$-\zeta ^{-1}$}{$-1$}{$-\zeta $}%
{$\zeta $}{$-1$}{$-\zeta $}{$-1$} & \\
  & \Dchainfour{}{$-\zeta $}{$-\zeta ^{-1}$}{$\zeta $}{$\zeta ^{-1}$}%
{$-1$}{$-\zeta ^{-1}$}{$-\zeta $} & \\
  & \Dchainfour{}{$-\zeta $}{$-\zeta ^{-1}$}{$-\zeta $}{$-\zeta ^{-1}$}%
{$-1$}{$-\zeta ^{-1}$}{$-\zeta $} & \\
\hline
18 & \Dchainfour{}{$\zeta ^{-1}$}{$\zeta $}{$\zeta ^{-1}$}{$\zeta $}%
{$\zeta $}{$\zeta ^{-1}$}{$-1$} & $\zeta \in R_3$ \\
  & \Dchainfour{}{$\zeta ^{-1}$}{$\zeta $}{$\zeta ^{-1}$}{$\zeta $}%
{$-1$}{$\zeta $}{$-1$} & \\
  & \Drightofway{}{$\zeta ^{-1}$}{$\zeta $}{$-1$}{$\zeta ^{-1}$}%
{$\zeta ^{-1}$}{$-1$}{$\zeta ^{-1}$}{$\zeta $}
  \Dthreefork{}{$\zeta ^{-1}$}{$\zeta $}{$-1$}{$\zeta ^{-1}$}%
{$\zeta $}{$-1$}{$-1$} & \\
  & \Dthreefork{}{$\zeta ^{-1}$}{$\zeta $}{$\zeta $}{$\zeta ^{-1}$}%
{$\zeta ^{-1}$}{$-1$}{$-1$}
  \Dthreefork{}{$\zeta ^{-1}$}{$\zeta $}{$\zeta ^{-1}$}{$\zeta $}%
{$\zeta $}{$-1$}{$-1$} & \\
\hline
\end{tabular}


\begin{tabular}{r|l|l}
row & gener. Dynkin diagrams & fixed param. \\
\hline \hline
19 & \Dchainfour{}{$-\zeta $}{$-\zeta ^{-1}$}{$-1$}{$-\zeta $}{$-\zeta ^{-1}$}%
{$-\zeta $}{$\zeta $} & $\zeta \in R_3$ \\
  & \Dchainfour{}{$-1$}{$-\zeta $}{$-1$}{$-\zeta ^{-1}$}{$-1$}%
{$-\zeta $}{$\zeta $} & \\
  & \Dchainfour{}{$-1$}{$-\zeta ^{-1}$}{$-\zeta $}{$-\zeta ^{-1}$}{$-1$}%
{$-\zeta $}{$\zeta $} & \\
  & \Dchainfour{}{$-1$}{$-\zeta $}{$-\zeta ^{-1}$}{$-\zeta $}{$-1$}%
{$-\zeta ^{-1}$}{$\zeta ^{-1}$} & \\
  & \Dchainfour{}{$-1$}{$-\zeta ^{-1}$}{$-1$}{$-\zeta $}{$-1$}%
{$-\zeta ^{-1}$}{$\zeta ^{-1}$} & \\
  & \Dchainfour{}{$-\zeta ^{-1}$}{$-\zeta $}{$-1$}{$-\zeta ^{-1}$}{$-\zeta $}%
{$-\zeta ^{-1}$}{$\zeta ^{-1}$} & \\
\hline
20 & \Dchainfour{}{$\zeta ^{-1}$}{$\zeta $}{$-1$}{$\zeta ^{-1}$}{$\zeta $}%
{$\zeta $}{$-1$} & $\zeta \in R_3$ \\
 & \Dchainfour{}{$\zeta ^{-1}$}{$\zeta $}{$-1$}{$\zeta ^{-1}$}{$-\zeta ^{-1}$}%
{$\zeta ^{-1}$}{$-1$} & \\
 & \Dchainfour{}{$-1$}{$\zeta ^{-1}$}{$-1$}{$\zeta $}{$-1$}%
{$\zeta $}{$-1$} & \\
 & \Dchainfour{}{$-1$}{$\zeta $}{$\zeta ^{-1}$}{$\zeta $}{$-1$}%
{$\zeta $}{$-1$} & \\
 & \Dchainfour{}{$-1$}{$\zeta ^{-1}$}{$-1$}{$\zeta $}{$\zeta $}%
{$\zeta ^{-1}$}{$-1$} & \\
 & \Dchainfour{}{$-1$}{$\zeta $}{$\zeta ^{-1}$}{$\zeta $}{$\zeta $}%
{$\zeta ^{-1}$}{$-1$} & \\
 & \Drightofway{}{$-1$}{$\zeta ^{-1}$}{$\zeta $}{$\zeta ^{-1}$}%
{$\zeta ^{-1}$}{$-1$}{$\zeta ^{-1}$}{$\zeta $}\
 \Drightofway{}{$-1$}{$\zeta $}{$-1$}{$\zeta ^{-1}$}%
{$\zeta ^{-1}$}{$-1$}{$\zeta ^{-1}$}{$\zeta $} & \\
 & \Dthreefork{}{$\zeta $}{$\zeta ^{-1}$}{$-1$}{$\zeta $}%
{$\zeta $}{$\zeta ^{-1}$}{$-1$}\
 \Dthreefork{}{$\zeta $}{$\zeta ^{-1}$}{$\zeta $}{$\zeta $}%
{$\zeta ^{-1}$}{$\zeta ^{-1}$}{$-1$} & \\
\hline
\end{tabular}


\begin{tabular}{r|l|l}
row & gener. Dynkin diagrams & fixed param. \\
\hline \hline
21  & \Dchainfour{}{$-1$}{$\zeta ^{-1}$}{$\zeta $}{$\zeta ^{-1}$}%
{$\zeta $}{$\zeta $}{$-1$} & $\zeta \in R_3$ \\
  & \Dchainfour{}{$-1$}{$\zeta ^{-1}$}{$\zeta $}{$\zeta ^{-1}$}%
{$-\zeta ^{-1}$}{$\zeta ^{-1}$}{$-1$} & \\
  & \Dchainfour{}{$-1$}{$\zeta $}{$-1$}{$\zeta ^{-1}$}{$\zeta $}%
{$\zeta $}{$-1$} & \\
  & \Dchainfour{}{$-1$}{$\zeta $}{$-1$}{$\zeta ^{-1}$}{$-\zeta ^{-1}$}%
{$\zeta ^{-1}$}{$-1$} & \\
  & \Dchainfour{}{$\zeta $}{$\zeta ^{-1}$}{$-1$}{$\zeta $}{$\zeta $}%
{$\zeta ^{-1}$}{$-1$} & \\
  & \Dchainfour{}{$\zeta $}{$\zeta ^{-1}$}{$-1$}{$\zeta $}{$-1$}%
{$\zeta $}{$-1$} & \\
  & \Drightofway{}{$\zeta $}{$\zeta ^{-1}$}{$\zeta $}{$\zeta ^{-1}$}%
{$\zeta ^{-1}$}{$-1$}{$\zeta ^{-1}$}{$\zeta $} & \\
\hline
22 & \Dchainfour{}{$-\zeta $}{$\zeta $}{$-1$}{$-\zeta $}{$\zeta $}%
{$\zeta $}{$-\zeta $} & $\zeta \in R_4$ \\
 & \Dchainfour{}{$-1$}{$-\zeta $}{$-1$}{$\zeta $}{$-1$}{$\zeta $}{$-\zeta $}
& \\
 & \Dchainfour{}{$-1$}{$\zeta $}{$-\zeta $}{$\zeta $}{$-1$}%
{$\zeta $}{$-\zeta $} & \\
 & \Drightofway{}{$-1$}{$-\zeta $}{$\zeta $}{$-\zeta $}%
{$-1$}{$-1$}{$-\zeta $}{$-1$}\
 \Drightofway{}{$-1$}{$\zeta $}{$-1$}{$-\zeta $}%
{$-1$}{$-1$}{$-\zeta $}{$-1$} & \\
& \Dchainfour{}{$-1$}{$-\zeta $}{$\zeta $}{$-1$}{$-1$}%
{$\zeta $}{$-\zeta $} & \\
& \Dchainfour{}{$-1$}{$\zeta $}{$-1$}{$-1$}{$-1$}{$\zeta $}{$-\zeta $} & \\
& \Drightofway{}{$-\zeta $}{$\zeta $}{$-1$}{$-\zeta $}%
{$-1$}{$\zeta $}{$-\zeta $}{$-1$} & \\
\hline
\end{tabular}

\newpage

\section{Connected arithmetic root systems of rank $d\ge 5$}
\label{a-rb4}

Recall the notation for simple chains in Section \ref{sec-simplechains}.

\vspace{\baselineskip}

\begin{tabular}{r|l|l}
row & gener. Dynkin diagrams & fixed parameters \\
\hline \hline
1 &
\rule[-4\unitlength]{0pt}{5\unitlength}
\begin{picture}(26,4)(0,3)
\put(13,3){\oval(26,6)}
\put(0,0){\makebox(26,6){\scriptsize $C(d,q;)$}}
\end{picture}
& $q\in k^\ast \setminus \{1\}$ \\
\hline
2 &
\begin{picture}(26,8)
\put(13,3){\oval(26,6)}
\put(0,0){\makebox(26,6){\scriptsize $C(d,q;i_1,\ldots ,i_j)$}}
\end{picture}
& $q\in k^\ast \setminus \{-1,1\}$ \\
 &
\begin{picture}(34,8)
\put(17,3){\oval(34,6)}
\put(0,0){\makebox(34,6){\scriptsize $C(d,q^{-1};i_1,\ldots ,i_{d+1-j})$}}
\end{picture}
& $1\le j\le (d+1)/2$ \\
\hline
3 &
\begin{picture}(48,8)
\put(17,3){\oval(34,6)}
\put(0,0){\makebox(34,6){\scriptsize $C(d-1,q^2;)$}}
\put(34,3){\line(1,0){10}}
\put(45,3){\circle{2}}
\put(40,4){\makebox[0pt]{\scriptsize $q^{-2}$}}
\put(45,5){\makebox[0pt]{\scriptsize $q$}}
\end{picture}
& $q\in k^\ast \setminus \{-1,1\}$ \\
\hline
4 &
\begin{picture}(44,8)
\put(15,3){\oval(30,6)}
\put(0,0){\makebox(30,6){\scriptsize $C(d-1,q^2;i_1,\ldots ,i_j)$}}
\put(30,3){\line(1,0){10}}
\put(41,3){\circle{2}}
\put(36,4){\makebox[0pt]{\scriptsize $q^{-2}$}}
\put(41,5){\makebox[0pt]{\scriptsize $q$}}
\end{picture}
& $q\in k^\ast $, $q^4\not=1$ \\
&
\begin{picture}(48,8)
\put(17,3){\oval(34,6)}
\put(0,0){\makebox(34,6){\scriptsize $C(d-1,q^{-2};i_1,\ldots ,i_{d-j})$}}
\put(34,3){\line(1,0){10}}
\put(45,3){\circle{2}}
\put(40,4){\makebox[0pt]{\scriptsize $q^2$}}
\put(45,5){\makebox[0pt]{\scriptsize $-q^{-1}$}}
\end{picture}
& $1\le j\le d-1$ \\
\hline
5 &
\begin{picture}(48,8)
\put(17,3){\oval(34,6)}
\put(0,0){\makebox(34,6){\scriptsize $C(d-1,-\zeta ^{-1};)$}}
\put(34,3){\line(1,0){10}}
\put(45,3){\circle{2}}
\put(40,4){\makebox[0pt]{\scriptsize $-\zeta $}}
\put(45,5){\makebox[0pt]{\scriptsize $\zeta $}}
\end{picture}
& $\zeta \in R_3$ \\
\hline
\end{tabular}


\begin{tabular}{r|l|l}
row & gener. Dynkin diagrams & fixed parameters \\
\hline \hline
10 & 
\rule[-4\unitlength]{0pt}{5\unitlength}
\begin{picture}(42,5)(0,3)
\put(15,3){\oval(30,6)}
\put(0,0){\makebox(30,6){\scriptsize $C(d-1,q;i_1,\ldots ,i_j)$}}
\put(30,3){\line(1,0){10}}
\put(41,3){\circle{2}}
\put(35,4){\makebox[0pt]{\scriptsize $q^{-2}$}}
\put(41,5){\makebox[0pt]{\scriptsize $q^2$}}
\end{picture}
& $q\in k^\ast \setminus \{-1,1\}$ \\
&
\rule[-10\unitlength]{0pt}{10\unitlength}
\begin{picture}(45,10)(0,10)
\put(17,10){\oval(34,6)}
\put(0,7){\makebox(34,6){\scriptsize $C(d-2,q;i_1,\ldots ,i_{j-1})$}}
\put(34,10){\line(1,1){7}}
\put(34,10){\line(1,-1){7}}
\put(41,18){\circle{2}}
\put(41,2){\circle{2}}
\put(41,3){\line(0,1){14}}
\put(39,15){\makebox[0pt][r]{\scriptsize $q^{-1}$}}
\put(37,5){\makebox[0pt][r]{\scriptsize $q^{-1}$}}
\put(43,17){\makebox[0pt][l]{\scriptsize $-1$}}
\put(42,10){\makebox[0pt][l]{\scriptsize $q^2$}}
\put(43,2){\makebox[0pt][l]{\scriptsize $-1$}}
\end{picture}
& $2\le j\le d-1$ \\
&
\rule[-10\unitlength]{0pt}{10\unitlength}
\begin{picture}(45,10)(0,10)
\put(17,10){\oval(34,6)}
\put(0,7){\makebox(34,6){\scriptsize $C(d-2,q^{-1};i_1,\ldots ,i_{d-j})$}}
\put(34,10){\line(1,1){7}}
\put(34,10){\line(1,-1){7}}
\put(41,18){\circle{2}}
\put(41,2){\circle{2}}
\put(38,15){\makebox[0pt][r]{\scriptsize $q$}}
\put(37,5){\makebox[0pt][r]{\scriptsize $q$}}
\put(43,17){\makebox[0pt][l]{\scriptsize $q^{-1}$}}
\put(43,2){\makebox[0pt][l]{\scriptsize $q^{-1}$}}
\end{picture}
& \\
\hline


%
%
%
%
11 &
\Dchainfive{$\zeta $}{$\zeta ^{-1}$}{$\zeta $}{$\zeta ^{-1}$}%
{$-1$}{$\zeta ^{-1}$}{$\zeta $}{$\zeta ^{-1}$}{$\zeta $}
& $\zeta \in R_3$ \\
 &
\rule[-4\unitlength]{0pt}{5\unitlength}
\begin{picture}(38,11)(0,3)
\put(1,1){\circle{2}}
\put(2,1){\line(1,0){10}}
\put(13,1){\circle{2}}
\put(13,2){\line(2,3){6}}
\put(14,1){\line(1,0){10}}
\put(25,1){\circle{2}}
\put(25,2){\line(-2,3){6}}
\put(19,12){\circle{2}}
\put(26,1){\line(1,0){10}}
\put(37,1){\circle{2}}
\put(1,4){\makebox[0pt]{\scriptsize $\zeta $}}
\put(7,3){\makebox[0pt]{\scriptsize $\zeta ^{-1}$}}
\put(12,4){\makebox[0pt]{\scriptsize $-1$}}
\put(19,3){\makebox[0pt]{\scriptsize $\zeta $}}
\put(26,4){\makebox[0pt]{\scriptsize $-1$}}
\put(31,3){\makebox[0pt]{\scriptsize $\zeta ^{-1}$}}
\put(37,4){\makebox[0pt]{\scriptsize $\zeta $}}
\put(15,8){\makebox[0pt]{\scriptsize $\zeta $}}
\put(23,8){\makebox[0pt]{\scriptsize $\zeta $}}
\put(22,11){\makebox[0pt]{\scriptsize $-1$}}
\end{picture}
& \\
\end{tabular}


\begin{tabular}{r|l|l}
row & gener. Dynkin diagrams & fixed parameters \\
\hline \hline
 &
\rule[-10\unitlength]{0pt}{12\unitlength}
\begin{picture}(36,10)(0,9)
\put(1,9){\circle{2}}
\put(2,9){\line(1,0){10}}
\put(13,9){\circle{2}}
\put(14,9){\line(1,0){10}}
\put(25,9){\circle{2}}
\put(26,9){\line(1,1){7}}
\put(26,9){\line(1,-1){7}}
\put(33,17){\circle{2}}
\put(33,1){\circle{2}}
\put(1,12){\makebox[0pt]{\scriptsize $\zeta $}}
\put(7,11){\makebox[0pt]{\scriptsize $\zeta ^{-1}$}}
\put(13,12){\makebox[0pt]{\scriptsize $-1$}}
\put(19,11){\makebox[0pt]{\scriptsize $\zeta $}}
\put(24,12){\makebox[0pt]{\scriptsize $\zeta ^{-1}$}}
\put(31,14){\makebox[0pt][r]{\scriptsize $\zeta ^{-1}$}}
\put(29,4){\makebox[0pt][r]{\scriptsize $\zeta $}}
\put(35,16){\makebox[0pt][l]{\scriptsize $\zeta $}}
\put(35,1){\makebox[0pt][l]{\scriptsize $-1$}}
\end{picture}
& \\
 &
\rule[-10\unitlength]{0pt}{12\unitlength}
\begin{picture}(36,10)(0,9)
\put(1,9){\circle{2}}
\put(2,9){\line(1,0){10}}
\put(13,9){\circle{2}}
\put(14,9){\line(1,0){10}}
\put(25,9){\circle{2}}
\put(26,9){\line(1,1){7}}
\put(26,9){\line(1,-1){7}}
\put(33,17){\circle{2}}
\put(33,1){\circle{2}}
\put(1,12){\makebox[0pt]{\scriptsize $-1$}}
\put(7,11){\makebox[0pt]{\scriptsize $\zeta $}}
\put(13,12){\makebox[0pt]{\scriptsize $-1$}}
\put(19,11){\makebox[0pt]{\scriptsize $\zeta ^{-1}$}}
\put(24,12){\makebox[0pt]{\scriptsize $-1$}}
\put(31,14){\makebox[0pt][r]{\scriptsize $\zeta ^{-1}$}}
\put(29,4){\makebox[0pt][r]{\scriptsize $\zeta $}}
\put(35,16){\makebox[0pt][l]{\scriptsize $\zeta $}}
\put(35,1){\makebox[0pt][l]{\scriptsize $-1$}}
\end{picture}
& \\
\hline
\end{tabular}


\begin{tabular}{r|l|l}
row & gener. Dynkin diagrams & fixed parameters \\
\hline \hline
13 &
\Dchainfive{$-1$}{$\zeta ^{-1}$}{$\zeta $}{$\zeta ^{-1}$}%
{$\zeta $}{$\zeta ^{-1}$}{$\zeta ^{-1}$}{$\zeta $}{$-1$}
& $\zeta \in R_3$ \\
 &
\Dchainfive{$-1$}{$\zeta $}{$-1$}{$\zeta ^{-1}$}%
{$\zeta $}{$\zeta ^{-1}$}{$\zeta ^{-1}$}{$\zeta $}{$-1$}
& \\
 &
\Dchainfive{$-1$}{$\zeta ^{-1}$}{$\zeta $}{$\zeta ^{-1}$}%
{$\zeta $}{$\zeta ^{-1}$}{$-1$}{$\zeta ^{-1}$}{$-1$}
& \\
 &
\Dchainfive{$\zeta $}{$\zeta ^{-1}$}{$-1$}{$\zeta $}%
{$-1$}{$\zeta ^{-1}$}{$\zeta ^{-1}$}{$\zeta $}{$-1$}
& \\
 &
\Dchainfive{$-1$}{$\zeta $}{$-1$}{$\zeta ^{-1}$}{$\zeta $}%
{$\zeta ^{-1}$}{$-1$}{$\zeta ^{-1}$}{$-1$}
& \\
 &
\rule[-10\unitlength]{0pt}{12\unitlength}
\begin{picture}(36,10)(0,9)
\put(1,9){\circle{2}}
\put(2,9){\line(1,0){10}}
\put(13,9){\circle{2}}
\put(14,9){\line(1,0){10}}
\put(25,9){\circle{2}}
\put(26,9){\line(1,1){7}}
\put(26,9){\line(1,-1){7}}
\put(33,17){\circle{2}}
\put(33,1){\circle{2}}
\put(33,2){\line(0,1){14}}
\put(1,12){\makebox[0pt]{\scriptsize $-1$}}
\put(7,11){\makebox[0pt]{\scriptsize $\zeta ^{-1}$}}
\put(13,12){\makebox[0pt]{\scriptsize $\zeta $}}
\put(19,11){\makebox[0pt]{\scriptsize $\zeta ^{-1}$}}
\put(24,12){\makebox[0pt]{\scriptsize $-1$}}
\put(31,14){\makebox[0pt][r]{\scriptsize $\zeta $}}
\put(29,4){\makebox[0pt][r]{\scriptsize $\zeta $}}
\put(35,16){\makebox[0pt][l]{\scriptsize $-1$}}
\put(34,9){\makebox[0pt][l]{\scriptsize $\zeta $}}
\put(35,1){\makebox[0pt][l]{\scriptsize $\zeta ^{-1}$}}
\end{picture}
& \\
 &
\Dchainfive{$\zeta $}{$\zeta ^{-1}$}{$\zeta $}{$\zeta ^{-1}$}{$-1$}%
{$\zeta $}{$-\zeta $}{$\zeta $}{$-1$}
& \\
\end{tabular}


\begin{tabular}{r|l|l}
row & gener. Dynkin diagrams & fixed parameters \\
\hline \hline
 &
\rule[-10\unitlength]{0pt}{12\unitlength}
\begin{picture}(36,10)(0,9)
\put(1,9){\circle{2}}
\put(2,9){\line(1,0){10}}
\put(13,9){\circle{2}}
\put(14,9){\line(1,0){10}}
\put(25,9){\circle{2}}
\put(26,9){\line(1,1){7}}
\put(26,9){\line(1,-1){7}}
\put(33,17){\circle{2}}
\put(33,1){\circle{2}}
\put(1,12){\makebox[0pt]{\scriptsize $\zeta ^{-1}$}}
\put(7,11){\makebox[0pt]{\scriptsize $\zeta $}}
\put(13,12){\makebox[0pt]{\scriptsize $-1$}}
\put(19,11){\makebox[0pt]{\scriptsize $\zeta ^{-1}$}}
\put(24,12){\makebox[0pt]{\scriptsize $\zeta $}}
\put(31,14){\makebox[0pt][r]{\scriptsize $\zeta ^{-1}$}}
\put(29,4){\makebox[0pt][r]{\scriptsize $\zeta ^{-1}$}}
\put(35,16){\makebox[0pt][l]{\scriptsize $\zeta $}}
\put(35,1){\makebox[0pt][l]{\scriptsize $-1$}}
\end{picture}
& \\
 &
\rule[-10\unitlength]{0pt}{12\unitlength}
\begin{picture}(36,10)(0,9)
\put(1,9){\circle{2}}
\put(2,9){\line(1,0){10}}
\put(13,9){\circle{2}}
\put(14,9){\line(1,0){10}}
\put(25,9){\circle{2}}
\put(26,9){\line(1,1){7}}
\put(26,9){\line(1,-1){7}}
\put(33,17){\circle{2}}
\put(33,1){\circle{2}}
\put(1,12){\makebox[0pt]{\scriptsize $-1$}}
\put(7,11){\makebox[0pt]{\scriptsize $\zeta ^{-1}$}}
\put(13,12){\makebox[0pt]{\scriptsize $-1$}}
\put(19,11){\makebox[0pt]{\scriptsize $\zeta $}}
\put(24,12){\makebox[0pt]{\scriptsize $\zeta ^{-1}$}}
\put(31,14){\makebox[0pt][r]{\scriptsize $\zeta ^{-1}$}}
\put(29,4){\makebox[0pt][r]{\scriptsize $\zeta $}}
\put(35,16){\makebox[0pt][l]{\scriptsize $\zeta $}}
\put(35,1){\makebox[0pt][l]{\scriptsize $-1$}}
\end{picture}
& \\
\hline
\end{tabular}


\begin{tabular}{r|l|l}
row & gener. Dynkin diagrams & fixed parameters \\
\hline \hline
14  &
\Dchainfive{$\zeta $}{$-\zeta $}{$\zeta $}{$-\zeta $}%
{$-1$}{$-1$}{$-1$}{$-\zeta $}{$\zeta $}
& $\zeta \in R_4$ \\
 &
\rule[-4\unitlength]{0pt}{5\unitlength}
\begin{picture}(38,11)(0,3)
\put(1,1){\circle{2}}
\put(2,1){\line(1,0){10}}
\put(13,1){\circle{2}}
\put(13,2){\line(2,3){6}}
\put(14,1){\line(1,0){10}}
\put(25,1){\circle{2}}
\put(25,2){\line(-2,3){6}}
\put(19,12){\circle{2}}
\put(26,1){\line(1,0){10}}
\put(37,1){\circle{2}}
\put(1,4){\makebox[0pt]{\scriptsize $\zeta $}}
\put(7,3){\makebox[0pt]{\scriptsize $-\zeta $}}
\put(12,4){\makebox[0pt]{\scriptsize $-1$}}
\put(19,3){\makebox[0pt]{\scriptsize $\zeta $}}
\put(26,4){\makebox[0pt]{\scriptsize $-1$}}
\put(31,3){\makebox[0pt]{\scriptsize $-\zeta $}}
\put(37,4){\makebox[0pt]{\scriptsize $\zeta $}}
\put(15,8){\makebox[0pt]{\scriptsize $\zeta $}}
\put(23,8){\makebox[0pt]{\scriptsize $-1$}}
\put(22,11){\makebox[0pt]{\scriptsize $-1$}}
\end{picture}
& \\
\hline
\end{tabular}


\begin{tabular}{r|l|l}
row & gener. Dynkin diagrams & fixed parameters \\
\hline \hline
16 &
\rule[-8\unitlength]{0pt}{9\unitlength}
\begin{picture}(41,8)(0,7)
\put(1,4){\circle{2}}
\put(2,4){\line(1,0){10}}
\put(13,4){\circle{2}}
\put(14,4){\line(1,0){10}}
\put(25,4){\circle{2}}
\put(26,4){\line(1,0){10}}
\put(37,4){\circle{2}}
\put(38,4){\line(1,0){10}}
\put(49,4){\circle{2}}
\put(25,5){\line(0,1){7}}
\put(25,13){\circle{2}}
\put(1,7){\makebox[0pt]{\scriptsize $q$}}
\put(7,6){\makebox[0pt]{\scriptsize $q^{-1}$}}
\put(13,7){\makebox[0pt]{\scriptsize $q$}}
\put(19,6){\makebox[0pt]{\scriptsize $q^{-1}$}}
\put(25,0){\makebox[0pt]{\scriptsize $q$}}
\put(31,1){\makebox[0pt]{\scriptsize $q^{-1}$}}
\put(37,7){\makebox[0pt]{\scriptsize $q$}}
\put(43,6){\makebox[0pt]{\scriptsize $q^{-1}$}}
\put(49,7){\makebox[0pt]{\scriptsize $q$}}
\put(26,8){\makebox[0pt][l]{\scriptsize $q^{-1}$}}
\put(27,12){\makebox[0pt][l]{\scriptsize $q$}}
\end{picture}
& $q\in k^\ast \setminus \{1\}$ \\
\hline
\end{tabular}


\begin{tabular}{r|l|l}
row & gener. Dynkin diagrams & fixed parameters \\
\hline \hline
18 &
\rule[-4\unitlength]{0pt}{5\unitlength}
\begin{picture}(62,4)(0,3)
\put(1,1){\circle{2}}
\put(2,1){\line(1,0){10}}
\put(13,1){\circle{2}}
\put(14,1){\line(1,0){10}}
\put(25,1){\circle{2}}
\put(26,1){\line(1,0){10}}
\put(37,1){\circle{2}}
\put(38,1){\line(1,0){10}}
\put(49,1){\circle{2}}
\put(50,1){\line(1,0){10}}
\put(61,1){\circle{2}}
\put(1,4){\makebox[0pt]{\scriptsize $\zeta $}}
\put(7,3){\makebox[0pt]{\scriptsize $\zeta ^{-1}$}}
\put(13,4){\makebox[0pt]{\scriptsize $\zeta $}}
\put(19,3){\makebox[0pt]{\scriptsize $\zeta ^{-1}$}}
\put(25,4){\makebox[0pt]{\scriptsize $\zeta $}}
\put(31,3){\makebox[0pt]{\scriptsize $\zeta ^{-1}$}}
\put(37,4){\makebox[0pt]{\scriptsize $\zeta $}}
\put(43,3){\makebox[0pt]{\scriptsize $\zeta ^{-1}$}}
\put(49,4){\makebox[0pt]{\scriptsize $\zeta ^{-1}$}}
\put(55,3){\makebox[0pt]{\scriptsize $\zeta $}}
\put(61,4){\makebox[0pt]{\scriptsize $-1$}}
\end{picture}
& $\zeta \in R_3$ \\
\end{tabular}


\begin{tabular}{r|l|l}
row & gener. Dynkin diagrams & fixed parameters \\
\hline \hline
 &
\rule[-4\unitlength]{0pt}{5\unitlength}
\begin{picture}(50,11)(0,3)
\put(1,1){\circle{2}}
\put(2,1){\line(1,0){10}}
\put(13,1){\circle{2}}
\put(14,1){\line(1,0){10}}
\put(25,1){\circle{2}}
\put(25,2){\line(2,3){6}}
\put(26,1){\line(1,0){10}}
\put(37,1){\circle{2}}
\put(37,2){\line(-2,3){6}}
\put(31,12){\circle{2}}
\put(38,1){\line(1,0){10}}
\put(49,1){\circle{2}}
\put(1,4){\makebox[0pt]{\scriptsize $\zeta $}}
\put(7,3){\makebox[0pt]{\scriptsize $\zeta ^{-1}$}}
\put(13,4){\makebox[0pt]{\scriptsize $-1$}}
\put(19,3){\makebox[0pt]{\scriptsize $\zeta $}}
\put(24,4){\makebox[0pt]{\scriptsize $\zeta ^{-1}$}}
\put(31,3){\makebox[0pt]{\scriptsize $\zeta $}}
\put(38,4){\makebox[0pt]{\scriptsize $-1$}}
\put(43,3){\makebox[0pt]{\scriptsize $\zeta ^{-1}$}}
\put(49,4){\makebox[0pt]{\scriptsize $\zeta $}}
\put(27,8){\makebox[0pt]{\scriptsize $\zeta $}}
\put(35,8){\makebox[0pt]{\scriptsize $\zeta $}}
\put(34,11){\makebox[0pt]{\scriptsize $-1$}}
\end{picture}
& \\
\hline
\end{tabular}


\begin{tabular}{r|l|l}
row & gener. Dynkin diagrams & fixed parameters \\
\hline \hline
19 &
\rule[-4\unitlength]{0pt}{5\unitlength}
\begin{picture}(62,4)(0,3)
\put(1,1){\circle{2}}
\put(2,1){\line(1,0){10}}
\put(13,1){\circle{2}}
\put(14,1){\line(1,0){10}}
\put(25,1){\circle{2}}
\put(26,1){\line(1,0){10}}
\put(37,1){\circle{2}}
\put(38,1){\line(1,0){10}}
\put(49,1){\circle{2}}
\put(50,1){\line(1,0){10}}
\put(61,1){\circle{2}}
\put(1,4){\makebox[0pt]{\scriptsize $\zeta $}}
\put(7,3){\makebox[0pt]{\scriptsize $-\zeta $}}
\put(13,4){\makebox[0pt]{\scriptsize $\zeta $}}
\put(19,3){\makebox[0pt]{\scriptsize $-\zeta $}}
\put(25,4){\makebox[0pt]{\scriptsize $\zeta $}}
\put(31,3){\makebox[0pt]{\scriptsize $-\zeta $}}
\put(37,4){\makebox[0pt]{\scriptsize $-1$}}
\put(43,3){\makebox[0pt]{\scriptsize $-1$}}
\put(49,4){\makebox[0pt]{\scriptsize $-1$}}
\put(55,3){\makebox[0pt]{\scriptsize $-\zeta $}}
\put(61,4){\makebox[0pt]{\scriptsize $\zeta $}}
\end{picture}
& $\zeta \in R_4$ \\
\hline
\end{tabular}


\begin{tabular}{r|l|l}
row & gener. Dynkin diagrams & fixed parameters \\
\hline \hline
20 &
\rule[-8\unitlength]{0pt}{9\unitlength}
\begin{picture}(62,8)(0,7)
\put(1,4){\circle{2}}
\put(2,4){\line(1,0){10}}
\put(13,4){\circle{2}}
\put(14,4){\line(1,0){10}}
\put(25,4){\circle{2}}
\put(26,4){\line(1,0){10}}
\put(37,4){\circle{2}}
\put(38,4){\line(1,0){10}}
\put(49,4){\circle{2}}
\put(50,4){\line(1,0){10}}
\put(61,4){\circle{2}}
\put(37,5){\line(0,1){7}}
\put(37,13){\circle{2}}
\put(1,7){\makebox[0pt]{\scriptsize $q$}}
\put(7,6){\makebox[0pt]{\scriptsize $q^{-1}$}}
\put(13,7){\makebox[0pt]{\scriptsize $q$}}
\put(19,6){\makebox[0pt]{\scriptsize $q^{-1}$}}
\put(25,7){\makebox[0pt]{\scriptsize $q$}}
\put(31,6){\makebox[0pt]{\scriptsize $q^{-1}$}}
\put(37,0){\makebox[0pt]{\scriptsize $q$}}
\put(43,1){\makebox[0pt]{\scriptsize $q^{-1}$}}
\put(49,7){\makebox[0pt]{\scriptsize $q$}}
\put(55,6){\makebox[0pt]{\scriptsize $q^{-1}$}}
\put(61,7){\makebox[0pt]{\scriptsize $q$}}
\put(38,8){\makebox[0pt][l]{\scriptsize $q^{-1}$}}
\put(39,12){\makebox[0pt][l]{\scriptsize $q$}}
\end{picture}
& $q\in k^\ast \setminus \{1\}$ \\
\hline
\end{tabular}

\end{appendix}


\textsc{Universit\"at Leipzig, Augustusplatz 10-11,
04109 Leipzig, Germany}

\textit{E-mail address:} \texttt{Istvan.Heckenberger@mathematik.uni-leipzig.de}

\end{document}